\documentclass[a4paper,twoside,10pt]{article}

\usepackage[a4paper,left=3cm,right=3cm, top=3cm, bottom=3cm]{geometry}
\usepackage[latin1]{inputenc}
\usepackage{mathrsfs}
\usepackage{graphicx}
\usepackage{epstopdf}
\graphicspath{{}}
\usepackage{amsmath}
\usepackage{amsthm}
\usepackage{amssymb}
\usepackage{stmaryrd}
\usepackage{float}
\usepackage{bigints}
\usepackage{cite}
\usepackage{color}
\usepackage[abs]{overpic}
\usepackage[font=footnotesize,labelfont=bf]{caption}
\usepackage{cases}
\usepackage{tikz}
\usepackage{algpseudocode}
\usepackage{rotating}
\usepackage{blkarray}
\usetikzlibrary{matrix,calc,arrows}
\usepackage[author={Lorenzo}]{pdfcomment}
\usepackage{soul,xcolor}
\usepackage{enumitem}  
\usepackage{comment}
\usepackage{verbatim}
\usepackage{graphicx}
\usepackage{subcaption}
\usepackage{mwe}
\usepackage{soul}
\usepackage{algorithm}
\usepackage{url}
\hypersetup{
    colorlinks=true,
    pdfborder={0 0 0}
}

\restylefloat{table}
\theoremstyle{plain}
\newtheorem{thm}{Theorem}[section]

\theoremstyle{definition}

\theoremstyle{remark}
\newtheorem{remark}{Remark}

\definecolor{darkgreen}{rgb}{0.0, 0.5, 0.0}
\newcommand{\blue}[1]{{\color{blue}{#1}}} %colors
\newcommand{\darkgreen}[1]{{\color{darkgreen}{#1}}} %colors
\newcommand{\red}[1]{{\color{red}{#1}}} %colors
\newcommand{\orange}[1]{{\color{orange}{#1}}} %colors

\newcommand{\eremk}{\hbox{}\hfill\rule{0.8ex}{0.8ex}}

\newcommand{\R}{\mathbb{R}}
\newcommand{\N}{\mathbb{N}}

\newcommand{\QT}{Q_T}
\newcommand{\SD}{\Sigma_D}
\newcommand{\SO}{\Sigma_0}
\newcommand{\Rbb}{\mathbb R}
\newcommand{\Nbb}{\mathbb N}
\newcommand{\Pbb}{\mathbb P}
\newcommand{\Pp}[2]{\Pbb_{#1}\left(#2\right)}
\newcommand{\dpt}[1]{\partial_t {#1}}

\newcommand{\x}{\mathbf x}
\newcommand{\nablax}{\nabla_{\x}}

\newcommand{\Deltax}{\Delta_{\x}}
\newcommand{\Norm}[1]{{\left\|{#1} \right\|}}
\newcommand{\SemiNorm}[1]{{\left|{#1} \right|}}
\newcommand{\dt}{\mbox{dt}}
\newcommand{\dx}{\mbox{d}\x}
\newcommand{\taun}{\mathcal T_h}

\newcommand{\E}{K}

\newcommand{\Ex}{\E_{\x}}
\newcommand{\tnmo}{t_{n-1}}

\newcommand{\h}{h}
\newcommand{\hE}{\h_\E}
\newcommand{\hEx}{\h_{\Ex}}

\newcommand{\F}{F}
\newcommand{\Fx}{\F_\x}
\newcommand{\Ftime}{\F_t}
\newcommand{\Ft}{\F_t}
\newcommand{\hF}{\h_\F}
\newcommand{\hFx}{\h_{\Fx}}
\newcommand{\Fcal}{\mathcal F}
\newcommand{\Fcalh}{\Fcal_\h}
\newcommand{\Fcalht}{\Fcalh^{\text{time}}}

\newcommand{\FcalEx}{\Fcal_\h^{\text{space}}}

\newcommand{\Vh}{V_\h}
\newcommand{\VhE}{\Vh(\E)}
\newcommand{\nbf}{\mathbf n}
\newcommand{\nbfF}{\nbf_\F}
\newcommand{\nbfFx}{\nbf_{\Fx}}
\newcommand{\malphaE}{m_{\alpha}^{\Kfrak}}

\newcommand{\malphaKfrakx}{m_{\gamma}^{\Kfrakx}}
\newcommand{\malphaF}{m_{\beta}^\F}
\newcommand{\calS}{\mathcal{S}}
\newcommand{\dS}{\mbox{d}S}

\newcommand{\uh}{u_h}

\newcommand{\vh}{v_h}

\newcommand{\Yh}{Y_\h}
\newcommand{\YE}{Y(\Kfrak)}
\newcommand{\Xtaun}{X(\taun)}
\newcommand{\Ytaun}{Y(\taun)}
\newcommand{\jump}[1]{\left[\!\left[#1\right]\!\right]}
\newcommand{\PiN}{\Pi^N_\p}

\newcommand{\Pistar}{\Pi^\star_\p}
\newcommand{\ah}{a_\h}

\newcommand{\bh}{b_\h}

\newcommand{\ahE}{\ah^{\E}}
\newcommand{\aE}{a^{\Kfrak}}
\renewcommand{\S}{S}
\newcommand{\SE}{\S^{\E}}

\newcommand{\Newton}{\mathfrak N}
\newcommand{\Newtonh}{\Newton_\h}

\newcommand{\gammaI}{\gamma_I}
\newcommand{\qp}{q_\p}

\newcommand{\qpmo}{q_{\p-1}}

\newcommand{\Pizpmo}{\Pi^{0}_{\p-1}}
\newcommand{\PizE}{\Pi^{0,\Kfrak}_{\p-1}}

\newcommand{\PizF}{\Pi^{0,\F}_{\p}}
\newcommand{\p}{p}

\newcommand{\phih}{\phi_\h}

\newcommand{\cH}{c_H}
\newcommand{\ctildeHE}{\widetilde c_H^{\Kfrak}}

\newcommand{\nutildeE}{\widetilde \nu^{\Kfrak}}

\newcommand{\etaE}{\eta_\E}
\newcommand{\NE}{N_{\taun}}
\newcommand{\pbf}{\mathbf \p}
\newcommand{\pbfFcalhEx}{\pbf^{\text{space}}}
\newcommand{\pbfFcalht}{\pbf^{\text{time}}}

\newcommand{\FcalExE}{\Fcal_{\Kfrak}^{\text{space}}}

\newcommand{\FcalEt}{\Fcal_{\Kfrak}^{\text{time}}}
\newcommand{\FcalEtj}{\Fcal_{\E_j}^{\text{time}}}
\newcommand{\Kfrak}{\mathcal K} 
\newcommand{\Kfrakx}{\Kfrak_{\x}}
\newcommand{\Kfrakt}{\Kfrak_t}

\newcommand{\at}{a_t}
\newcommand{\bt}{b_t}

\newcommand{\hfrakx}{\h_{\Kfrakx}}
\newcommand{\hfrakt}{\h_{\Kfrakt}}
\newcommand{\NC}{\mathcal{NC}}
\newcommand{\UW}{\mathcal {U}}
\newcommand{\UWEx}{\UW^{\Ex}}
\newcommand{\etaFEM}{\widetilde \eta}

\DeclareMathOperator{\card}{card}
\newcommand{\utildeh}{\widetilde u_\h}
\newcommand{\vtildeh}{\widetilde v_\h}
\newcommand{\Xtildeh}{\widetilde X_\h}
\newcommand{\Ytildeh}{\widetilde Y_\h}
\newcommand{\EcalY}{\mathcal E^Y}
\newcommand{\EcalU}{\mathcal E^U}
\newcommand{\EcalN}{\mathcal E^N}
\newcommand{\EcalX}{\mathcal E^X}
\newcommand{\EY}{E^Y}

\title{\Large{Design and performance of a space--time virtual element method for the heat equation on prismatic meshes}
\thanks{The authors have been funded by the Austrian Science Fund (FWF) through the projects F~65
and P~33477 (I. Perugia),
and the Italian Ministry of University and Research 
through the PRIN project ``NA-FROM-PDEs'' and PNRR-M4C2-I1.4-NC-HPC-Spoke6 (S. G\'omez).
S. G\'omez acknowledges the kind hospitality of the Erwin Schr\"odinger International Institute for Mathematics and Physics (ESI), where part of this research was developed.}}

\author{\large{Sergio G\'omez\thanks{Department of Mathematics, University of Pavia, 27100 Pavia, Italy (sergio.gomez01@universitadipavia.it)}
\thanks{Erwin Schr\"odinger Institute for Mathematics and Physics, University of Vienna, Austria},\; 
Lorenzo Mascotto\thanks{Department of Mathematics and Applications, University of Milano-Bicocca, 20125 Milan, Italy (lorenzo.mascotto@unimib.it)}
\thanks{Faculty of Mathematics, University of Vienna, 1090 Vienna, Austria (ilaria.perugia@univie.ac.at)}
\thanks{IMATI-CNR, Pavia, Italy},\;
Ilaria Perugia\footnotemark[5]}}
\date{}

\begin{document}
%%%%%%%%%%%%%%%%%%%%%%%%%%%%%%%%%%%%%
\maketitle

\begin{abstract}
\noindent
We present a space--time virtual element method for the discretization of
the heat equation,
which is defined on general prismatic meshes and variable degrees of accuracy.
Strategies to handle efficiently the space--time mesh structure
are discussed.
We perform convergence tests for the $\h$- and $\h\p$-versions 
of the method in case of smooth and singular solutions,
and test space--time adaptive mesh refinements
driven by a residual-type error indicator.


\end{abstract}

%%%%%%%%%%%%%%%%%%%%%%%%%%%%%%%%%%%%%%%%%%%%%%%%%%%%%%%%%%%%%%%%%%%%%%%%%%
\section{Introduction} \label{section:introduction}
%%%%%%%%%%%%%%%%%%%%%%%%%%%%%%%%%%%%%%%%%%%%%%%%%%%%%%%%%%%%%%%%%%%%%%%%%%
Space--time Galerkin methods aim at approximating solutions to time-dependent partial differential equations
treating the time variable as an additional space variable.
Even though the foundation of space--time Galerkin methods
traces back to the 70ies of the last century~\cite{Jamet:1978}
and a few contributions were developed
in the twenty years to follow
\cite{Babuska-Janik:1989, Babuska-Janik:1990, Ericksson-Johnson-Thomee-1985, Hughes-Stewart:1996, French-Peterson:1996},
only in the last two decades
there has been a growing attention on this topic,
mainly due to the improved performance of computers.

Compared to the time-stepping approach,
space--time Galerkin methods have some important
upsides and features:
the discrete solution can be evaluated on
the whole space--time domain and not only at a finite number of times
without additional post-processing of the discrete solution;
such methods allow for space--time adaptivity;
it is possible to design space--time parallel solvers.

In this paper, we focus on the approximation of solutions to the heat equation.
Several space--time methods have been designed to this aim
and they can be classified into two main groups. 
The first one is based on the discretization of a standard Petrov-Galerkin formulation~\cite{Dautray-Lions:1992};
see~\cite{Aziz_Monk:1989,Steinbach:2015} for continuous finite element methods,
\cite{Schwab-Stevenson:2009} for a wavelet method, \cite{Sudirham-VanderVegt:2006,Cangiani-Dong-Georgoulis:2017} for discontinuous Galerkin methods,
\cite{Langer-Moore-Neumuller:2016} for an isogeometric method,
\cite{Steinbach-Zank:2019} for a coercive method based on a Hilbert transformation of the test space,
and~\cite{Stevenson-Westerdiep:2021}
for a mixed finite element method. Residual-type error indicators for the method of~\cite{Steinbach:2015} were considered in~\cite{Steinbach-Yang:2018,Steinbach-Yang:2019}. 
Possible drawbacks of employing continuous finite elements are that
suboptimal convergence rates are obtained for some singular solutions,
and incompatible boundary and initial conditions
cannot be naturally handled.

The second group is based on first order system least squares discretizations (FOSLS);
see \cite{Fuhrer-Karkulik:2021, Gantner-Stevenson:2021, Voronin-Lee-Neumuller-Sepulveda-Vassilevski:2018} for finite element methods
and~\cite{Montardini-Negri-Sangalli-Tani:2020}
for isogeometric methods. Space--time FOSLS finite elements naturally provide reliable and efficient error indicators;
see~\cite{Fuhrer-Karkulik:2021, Schafelner-Vassilevski:2021, Gantner-Stevenson:2021}.
However, they require the computation of an additional vector-valued flux variable.

In this paper, we extend the nonconforming space--time virtual element method of~\cite{Gomez-Mascotto-Moiola-Perugia:2022}
to the case of general prismatic meshes and nonuniform degrees of accuracy.
This method allows for the use of space--time meshes with hanging facets (nodes, edges, faces)
and is based on discontinuous in time test and trial functions that are solutions to local space--time problems with polynomial data.
This approach provides a natural framework for space--time adaptivity without need of re-meshing neighbouring elements.
If the space--time mesh is decomposed into separate time-slabs,
the global linear system can be split into much smaller systems that can be solved sequentially.
Due to the nonconformity across time-like facets,
the design and implementation of the method
are independent of the spatial dimension.
No artificial compatibility of
initial and boundary conditions
is enforced on the trial virtual element space,
and optimal convergence rates have been achieved also for singular solutions.
However, the stability analysis of the method
relies on a discrete inf-sup condition 
and the method requires a stabilization term for the spatial Laplacian
that needs to be carefully designed;
see~\cite{Gomez-Mascotto-Moiola-Perugia:2022}.

The main goals of this manuscript are:
\begin{itemize}
    \item to construct a space--time virtual element method on general prismatic meshes, possibly with hanging facets and variable degrees of accuracy;
    \item to discuss suitable strategies to handle efficiently the space--time mesh structure;    
    \item to investigate~$\h$- and~$\h\p$- refinements numerically;
    \item to test space--time adaptive mesh refinements driven by a residual-type error indicator.
\end{itemize}

In the remainder of this section, we introduce some notation and the model problem.
Finally, we outline the structure of the paper.

\subsection{Notation} \label{subsection:intro-notation}
We denote the first partial derivative with respect to the time variable~$t$ by~$\dpt{}$,
and the spatial gradient and Laplacian operators by~$\nablax{}, \ \Deltax{} $, respectively.

Standard notation for Sobolev spaces is employed.
For a given domain~$D$ in~$\R^d$, $d \in \N$,
$H^s(D)$ represents the standard Sobolev space of order $s$ in~$\Nbb$
endowed with the standard inner product $(\cdot, \cdot)_{s,D}$,
the seminorm $\SemiNorm{\cdot}_{s,D}$,
and the norm~$\Norm{\cdot}_{s, D}$.
In particular, we let $H^0(D)$ be
the space~$L^2(D)$ of Lebesgue square integrable functions over~$D$
and $H_0^1(D)$ be the subspace of functions in~$H^1(D)$ with zero trace on~$\partial D$. 
If~$s$ is a fractional or negative number,
then the Sobolev space~$H^{s}(D)$ is defined by means of interpolation and duality, respectively.
We denote the duality product between $H^{-1}(D)$ and $H^1_0(D)$ by $\langle\cdot,\cdot\rangle$.
The Sobolev spaces on~$\partial D$ are defined analogously
and denoted by~$H^s(\partial D)$.

Given~$s$ in~$\Rbb$,
a time interval~$(a,b)$,
and a Banach space $(Z, \Norm{\cdot}_Z)$,
we introduce the Bochner space~$H^s(a,b; Z)$.
In particular, we define
\begin{equation}\label{X-Y}
Y := L^2(0, T; H_0^1(\Omega)),
\quad
X := \left\{v \in Y \cap H^1(  0, T; H^{-1}(\Omega)) 
        \mid v(\x, 0) = 0 \; \forall \x \in \Omega\right\}.
\end{equation}
We endow $Y$ and $X$ with the norms
\begin{equation} \label{Y-X-norms}
\Norm{v}_{Y}^2 := \nu \int_0^T \SemiNorm{v(\cdot, t)}_{1,\Omega}^2 \dt,
\qquad\qquad
\Norm{v}_{X}^2 := \cH \Norm{\dpt{v}}_{L^2(0,T;H^{-1}(\Omega))}^2
                  + \Norm{v}_{Y}^2,
\end{equation}
respectively, where
$\nu$ and $\cH$ are the problem coefficients defined at the beginning of Section~\ref{subsection:intro-model-problem} below,
and we have set 
%for all~$\phi$ in~$L^2(0,T;H^{-1}(\Omega))$,
\[
\Norm{\phi}_{L^2(0,T;H^{-1}(\Omega))}
:= \sup_{0\ne v \in Y} \frac{\int_0^T \langle \phi,v \rangle \dt}{\Norm{v}_{Y}}.
\]
We denote the space of polynomials in $d$ variables of degree at most~$\p$
on a domain~$D\subset\Rbb^d$ by~$\Pbb_\p(D)$.

\subsection{Model problem} \label{subsection:intro-model-problem}
We consider the heat equation 
on the space--time domain~$\QT := \Omega \times (0, T)$,
where~$\Omega\subset\Rbb^d$, $d = 1,\ 2,\ 3$,
and~$T > 0$
are the (bounded) spatial domain and the final time, respectively.

Associated with~$\QT$,
we introduce the space-like surface~$\SO := \Omega \times \left\{0\right\}$
and the time-like surface~$\SD := \partial \Omega \times (0, T)$.

Let~$f:\QT\to\Rbb$, $u_0:\SO\to\Rbb$, and~$g:\SD \to \Rbb$
denote the source term, initial condition, and Dirichlet boundary condition, respectively.
We assume $f\in L^2(0,T; H^{-1}(\Omega))$.
We further consider~$\cH>0$ and~$\nu>0$
given positive constant volumetric heat capacity and thermal conductivity, respectively.

The heat equation in strong formulation reads:
find a function~$u:\QT \to \Rbb$ such that
\begin{equation} \label{continuous-strong}
\begin{cases}
\cH \dpt{u} - \nu \Deltax u = f & \text{in } \QT, \\
u = u_0                         & \text{on } \SO, \\
u = g                           & \text{on } \SD.
\end{cases}
\end{equation}
To simplify the presentation, we henceforth assume that~$u_0=0$ and~$g=0$.
Nonetheless,
in Sections~\ref{section:numerical-a-priori} and~\ref{section:residual-type},
we present numerical results for solutions with inhomogeneous initial and boundary conditions.
Inhomogeneous boundary conditions can be treated as in the elliptic case
and inhomogeneous initial conditions can be dealt with
by standard lifting arguments.

Define the space--time bilinear form
$b(\cdot,\cdot):X\times Y \to \Rbb$ as
\begin{equation} \label{continuous:bf:global}
b(u, v) := \int_0^T \langle\cH\dpt u,v\rangle \dt
                                + \int_0^T\int_\Omega \nu \nablax u \cdot \nablax v\,\dx \dt.
\end{equation}
We consider the Petrov-Galerkin weak formulation of~\eqref{continuous-strong} 
as in~\cite{Dautray-Lions:1992}:
\begin{equation} \label{continuous-weak}
\begin{cases}
\text{find } u \in X \text{ such that}\\
b(u,v)= \int_0^T \langle f,v \rangle \dt
\qquad \forall v \in Y.
\end{cases}
\end{equation}
Problem~\eqref{continuous-weak} is well posed;
see, e.g., \cite[Corollary~$2.3$]{Steinbach:2015}.

\paragraph*{Structure of the paper.}
In Section~\ref{section:VEM}, we design a space--time virtual element method
for the approximation of solutions to~\eqref{continuous-weak}
on general prismatic meshes;
in Section~\ref{section:handling-meshes-degrees},
strategies to handle efficiently the space--time mesh structure
are discussed,
including flagging strategies that might be used for other polytopic methods as well;
in Section~\ref{section:numerical-a-priori},
we assess the convergence of the $\h$- and $\h\p$-versions of the method
on some singular solutions;
in Section~\ref{section:residual-type},
a computable residual-type error indicator
is defined and used
to lead adaptive mesh refinements.

%%%%%%%%%%%%%%%%%%%%%%%%%%%%%%%%%%%%%%%%%%%%%%%%%%%%%%%%%%%%%%%%%%%%%%%%%%
\section{The space--time virtual element method} \label{section:VEM}
%%%%%%%%%%%%%%%%%%%%%%%%%%%%%%%%%%%%%%%%%%%%%%%%%%%%%%%%%%%%%%%%%%%%%%%%%%
We describe a nonconforming space--time virtual element method for the approximation
of solutions to~\eqref{continuous-weak}.
We proceed in several steps:
we introduce general prismatic meshes in Section~\ref{subsection:intro-meshes};
design local virtual element spaces
and describe their degrees of freedom (DoFs)
in Section~\ref{subsection:local-VE-spaces};
show that the choice of the DoFs allows for the computation of orthogonal projectors onto polynomial spaces
in Section~\ref{subsection:polynomial-projections};
define global nonconforming space--time virtual element spaces in Section~\ref{subsection:global-VE-spaces};
detail the discrete bilinear forms in Section~\ref{subsection:discrete-bilinear-forms};
present the method in Section~\ref{subsection:method}.

\subsection{General prismatic meshes} \label{subsection:intro-meshes}

We consider sequences of meshes~$\{ \taun \}$
consisting of nonoverlapping prismatic-type polytopes covering~$\QT$ in the following sense:
each element~$\Kfrak$ in~$\taun$ can be written as~$\Kfrakx \times \Kfrakt$
for some open~$d$-dimensional polytope~$\Kfrakx$ contained in~$\Omega$ with boundary~$\partial \Kfrakx$,
and some time subinterval~$\Kfrakt := (a_t, b_t)$ of~$(0, T)$.
We denote the diameter of~$\Kfrakx$ and the length of~$\Kfrakt$
by~$\hfrakx$ and~$\hfrakt$, respectively.

We call ``mesh facet" any intersection~$\partial \Kfrak^1 \cap \partial \Kfrak^2$,
$\partial \Kfrak^1  \cap (\partial \Omega \times (0, T))$,
$\partial \Kfrak^1 \cap (\Omega \times \{0\})$,
or~$\partial \Kfrak^1 \cap (\Omega \times \{T\})$,
for given~$\Kfrak^1, \Kfrak^2$ in~$\taun$,
that has positive~$d$-dimensional measure
and is contained in a~$d$-dimensional hyperplane.
For each element~$\Kfrak$ in~$\taun$,
we identify two types of nonoverlapping mesh facets:
space-like facets~$\Ex$,
whose union gives~$\Kfrakx \times \{\at\}$, and
time-like facets~$\F := \Fx \times \Ftime$,
where~$\Fx \subset \partial \Kfrakx$ is a facet of~$\Kfrakx$ and~$\Ftime \subset \Kfrakt$.
\footnote{For our purposes,
we do not need to classify the facets of~$\Kfrak$ on~$\Kfrakx \times \{\bt\}$.
One may consider that the facet~$\Kfrakx \times \{\bt\}$ is a unique facet of $\Kfrak$.}
We collect the space-like and time-like facets
of~$\Kfrak$ into the sets~$\FcalExE$ and~$\FcalEt$. 

Associated with each element~$\Kfrak$ in~$\taun$, we define~$\E$ as the~$(d+1)$-dimensional closed polytope~$\E$ whose interior is~$\Kfrak$
and whose boundary facets are those
in the set~$\FcalExE \cup \FcalEt \cup \{\Kfrakx \times \{b_t\}\}$.

\begin{figure}[!ht]
\centering
\begin{minipage}{0.32\textwidth}
\begin{center}
\begin{tikzpicture}[scale=0.04]
\draw[black, thick, -] (0,0) -- (100,0) -- (100,100) -- (0,100) -- (0,0);
\draw[black, thick, -, dashed] (0, 50) -- (100, 50);
\draw[black, thick, -, dashed] (50, 0) -- (50, 100);
\draw[black, thick, -, dashed] (25, 0) -- (25, 50);
\draw[black, thick, -, dashed] (50, 75) -- (100, 75);
\draw[black, thick, -, dashed] (75, 75) -- (75, 100);
\draw[fill=yellow, opacity=0.15] (0, 50) -- (50,50) -- (50, 100) -- (0, 100) -- (0, 50);
\draw (25, 75) node[black] {$\Kfrak$};
\draw[fill=magenta, opacity=0.15] (50, 50) -- (100,50) -- (100, 75) -- (50, 75) -- (50, 50);
\draw[fill=magenta, opacity=0.15] (50, 75) -- (75, 75) -- (75, 100) -- (50, 100) -- (50, 75);
\draw (75, 62.5) node[black] {$\widetilde{\Kfrak}^1$}; 
\draw (62.5, 87.5) node[black] {$\widetilde{\Kfrak}^2$}; 
%\draw[red, thick, -] (0.5, 50) -- (24.5, 50);
%\draw[red, thick, -] (25.5, 50) -- (49.5, 50);
%\draw (12.5, 56) node[red] {\tiny{$\Ex^1$}};
%\draw (37.5, 56) node[red] {\tiny{$\Ex^2$}};
%\draw[blue, thick, -] (50, 50.5) -- (50, 74.5);
%\draw[blue, thick, -] (50, 75.5) -- (50, 99.5);
%\draw[blue, thick, -] (0, 50.5) -- (0, 99.5);
%\draw (46, 87.5) node[blue] {\tiny{$F^2$}};
%\draw (46, 62.5) node[blue] {\tiny{$F^1$}};
%\draw (5, 75) node[blue] {\tiny{$F^3$}};
%\draw[black, thick, -] (75, 75) -- (75, 100);
%%
\end{tikzpicture}
\end{center}
\end{minipage}
\begin{minipage}{0.32\textwidth}
\begin{center}
\begin{tikzpicture}[scale=0.04]
\draw[black, thick, dashed, -] (100,0) -- (100,100) -- (0,100);
\draw[fill=yellow, opacity=0.15] (0.5, 0.5) -- (99.5, 0.5) -- (99.5, 99.5) -- (0.5, 99.5) -- (0.5, 0.5);
\draw (50, 50) node[black] {$\Kfrak$};
\draw[blue, thick, dashed, -] (0.5, 0) -- (99.5, 0);
\draw (50, 10) node[blue] {{$\Kfrakx$}};
\draw[red, thick, dashed, -] (0, 0.5) -- (0, 99.5);
\draw (10, 50) node[red] {{$\Kfrakt$}};
\end{tikzpicture}
\end{center}
\end{minipage}
\begin{minipage}{0.32\textwidth}
\begin{center}
\begin{tikzpicture}[scale=0.04]
\draw[black, thick, -] (0,0) -- (100,0) -- (100,100) -- (0,100) -- (0,0);
\draw[fill=yellow, opacity=0.15] (0, 0) -- (100, 0) -- (100, 100) -- (0, 100) -- (0, 0);
\draw (50, 50) node[black] {$K$};
\draw[blue, thick, -] (0.5, 0) -- (49.5, 0);
\draw[blue, thick, -] (50.5, 0) -- (99.5, 0);
%\draw[blue, thick, -] (0.5, 100) -- (99.5, 100);
\draw (25, 10) node[blue] {{$\Ex^1$}};
\draw (75, 10) node[blue] {{$\Ex^2$}};
%\draw (50, 90) node[blue] {{$\Ex^3$}};
\draw[red, thick, -] (100, 0.5) -- (100, 49.5);
\draw[red, thick, -] (100, 50.5) -- (100, 99.5);
\draw[red, thick, -] (0, 0.5) -- (0, 99.5);
\draw (90, 75) node[red] {{$F^2$}};
\draw (90, 25) node[red] {{$F^1$}};
\draw (10, 50) node[red] {{$F^3$}};
\foreach \point in{ (0, 0), (50, 0), (100, 0), (100, 50), (100, 100), (0, 100)}{
\draw[fill = black] \point circle (30pt);
}
\end{tikzpicture}
\end{center}
\end{minipage}
\caption{Example of a prismatic space--time mesh. 
\textbf{Left panel:} The prismatic partition~$\taun$ of the space--time domain.
\textbf{Central panel:} Zoom of the element~$\Kfrak \in \taun$.
\textbf{Right panel:} The associated closed hexagon~$\E$ with two space-like facets~$\Ex^1$ and~$\Ex^2$,
and three time-like facets~$F^1, F^2$, and~$F^3$.
By definition, $h_{F_\x^1} = \min\{\hfrakx, h_{\widetilde{\Kfrak}^1_\x}\}$,
$h_{F_\x^2} = \min\{\hfrakx, h_{\widetilde{\Kfrak}^2_\x}\}$,
and~$h_{F_\x^3} = \hfrakx$.
\label{FIG::EXAMPLE-PRISMATIC-ELEMENTS}}
\end{figure}
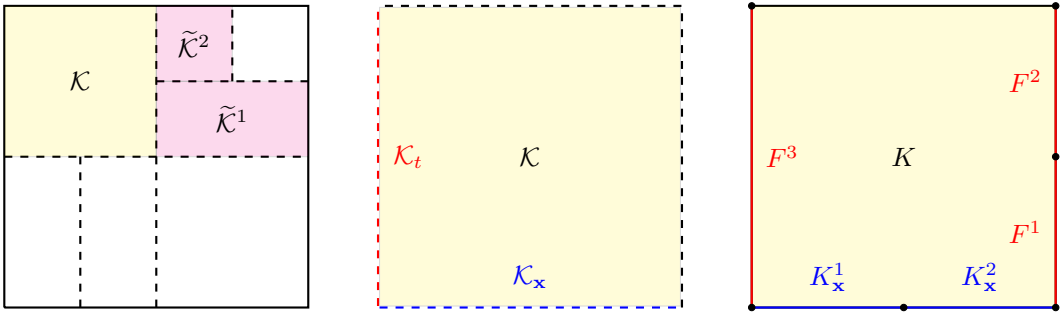

For any time-like facet~$\F = \Fx \times \Ft$ in~$\FcalEt$, we define
\begin{equation*}
\hFx:=\begin{cases}
\min\{\hfrakx,h_{\widetilde{\Kfrak}_\x}\}
&\text{if } \F = \Kfrak \cap \widetilde{\Kfrak} 
\text{ for some } 
\widetilde{\Kfrak} = \widetilde{\Kfrak}_{\x} \times \widetilde{\Kfrak}_t\in\taun,\\
\hEx & \text{if } \Fx\subset \partial\Omega.
\end{cases}
\end{equation*}
We define the sets of \emph{space-like}
and \emph{time-like} facets of~$\taun$ as
\[
\begin{split}
\FcalEx  := \bigcup\limits_{\E \in \taun} \FcalExE,
\qquad\qquad
\Fcalht  := \bigcup\limits_{\E \in \taun} \FcalEt.
\end{split}
\]
Hanging nodes ($1+1$ dimensional case),
edges ($2+1$ dimensional case),
and faces ($3+1$ dimensional case)
are included within this structure of the mesh.
In the absence of hanging facets,
i.e., for tensor-product-in-time meshes,
it is not necessary to distinguish between~$\Kfrak$ and~$\E$,
and all the definitions in this work reduce to those in~\cite{Gomez-Mascotto-Moiola-Perugia:2022}.
For each~$\F$ in~$\Fcalht$,
we fix a normal unit vector~$\nbfF$ in~$\Rbb^{d+1}$ of the form~$(\nbfFx,0)$, $\nbfFx\in\Rbb^d$.
In Figure~\ref{FIG::EXAMPLE-PRISMATIC-ELEMENTS}, we illustrate the definitions in this section.

For a given mesh~$\taun$, we define the broken Sobolev space of order~$s \in \Rbb^+$ as
\[
H^s(\taun)
:= \left\{
v \in L^2(\Omega)
\ \middle| \ 
v{}_{|\E} \in H^s(\E) \; \forall \E \in \taun
\right\} ,
\]
and endow it with the standard broken norm~$\Norm{\cdot}_{s,\taun}$
and seminorm~$\SemiNorm{\cdot}_{s,\taun}$.
We denote the space of piecewise polynomials of maximum degree~$\p$ in~$\N$ over~$\taun$ by~$\calS_\p(\taun)$.

%%%%%%%%%%%%%
\subsection{Local space--time virtual element spaces} \label{subsection:local-VE-spaces}
%%%%%%%%%%%%%
We introduce local space--time virtual element spaces,
extending the construction in~\cite{Gomez-Mascotto-Moiola-Perugia:2022}
to more general meshes as in Section~\ref{subsection:intro-meshes}.

Introduce the scaling factors
\begin{equation} \label{ctilde-nutilde}
\ctildeHE := \hfrakt ,
\qquad \qquad
\nutildeE := \hfrakx^2.
\end{equation}
Given~$\p$ in~$\Nbb$ and~$\Kfrak$ an element in~$\taun$,
we define the local space--time virtual element space as
\begin{equation} \label{local-VE-space}
\begin{split}
\VhE := \Big\{v \in L^2(\Kfrak) \mid\ 
& \ctildeHE \dpt \vh - \nutildeE \Deltax \vh  \in \Pp{p-1}{\Kfrak};\\[0.2em]
& \vh{}_{|\Kfrakx} \in \Pp{p}{\Kfrakx};\;
  \nbfFx \cdot \nablax \vh{}_{|\F} \in \Pp{\p}{\F} \, \forall \F \in \FcalEt \Big\}.
\end{split}
\end{equation}
Functions in~$\VhE$ are not known in closed form.
Yet, $\Pbb_{\p}(\Kfrak)$ is contained in~$\VhE$.

Given an element~$\Kfrak$ of~$\taun$
and any of its time-like facets~$\F$,
let~$\{ \malphaE \}_{\alpha=1}^{\dim(\Pbb_{\p-1}(\Kfrak))}$,
$\{ \malphaF{} \}_{\beta=1}^{\dim(\Pbb_{\p}(\F))}$,
and~$\{ \malphaKfrakx \}_{\gamma=1}^{\dim(\Pbb_{\p}(\Kfrakx))}$
be given bases of~$\Pbb_{p-1}(\Kfrak)$, $\Pbb_{\p}(\F)$, and~$\Pbb_{\p}(\Kfrakx)$, respectively.
Those basis elements are assumed to be invariant with respect to translations and dilations.

Introduce the following set of linear functionals:
\begin{itemize}
\item the \emph{bulk} moments
\begin{equation} \label{bulk-dofs}
\frac{1}{\vert \Kfrak \vert}
\int_{\E} \vh \ \malphaE \dx\ \dt,
\qquad  \alpha=1,\dots, \dim(\Pbb_{p-1}(\Kfrak));
\end{equation}
\item for all \emph{time-like} facets~$\F \in \FcalEt$, the \emph{time-like} moments
\begin{equation} \label{vertical-dofs}
\frac{1}{\vert\F\vert}
\int_{\F} \vh \ \malphaF \dS\ \dt,
\qquad  \beta=1,\dots,\dim(\Pbb_{\p}(\F));
\end{equation}
\item the \emph{space-like} moments
\begin{equation} \label{bottom-dofs} 
\frac{1}{\vert \Kfrakx \vert} 
\int_{\Kfrakx} \vh(\cdot, \tnmo) \malphaKfrakx \dx,
\qquad  \gamma=1,\dots,\dim(\Pbb_{p}(\Kfrakx)).
\end{equation}
\end{itemize}
These linear functionals constitute a set of unisolvent degrees of freedom
for the space~$\VhE$.
To see this, it suffices to extend \cite[Lemma~$2.1$]{Gomez-Mascotto-Moiola-Perugia:2022}
to the case of the general prismatic-type elements introduced in Section~\ref{subsection:intro-meshes}.

%%%%%%%%%%%%%
\subsection{Polynomial projections} \label{subsection:polynomial-projections}
%%%%%%%%%%%%%
We define four orthogonal projectors onto polynomial spaces as in~\cite{Gomez-Mascotto-Moiola-Perugia:2022}:
for all elements~$\Kfrak$ in~$\taun$ and all~$\varepsilon>0$,
\begin{enumerate}
\item we introduce
$ \PiN:   H^{\frac12+\varepsilon}(\Kfrakt; L^2(\Kfrakx))
    \cap L^2(\Kfrakt; H^1(\Kfrakx)) \to \Pbb_p(\Kfrak)$ as
\begin{subequations} \label{PiN}
\begin{align}
& \int_{\Kfrakt}\int_{\Kfrakx} \nablax \qp^{\Kfrak} \cdot \nablax \left(\PiN{v} - v \right) \ \dx \dt = 0 \quad \forall \qp^{\Kfrak} \in \Pp{p}{\Kfrak} \backslash \Pp{p}{\Kfrakt},    \label{PiN-1} \\
& \int_{\Kfrakt} \int_{\Kfrakx} \qpmo(t) \left(\PiN{v} - v \right) \dx \dt = 0 \quad \forall \qpmo \in \Pp{p-1}{\Kfrakt}, \label{PiN-2} \\
& \int_{\Kfrakx} \left(\PiN{v}(\x,\tnmo) - v(\x, \tnmo) \right) \dx = 0;    \label{PiN-3}
\end{align}
\end{subequations}
\item we introduce
$\Pistar: H^{\frac12+\varepsilon}(\Kfrakt; L^2(\Kfrakx)) \to \Pbb_p(\Kfrak)$ as
\begin{subequations} \label{Pistar}
\begin{align}
& \int_{\Kfrakt} \int_{\Kfrakx} \qpmo^{\Kfrak} \left(\Pistar{v} - v \right) \dx \dt 
    = 0\quad \forall \qpmo^{\Kfrak} \in \Pp{p-1}{\Kfrak}, \label{Pistar-1} \\
&\int_{\Kfrakx} \qp^{\Kfrakx} \left(\Pistar{v}(\x, \tnmo)  - v(\x, \tnmo) \right) 
    \dx  = 0 \quad \forall \qp^{\Kfrakx} \in \Pp{p}{\Kfrakx} \label{Pistar-2};
\end{align}
\end{subequations}
\item we introduce~$\PizE: L^2(\Kfrak) \to \Pbb_{\p-1}(\Kfrak)$ as
\begin{equation} \label{L2-projection:bulk}
(\qpmo^{\Kfrak}, v-\PizE v)_{0, \Kfrak} = 0
\qquad
\forall \qpmo^{\Kfrak} \in \Pbb_{\p-1}(\Kfrak);
\end{equation}
\item for all \emph{time-like} facets~$\F$ in~$\FcalEt$,
we introduce~$\PizF: L^2(\F) \to \Pbb_{\p}(\F)$ as
\begin{equation} \label{L2-projection:vertical-face}
(\qp^{\F}, v-\PizF v)_{0,\F} = 0 
\qquad
\forall \qp^{\F} \in \Pbb_{\p}(\F).
\end{equation}
\end{enumerate}
Let~$v$ be a function in~$\VhE$ with given degrees of freedom \eqref{bulk-dofs}--\eqref{bottom-dofs}. 
Then, the proof of the well posedness and computability of the above four projectors
follows by extending the results in~\cite[Section~$2.2$]{Gomez-Mascotto-Moiola-Perugia:2022}
to the case of the general prismatic-type elements in Section~\ref{subsection:intro-meshes}.

%%%%%%%%%%%%%
\subsection{Global nonconforming space--time virtual element spaces} \label{subsection:global-VE-spaces}
%%%%%%%%%%%%%
We design a global virtual element space~$\Yh$,
consisting of functions that are discontinuous in time and
nonconforming in space.

To this aim, we introduce the jump operator
on each time-like facet~$\F$
as the functional~$\jump{\cdot}_\F: H^{\frac12+\varepsilon}(\taun) \to  [L^2(\F)]^{d+1}$,
$\varepsilon >0$,
given by
\[
\jump{v}_{\F}:=
\begin{cases}
v_{|\Kfrak_1} \nbf_{\E_1}^\F + v_{|\Kfrak_2} \nbf_{\E_2}^\F 
& \text{if } \F \subset \partial \Kfrak_1 \cap \partial \Kfrak_2 \text{ is an internal facet, for } \Kfrak_1, \Kfrak_2 \in \taun\\
 v_{|\Kfrak_3}  \nbf_{\E_3}^\F
 & \text{if } \F \subset \partial \Kfrak_3 \text{ is a  boundary facet, for } \Kfrak_3\in\taun.
\end{cases}     
\]
%Recall that any element~$\E$ shares its volume with a prism~$\Kfrak = \Kfrakx \times \Kfrakt$;
%see Section~\ref{subsection:intro-meshes}.
We define the Sobolev nonconforming space of order~$\p$
associated with the mesh~$\taun$ as
\begin{equation} \label{nonconforming-condition}
\begin{split}
\NC_\p(\taun)
:=
\Big\{
v\in L^2 (\QT) \Big| 
& \,v_{|\Kfrak} \in L^2(\Kfrakt, H^1(\Kfrakx))\;
    \forall \Kfrak \in \taun;    \\
& \int_{\F} \qp^\F \jump{v}_{\F} \cdot \nbfF \dS\ = 0 \quad \forall \qp^\F \in \Pp{\p}{\F},\;
\forall \F \in \Fcalht \Big\}.
\end{split}
\end{equation}
We define the global virtual element space as follows:
\[
\begin{split}
\Yh 
:= \Big\{ \vh \in L^2(\QT) \ \Big| \
& \vh{}_{|\Kfrak} \in \Vh(\E) \;\; \forall\E \in \taun;\ 
    \vh \in \NC_\p(\taun) \Big\}.
\end{split}
\]

%%%%%%%%%%%%%
\subsection{The discrete bilinear forms} \label{subsection:discrete-bilinear-forms}
%%%%%%%%%%%%%
Since functions in the space--time virtual element space~$\Yh$ are not available in closed form,
we discretize the bilinear forms by computable counterparts.
On each element~$\Kfrak$, define the local continuous bilinear form
and the induced norm
\[
\aE(\uh,\vh) := \nu (\nablax \uh, \nablax \vh)_{0, \Kfrak},
\qquad\qquad
\SemiNorm{\vh}_{\YE}^2 := \aE(\vh, \vh)
= \nu \Norm{\nablax \vh}_{0, \Kfrak}^2.
\]
Let~$\SE:\VhE\times\VhE \to \Rbb$ be  any symmetric bilinear form,
which is computable via the degrees of freedom and satisfies the following property:
there exist~$0 < c_* < c^*$ independent of~$\E$ such that
\begin{equation} \label{stab-prop-0}
    c_* \SemiNorm{\vh}_{\YE}^2
    \le \nu \SE(\vh,\vh)
    \le c^* \SemiNorm{\vh}_{\YE}^2
    \qquad\qquad \forall \vh \in \VhE \cap \ker(\PiN).
\end{equation}
Define
\begin{equation} \label{ahtildeE}
\ahE(\uh,\vh)
:=  \aE(\PiN\uh, \PiN \vh) 
    + \nu \SE( (I-\PiN)\uh, (I-\PiN)\vh ).
\end{equation}
Extending \cite[Lemma~2.7]{Gomez-Mascotto-Moiola-Perugia:2022}
to the case of general prismatic-type elements
as in Section~\ref{subsection:intro-meshes},
it is possible to prove that there exist~$0<\alpha_*<\alpha^*$ independent of~$\E$ such that
the following local stability bounds are valid:
\begin{equation} \label{stability-bounds:local}
\alpha_* \SemiNorm{\vh}_{\YE}^2
\le \ahE(\vh,\vh)
\le \alpha^* \SemiNorm{\vh}_{\YE}^2
\qquad\qquad 
\forall \vh\in \VhE.
\end{equation}
Besides, we introduce upwind-type terms,
which allows for imposing weakly the continuity in time of the trial functions:
for all space-like facets~$\Ex \subset \Omega \times \{ t^*\}$, $t^* \in [0,T]$,
\begin{equation} \label{upwind-term}
\UWEx(\uh)
:= 
\begin{cases}
\cH \Pistar \uh{}_{|\Kfrak}(\cdot,t^*)
& \text{if }  t^*=0 \\
\cH \left(\Pistar \uh{}_{|\Kfrak^+}(\cdot, t^*) - \Pistar \uh{}_{|\Kfrak^-}(\cdot,t^*)\right)
& \text{if } t^*>0,
\end{cases}
\end{equation}
where~$\Ex \in \FcalExE$ if~$t^*=0$,
and~$\Ex \in F^{\text{space}}_{\Kfrak^+}$ is a subset of~$\partial \Kfrak^+\cap \partial \Kfrak^-$ if~$t^*>0$. 

\begin{remark}{\label{rem:U0}}
The case of nonzero initial conditions can be dealt with
by modifying the definition of the upwind functional at time~$t = 0$ as follows: for all space-like facets~$\Ex \subset \Omega \times \{ 0\}$,
\begin{equation*}
\UWEx(\uh)
:= 
\cH \left(\Pistar \uh{}_{|\Kfrak}(\cdot, 0) - u_0\right).
\end{equation*}\eremk
\end{remark}

The computability of the upwind terms follows from the definition
of the space-like moment degrees of freedom~\eqref{bottom-dofs}
and the computability of the projector~$\Pistar$.

We define a discrete counterpart of the bilinear
form~$b(\cdot,\cdot)$ in~\eqref{continuous:bf:global} as follows:
\begin{equation} \label{bh}
\bh(\uh,\vh)
:= \sum_{\Kfrak \in \taun} \left[ \cH(\dpt \Pistar \uh,\vh)_{0, \Kfrak}
                       + \ahE(\uh,\vh) \right]
    + \sum_{\Ex \in \FcalEx} \left(\UWEx(\uh),\vh{}_{|\Kfrak^+}\right)_{0, \Ex}.
\end{equation}
A computable stabilization satisfying~\eqref{stability-bounds:local} is given by

\begin{equation} 
\label{explicit:stabilization}
\begin{split}
\SE(\uh,\vh)
& := p^2\hfrakx^{-2} (\PizE \uh, \PizE \vh)_{0, \Kfrak}
   + \sum_{\F \in \FcalEt}
 p \hFx^{-1} (\PizF \uh, \PizF \vh)_{0,\F}\\
&\quad + p\hfrakt \hfrakx^{-2} (\uh, \vh)_{0,\Kfrakx}.
\end{split}
\end{equation}

%%%%%%%%%%%%%
\subsection{The method} \label{subsection:method}
%%%%%%%%%%%%%
Henceforth, we assume that~$f$ belongs to~$L^2(\QT)$.
The space--time virtual element method for the approximation of solutions to~\eqref{continuous-weak} reads
\begin{equation} \label{VEM}
\begin{cases}
\text{find } \uh \in \Yh \text{ such that}\\
\bh(\uh,\vh) = (f, \Pizpmo \vh)_{0,\Omega} \qquad \forall \vh \in \Yh.
\end{cases}
\end{equation}
For Cartesian-product-in-time meshes, the well posedness of the method and \emph{a priori} error estimates
were proven in~\cite[Sections~3 and~4]{Gomez-Mascotto-Moiola-Perugia:2022} based on a discrete inf-sup argument.
That analysis can be extended to the meshes introduced in Section~\ref{subsection:intro-meshes}, using the norms 
\begin{align*}
\Norm{v}_{\Ytaun}^2 & := \sum_{\Kfrak \in \taun} \Norm{v}_{\YE}^2, \\
\Norm{v}_{\Xtaun}^2 & := \Norm{v}_{\Ytaun}^2 + \Norm{\Newtonh \Pistar v}_{\Ytaun}^2 \\
& \quad + \frac{\cH}{2} 
\Big(
\Norm{\Pistar v(\cdot, 0)}_{0, \Omega}^2 +
\!\!\!\!\!\!\!\!\sum_{
\scriptsize
\begin{tabular}{c}
$\Ex \in \FcalExE$ \\
$\Ex \not\subset \Omega \times\{0\}$
\end{tabular}
}
\Norm{\UWEx (v)}_{0, \Ex}^2 
+ \Norm{\Pistar v(\cdot, T)}_{0, \Omega}^2 
\Big),
\end{align*}
where the discrete Newton potential~$\Newtonh : \calS_p(\taun) \rightarrow \Yh$ is defined as follows: for all~$\vh \in \Yh$, 
\begin{equation*}
\begin{split}
    \ah(\Newtonh \phih, \vh) = \cH \Big( (\dpt \phih, \vh)_{0, \QT} + \left(\phih(\cdot, 0), \vh(\cdot, 0)\right)_{0, \Omega} + 
\!\!\!\!\!\!\!\!\!\! \sum_{
\scriptsize
\begin{tabular}{c}
$\Ex \in \FcalExE$ \\
$\Ex \not\subset \Omega \times\{0\}$
\end{tabular}
}
\Big(\UWEx(\phih), \vh{}_{|_{\Kfrak^+}}\Big)_{0, \Ex} \Big).
\end{split}
\end{equation*}
%%%%%
In particular, we have the following result.
\begin{thm} \label{theorem:inf-sup-well-posedness}
There exists~$\gammaI>0$ independent of~$\taun$
such that
\[
\sup_{0\ne \vh \in \Yh} \frac{\bh(\uh,\vh)}{\Norm{\vh}_{\Ytaun}}
\ge \gammaI \Norm{\uh}_{\Xtaun}
\qquad\qquad \forall \uh \in \Yh.
\]
Therefore, method~\eqref{VEM} is well posed.
An inspection of \cite[Proposition~3.2]{Gomez-Mascotto-Moiola-Perugia:2022} reveals that
the constant~$\gammaI$ only depends on the stability constants~$\alpha_*$ and~$\alpha^*$
in~\eqref{stability-bounds:local}.
\end{thm}

%%%%%%%%%%%%%%%%%%%%%%%%%%%%%%%%%%%%%%%%%%%%%%%%%%%%%%%%%%%%%%%%%%%%%%%%%%
\section{Handling space--time meshes and variable degrees} \label{section:handling-meshes-degrees}
%%%%%%%%%%%%%%%%%%%%%%%%%%%%%%%%%%%%%%%%%%%%%%%%%%%%%%%%%%%%%%%%%%%%%%%%%%
We focus on the handling of the general prismatic mesh structure.
More precisely, we describe refinements procedures in Section~\ref{subsection:mesh-refinement};
a time-slab flagging strategy to split method~\eqref{VEM}
into smaller linear systems in Section~\ref{subsection:flag}; 
an element-topology flagging strategy in Section~\ref{subsection:flagging}.
Further, we discuss space--time virtual elements with variable degrees of accuracy
in Section~\ref{subsection:hp-spaces}.

%%%%%%%%%%%%%
\subsection{Mesh refinements} \label{subsection:mesh-refinement}
%%%%%%%%%%%%%
We describe a procedure to refine general prismatic space--time elements.
For the sake of presentation,
we consider the $(1+1)$ dimensional case only;
the extension to any spatial dimensions
follows with a minor effort.

As discussed in Section~\ref{subsection:intro-meshes},
each mesh consists of rectangular elements~$\Kfrak$,
with boundary given by the union of four (two space-like, two time-like) straight segments.
Each straight segment may be the union of aligned edges of the element,
which constitute the boundary of an associated closed polytope~$K$; see Section~\ref{subsection:intro-meshes}.
Similar to~\cite[Figure~7]{Cangiani-Georgulis-Pryer-Sutton:2017},
regardless of the number of existing hanging nodes from previous refinements,
a given~$\E$ to be refined is split into four siblings
by connecting the centroid of~$\E$ with the midpoints of each straight segment of the boundary.

\begin{figure}[h]
\centering
\begin{minipage}{0.32\textwidth}
\begin{center}
\begin{tikzpicture}[scale=1]
\draw[black, thick, -] (0,0) -- (4,0) -- (4,4) -- (0,4) -- (0,0);
\draw[fill=red] (0,0) circle (2pt);
\draw[fill=red] (4,0) circle (2pt);
\draw[fill=red] (4,4) circle (2pt);
\draw[fill=red] (0,4) circle (2pt);
\draw[fill=red] (4,2) circle (2pt);
\draw[fill=red] (4,1) circle (2pt);
\draw[fill=red] (4,.5) circle (2pt);
\end{tikzpicture}
\end{center}
\end{minipage}
\qquad\qquad\qquad
\begin{minipage}{0.32\textwidth}
\begin{center}
\begin{tikzpicture}[scale=1]
\draw[black, thick, -] (0,0) -- (4,0) -- (4,4) -- (0,4) -- (0,0);
\draw[black, thick, -] (2,0) -- (2,4);
\draw[black, thick, -] (0,2) -- (4,2);
\draw[fill=blue] (2,0) circle (2pt);
\draw[fill=blue] (0,2) circle (2pt);
\draw[fill=blue] (2,4) circle (2pt);
\draw[fill=blue] (2,2) circle (2pt);
\draw[fill=red] (0,0) circle (2pt);
\draw[fill=red] (4,0) circle (2pt);
\draw[fill=red] (4,4) circle (2pt);
\draw[fill=red] (0,4) circle (2pt);
\draw[fill=red] (4,2) circle (2pt);
\draw[fill=red] (4,1) circle (2pt);
\draw[fill=red] (4,.5) circle (2pt);
\end{tikzpicture}
\end{center}
\end{minipage}
\caption{Element refinement strategy.
We connect the centroid of the element with the midpoints of each straight segment of the boundary,
regardless of the presence of previously generated hanging nodes.
The red dots denote the nodes of~$\E$;
the blue dots denote the newly created nodes.}
\label{figure:element-refinement}
\end{figure}
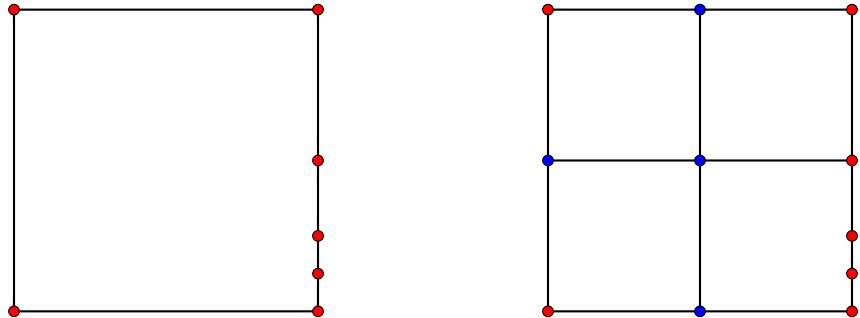

In Figure~\ref{figure:element-refinement}, we show an example of an element refinement.
By this procedure, at most new five nodes are generated, 
fewer in presence of previously generated hanging nodes.

%%%%%%%%%%%%%
\subsection{Time-slab flagging strategy} \label{subsection:flag}
%%%%%%%%%%%%%
We present a flagging strategy that allows for the decomposition
of the algebraic linear system stemming from~\eqref{VEM}
into smaller linear systems.
To this aim, we assume that the first mesh of a sequence~$\{\taun\}$
is a ``tensor-product-in-time'' mesh,
which can be arranged into time-slabs.
On this mesh, method~\eqref{VEM} can be assembled and solved sequentially with respect to the time-slabs,
i.e., it can be interpreted as a time-stepping scheme.

Assume now that we are given a mesh at the refinement step~$n-1$,
which is split into time-slabs.
We explain here how to identify a time-slab partition of the refined mesh at step~$n$.
The time interval~$[0,T]$
can be partitioned as $0 =: t_0 < \dots < t_{\ell} := T$,
so that there exists at least one element~$\Kfrak = \Kfrakx \times \Kfrakt$ of the mesh at step~$n$
with~$\Kfrakx$ contained in~$\Omega \times \{t_j\}$
for a~$j=0,\dots,\ell-1$, and with~$\Kfrakt = (t_{i},t_{j})$
for given~$i,j=0,\dots,\ell$, $i<j$.
The time-slabs may be generated cycling on the~$t_j$, $j=1,\dots,\ell-1$:
if all the elements in the mesh at step $n$
are either below or above~$t_j$,
then a new time-slab is created by flagging all the elements that lie
below the given~$t_j$
but have not been already allocated in a previously identified time-slab.

\begin{figure}[H]
\centering
\begin{minipage}{0.30\textwidth}
\begin{center}
\begin{tikzpicture}[scale=.7]
\draw[black, ->] (-0.5,0) -- (5,0);
\draw[black, ->] (0,-0.5) -- (0,5);
\draw (5,-.2) node[black, left] {\small x}; 
\draw (-.1,4.5) node[black, left] {\small t}; 
\draw[black, very thick, -] (0,0) -- (4,0) -- (4,4) -- (0,4) -- (0,0);
\draw[black, very thick, -] (2,0) -- (2,4);
\draw[black, very thick, -] (0,2) -- (4,2);
\draw[black, very thick, -] (1,0) -- (1,2);
\draw[black, very thick, -] (0,1) -- (2,1);
\node[align=left] at (.5,.5) {$i$};
\node[align=left] at (1.5,.5) {$i$};
\node[align=left] at (.5,1.5) {$i$};
\node[align=left] at (1.5,1.5) {$i$};
\node[align=left] at (3,1) {$i$};
\node[align=left] at (1,3) {$ii$};
\node[align=left] at (3,3) {$ii$};
\end{tikzpicture}
\end{center}
\end{minipage}
%%%%%%%%%%%%%%%%%%%%%%%%%%%%%%%%
\begin{minipage}{0.30\textwidth}
\begin{center}
\begin{tikzpicture}[scale=.7]
\draw[black, ->] (-0.5,0) -- (5,0);
\draw[black, ->] (0,-0.5) -- (0,5);
\draw (5,-.2) node[black, left] {\small x}; 
\draw (-.1,4.5) node[black, left] {\small t}; 
\draw[black, very thick, -] (0,0) -- (4,0) -- (4,4) -- (0,4) -- (0,0);
\draw[black, very thick, -] (2,0) -- (2,4);
\draw[black, very thick, -] (0,2) -- (4,2);
\draw[black, very thick, -] (1,0) -- (1,2);
\draw[black, very thick, -] (3,0) -- (3,2);
\draw[black, very thick, -] (0,1) -- (4,1);
\node[align=left] at (.5,.5) {$i$};
\node[align=left] at (1.5,.5) {$i$};
\node[align=left] at (.5,1.5) {$ii$};
\node[align=left] at (1.5,1.5) {$ii$};
\node[align=left] at (.5+2,.5) {$i$};
\node[align=left] at (1.5+2,.5) {$i$};
\node[align=left] at (.5+2,1.5) {$ii$};
\node[align=left] at (1.5+2,1.5) {$ii$};
\node[align=left] at (1,3) {$iii$};
\node[align=left] at (3,3) {$iii$};
\end{tikzpicture}
\end{center}
\end{minipage}
\caption{
\textbf{Left panel:} we start with a given space--time mesh
with prescribed time-slab structure and flagging.
\textbf{Right panel:}
time-slab flagging based on the proposed strategy after one refinement step.}
\label{figure:1st-flagging}
\end{figure}
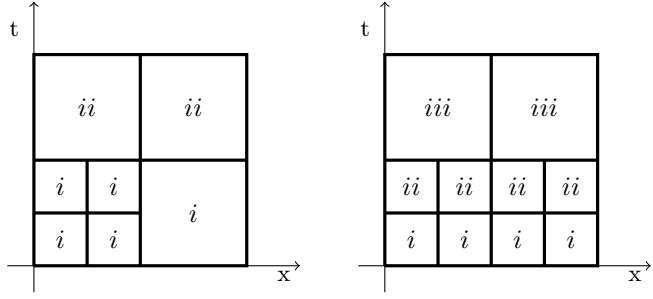

Based on the new flagging, at the $n$-th refinement step,
we can assemble and solve method~\eqref{VEM} sequentially
on the newly created time-slabs.
In Figure~\ref{figure:1st-flagging}, we illustrate this flagging procedure with an example.

%%

%%%%%%%%%%%%
\subsection{Element-topology flagging strategy} \label{subsection:flagging}
%%%%%%%%%%%%
A common issue in the implementation of polytopic methods is the lack of a reference element;
this renders the computation of the local matrices
more expensive than for methods
based on simplicial or Cartesian meshes.

Assume that the elements of a given mesh~$\taun$
can be grouped into a uniformly finite number of equivalence classes (up to dilations and translations),
i.e., there exists a set~$\mathcal C$ of ``reference elements'',
with~$\card(\mathcal C)$ bounded uniformly for all meshes,
such that each element~$\E$ of~$\taun$
is equivalent to a ``reference element'' in~$\mathcal C$.
In this case, local matrices need to be computed only for the ``reference elements'',
thus hastening the assembling of the final system.
In adaptive mesh refinements,
different ``reference elements'' may appear.
For this reason, we introduce a flag associated with the element topology that identifies the corresponding ``reference element''.
Similar ideas were used in~\cite{Frittelli-Madzvamuse-Sgura:2023} in a space-only VEM context.

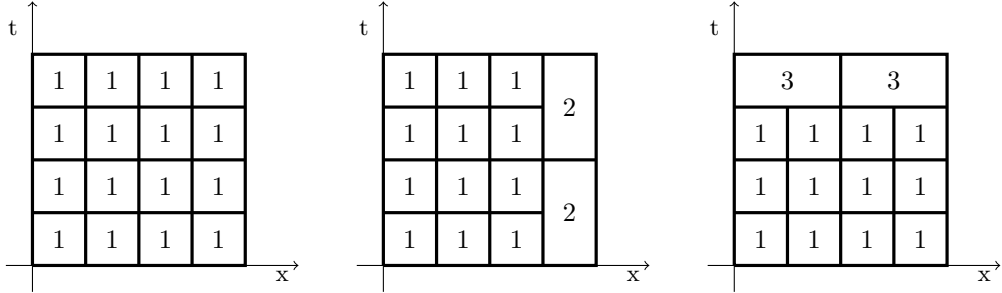
\begin{figure}[H]
\centering
\begin{minipage}{0.30\textwidth}
\begin{center}
\begin{tikzpicture}[scale=.7]
\draw[black, ->] (-0.5,0) -- (5,0);
\draw[black, ->] (0,-0.5) -- (0,5);
\draw (5,-.2) node[black, left] {\small x}; 
\draw (-.1,4.5) node[black, left] {\small t}; 
\draw[black, very thick, -] (0,0) -- (4,0) -- (4,4) -- (0,4) -- (0,0);
\draw[black, very thick, -] (1,0) -- (1,4);
\draw[black, very thick, -] (2,0) -- (2,4);
\draw[black, very thick, -] (3,0) -- (3,4);
\draw[black, very thick, -] (0,1) -- (4,1);
\draw[black, very thick, -] (0,2) -- (4,2);
\draw[black, very thick, -] (0,3) -- (4,3);
\node[align=left] at (.5,.5) {$1$}; \node[align=left] at (1.5,.5) {$1$}; \node[align=left] at (2.5,.5) {$1$}; \node[align=left] at (3.5,.5) {$1$};
\node[align=left] at (.5,1.5) {$1$}; \node[align=left] at (1.5,1.5) {$1$}; \node[align=left] at (2.5,1.5) {$1$}; \node[align=left] at (3.5,1.5) {$1$};
\node[align=left] at (.5,2.5) {$1$}; \node[align=left] at (1.5,2.5) {$1$}; \node[align=left] at (2.5,2.5) {$1$}; \node[align=left] at (3.5,2.5) {$1$};
\node[align=left] at (.5,3.5) {$1$}; \node[align=left] at (1.5,3.5) {$1$}; \node[align=left] at (2.5,3.5) {$1$}; \node[align=left] at (3.5,3.5) {$1$};
\end{tikzpicture}
\end{center}
\end{minipage}
\begin{minipage}{0.30\textwidth}
\begin{center}
\begin{tikzpicture}[scale=.7]
\draw[black, ->] (-0.5,0) -- (5,0);
\draw[black, ->] (0,-0.5) -- (0,5);
\draw (5,-.2) node[black, left] {\small x}; 
\draw (-.1,4.5) node[black, left] {\small t}; 
\draw[black, very thick, -] (0,0) -- (4,0) -- (4,4) -- (0,4) -- (0,0);
\draw[black, very thick, -] (1,0) -- (1,4);
\draw[black, very thick, -] (2,0) -- (2,4);
\draw[black, very thick, -] (3,0) -- (3,4);
\draw[black, very thick, -] (0,1) -- (3,1);
\draw[black, very thick, -] (0,2) -- (4,2);
\draw[black, very thick, -] (0,3) -- (3,3);
\node[align=left] at (.5,.5) {$1$}; \node[align=left] at (1.5,.5) {$1$}; \node[align=left] at (2.5,.5) {$1$};
\node[align=left] at (.5,1.5) {$1$}; \node[align=left] at (1.5,1.5) {$1$}; \node[align=left] at (2.5,1.5) {$1$};
\node[align=left] at (.5,2.5) {$1$}; \node[align=left] at (1.5,2.5) {$1$}; \node[align=left] at (2.5,2.5) {$1$};
\node[align=left] at (.5,3.5) {$1$}; \node[align=left] at (1.5,3.5) {$1$}; \node[align=left] at (2.5,3.5) {$1$};
\node[align=left] at (3.5,1) {$2$}; \node[align=left] at (3.5,3) {$2$};
\end{tikzpicture}\end{center}
\end{minipage}
\begin{minipage}{0.30\textwidth}
\begin{center}
\begin{tikzpicture}[scale=.7]
\draw[black, ->] (-0.5,0) -- (5,0);
\draw[black, ->] (0,-0.5) -- (0,5);
\draw (5,-.2) node[black, left] {\small x}; 
\draw (-.1,4.5) node[black, left] {\small t}; 
\draw[black, very thick, -] (0,0) -- (4,0) -- (4,4) -- (0,4) -- (0,0);
\draw[black, very thick, -] (1,0) -- (1,3);
\draw[black, very thick, -] (2,0) -- (2,4);
\draw[black, very thick, -] (3,0) -- (3,3);
\draw[black, very thick, -] (0,1) -- (4,1);
\draw[black, very thick, -] (0,2) -- (4,2);
\draw[black, very thick, -] (0,3) -- (4,3);
\node[align=left] at (.5,.5) {$1$}; \node[align=left] at (1.5,.5) {$1$}; \node[align=left] at (2.5,.5) {$1$}; \node[align=left] at (3.5,.5) {$1$};
\node[align=left] at (.5,1.5) {$1$}; \node[align=left] at (1.5,1.5) {$1$}; \node[align=left] at (2.5,1.5) {$1$}; \node[align=left] at (3.5,1.5) {$1$};
\node[align=left] at (.5,2.5) {$1$}; \node[align=left] at (1.5,2.5) {$1$}; \node[align=left] at (2.5,2.5) {$1$}; \node[align=left] at (3.5,2.5) {$1$};
\node[align=left] at (1,3.5) {$3$}; \node[align=left] at (3,3.5) {$3$};
\end{tikzpicture}
\end{center}
\end{minipage}
\caption{\textbf{Left panel:} a mesh consisting of equivalent elements.
\textbf{Central panel:} a mesh with nonmatching time-like facets.
\textbf{Right panel:} a mesh with nonmatching space-like facets.}
\label{figure:meshes}
\end{figure}

For instance, in Figure~\ref{figure:meshes} (left panel),
we only have one ``reference element'' (a space--time square);
in Figure~\ref{figure:meshes} (central panel),
we have two types of elements (a space--time square and a space--time pentagon with two time-like facets on the left;
in Figure~\ref{figure:meshes} (right panel),
we have two types of elements (a space--time square and a space--time pentagon with two space-like facets at the bottom.
Inside each element, we denote
the associated ``topology-flag'' with a natural number.

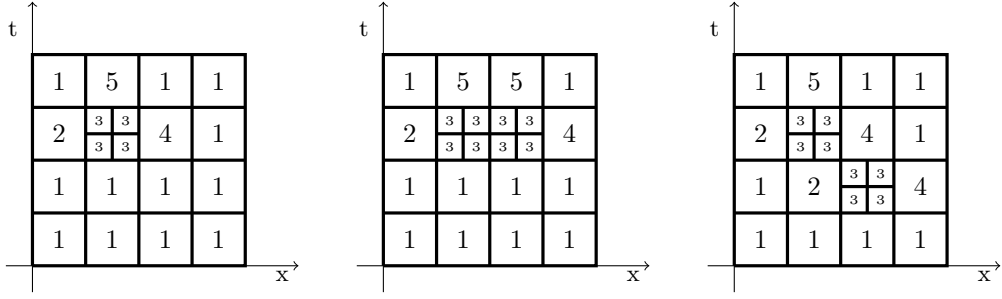
\begin{figure}[H]
\centering
\begin{minipage}{0.30\textwidth}
\begin{center}
\begin{tikzpicture}[scale=.7]
\draw[black, ->] (-0.5,0) -- (5,0);
\draw[black, ->] (0,-0.5) -- (0,5);
\draw (5,-.2) node[black, left] {\small x}; 
\draw (-.1,4.5) node[black, left] {\small t}; 
\draw[black, very thick, -] (0,0) -- (4,0) -- (4,4) -- (0,4) -- (0,0);
\draw[black, very thick, -] (1,0) -- (1,4);
\draw[black, very thick, -] (2,0) -- (2,4);
\draw[black, very thick, -] (3,0) -- (3,4);
\draw[black, very thick, -] (0,1) -- (4,1);
\draw[black, very thick, -] (0,2) -- (4,2);
\draw[black, very thick, -] (0,3) -- (4,3);
\draw[black, very thick, -] (1,2.5) -- (2,2.5);
\draw[black, very thick, -] (1.5,2) -- (1.5,3);
\node[align=left] at (.5,.5) {$1$}; \node[align=left] at (1.5,.5) {$1$}; \node[align=left] at (2.5,.5) {$1$}; \node[align=left] at (3.5,.5) {$1$};
\node[align=left] at (.5,1.5) {$1$}; \node[align=left] at (1.5,1.5) {$1$}; \node[align=left] at (2.5,1.5) {$1$}; \node[align=left] at (3.5,1.5) {$1$};
\node[align=left] at (.5,2.5) {$2$}; \node[align=left] at (2.5,2.5) {$4$}; \node[align=left] at (3.5,2.5) {$1$};
\node[align=left] at (.5,3.5) {$1$}; \node[align=left] at (1.5,3.5) {$5$}; \node[align=left] at (2.5,3.5) {$1$}; \node[align=left] at (3.5,3.5) {$1$};
\node[align=left] at (1.25,2.25) {\tiny{$3$}}; \node[align=left] at (1.25,2.75) {\tiny{$3$}}; \node[align=left] at (1.75,2.25) {\tiny{$3$}}; \node[align=left] at (1.75,2.75) {\tiny{$3$}};
\end{tikzpicture}
\end{center}
\end{minipage}
\begin{minipage}{0.30\textwidth}
\begin{center}
\begin{tikzpicture}[scale=.7]
\draw[black, ->] (-0.5,0) -- (5,0);
\draw[black, ->] (0,-0.5) -- (0,5);
\draw (5,-.2) node[black, left] {\small x}; 
\draw (-.1,4.5) node[black, left] {\small t}; 
\draw[black, very thick, -] (0,0) -- (4,0) -- (4,4) -- (0,4) -- (0,0);
\draw[black, very thick, -] (1,0) -- (1,4);
\draw[black, very thick, -] (2,0) -- (2,4);
\draw[black, very thick, -] (3,0) -- (3,4);
\draw[black, very thick, -] (0,1) -- (4,1);
\draw[black, very thick, -] (0,2) -- (4,2);
\draw[black, very thick, -] (0,3) -- (4,3);
\draw[black, very thick, -] (1,2.5) -- (2,2.5);
\draw[black, very thick, -] (1.5,2) -- (1.5,3);
\draw[black, very thick, -] (2,2.5) -- (3,2.5);
\draw[black, very thick, -] (2.5,2) -- (2.5,3);
\node[align=left] at (.5,.5) {$1$}; \node[align=left] at (1.5,.5) {$1$}; \node[align=left] at (2.5,.5) {$1$}; \node[align=left] at (3.5,.5) {$1$};
\node[align=left] at (.5,1.5) {$1$}; \node[align=left] at (1.5,1.5) {$1$}; \node[align=left] at (2.5,1.5) {$1$}; \node[align=left] at (3.5,1.5) {$1$};
\node[align=left] at (.5,2.5) {$2$}; \node[align=left] at (3.5,2.5) {$4$};
\node[align=left] at (.5,3.5) {$1$}; \node[align=left] at (1.5,3.5) {$5$}; \node[align=left] at (2.5,3.5) {$5$}; \node[align=left] at (3.5,3.5) {$1$};
\node[align=left] at (1.25,2.25) {\tiny{$3$}}; \node[align=left] at (1.25,2.75) {\tiny{$3$}}; \node[align=left] at (1.75,2.25) {\tiny{$3$}}; \node[align=left] at (1.75,2.75) {\tiny{$3$}};
\node[align=left] at (2.25,2.25) {\tiny{$3$}}; \node[align=left] at (2.25,2.75) {\tiny{$3$}}; \node[align=left] at (2.75,2.25) {\tiny{$3$}}; \node[align=left] at (2.75,2.75) {\tiny{$3$}};
\end{tikzpicture}\end{center}
\end{minipage}
\begin{minipage}{0.30\textwidth}
\begin{center}
\begin{tikzpicture}[scale=.7]
\draw[black, ->] (-0.5,0) -- (5,0);
\draw[black, ->] (0,-0.5) -- (0,5);
\draw (5,-.2) node[black, left] {\small x}; 
\draw (-.1,4.5) node[black, left] {\small t}; 
\draw[black, very thick, -] (0,0) -- (4,0) -- (4,4) -- (0,4) -- (0,0);
\draw[black, very thick, -] (1,0) -- (1,4);
\draw[black, very thick, -] (2,0) -- (2,4);
\draw[black, very thick, -] (3,0) -- (3,4);
\draw[black, very thick, -] (0,1) -- (4,1);
\draw[black, very thick, -] (0,2) -- (4,2);
\draw[black, very thick, -] (0,3) -- (4,3);
\draw[black, very thick, -] (1,2.5) -- (2,2.5);
\draw[black, very thick, -] (1.5,2) -- (1.5,3);
\draw[black, very thick, -] (2,1.5) -- (3,1.5);
\draw[black, very thick, -] (2.5,1) -- (2.5,2);
\node[align=left] at (.5,.5) {$1$}; \node[align=left] at (1.5,.5) {$1$}; \node[align=left] at (2.5,.5) {$1$}; \node[align=left] at (3.5,.5) {$1$};
\node[align=left] at (.5,1.5) {$1$}; \node[align=left] at (1.5,1.5) {$2$}; \node[align=left] at (3.5,1.5) {$4$};
\node[align=left] at (.5,2.5) {$2$}; \node[align=left] at (2.5,2.5) {$4$}; \node[align=left] at (3.5,2.5) {$1$};
\node[align=left] at (.5,3.5) {$1$}; \node[align=left] at (1.5,3.5) {$5$}; \node[align=left] at (2.5,3.5) {$1$}; \node[align=left] at (3.5,3.5) {$1$};
\node[align=left] at (1.25,2.25) {\tiny{$3$}}; \node[align=left] at (1.25,2.75) {\tiny{$3$}}; \node[align=left] at (1.75,2.25) {\tiny{$3$}}; \node[align=left] at (1.75,2.75) {\tiny{$3$}};
\node[align=left] at (2.25,1.25) {\tiny{$3$}}; \node[align=left] at (2.25,1.75) {\tiny{$3$}}; \node[align=left] at (2.75,1.25) {\tiny{$3$}}; \node[align=left] at (2.75,1.75) {\tiny{$3$}};
\end{tikzpicture}
\end{center}
\end{minipage}
\caption{Three different meshes with element-topology flagging.}
\label{figure:meshes2}
\end{figure}

%%%%%
\noindent In Figure~\ref{figure:meshes2},
we present three situations that
are more elaborated than those in Figure~\ref{figure:meshes}.

The above definition of equivalence classes refers to the geometry of the mesh elements. This is enough to define a flagging strategy in case of a uniform degree of accuracy. 
For the case of variable degrees of accuracy in Section~\ref{subsection:hp-spaces} below,
this definition 
can be 
extended taking into account the distribution of degrees of accuracy. 
Namely, we flag two elements in the same way if the following conditions are satisfied:
they have the same geometric flag;
they have the same degree of accuracy;
the degrees of accuracy assigned to their time-like facets are the same.
The existence of a maximum degree of accuracy
and a maximum number of time-like facets for all elements
would allow for a uniformly  bounded number of reference elements in the extended sense. 
Without a fixed maximum degree of accuracy
(e.g., in the $p$-version of the method, where $p\to\infty$),
one cannot expect to have a finite number of reference elements;
this is also the case for the~$p$-version of any standard
finite element discretization.

%%%

%%%%%%%%%%
\subsection{Space--time virtual element spaces with variable degrees of accuracy} \label{subsection:hp-spaces}
%%%%%%%%%%
Let~$\taun$ be a given space--time polytopic mesh consisting of~$\NE$ 
elements and~$\pbf \in \Nbb^{\NE}$ be a given distribution of degrees of accuracy.
More precisely, we sort the elements of~$\taun$ as~$\{\Kfrak_j\}_{j=1}^{\NE}$,
and denote the degree of accuracy in each element~$\Kfrak_j$ by~$\p_{j}$,
$j=1,\dots,\NE$.

Given the vector~$\pbf$,
we fix the degrees of freedom associated with each element
according to the following maximum strategy:
\begin{itemize}
    \item in each element~$\Kfrak_j$, $j=1,\dots,\NE$,
    we take bulk moments~\eqref{bulk-dofs}
    up to degree~$\p_j-1$;
    \item on each internal \emph{time-like} facet~$\F$ shared by two different elements~$\Kfrak_j$ and~$\Kfrak_\ell$
    for given~$j,\ell=1,\dots,\NE$,
    we take \emph{time-like} moments~\eqref{vertical-dofs} up to degree~$\max(\p_j,\p_\ell)$;
    \item on each boundary \emph{time-like} facet~$\F$ on the boundary of the element~$\Kfrak_j$
    for a given $j=1,\dots,\NE$,
    we take \emph{time-like} moments~\eqref{vertical-dofs}
    up to degree~$\p_j$;
    \item if~$\Kfrak_j = \Kfrakx{}_{,j} \times \Kfrakt{}_{, j}$
    for a given $j =1,\dots,\NE$,
    we take the \emph{space-like} moments~\eqref{bottom-dofs} on~$\Kfrakx{}_{,j}$ up to degree~$\p_j$.
\end{itemize}
We collect the \emph{time-like} polynomial degrees in the vector~$\pbfFcalht \in \Nbb^{\card(\Fcalht)}$
and the \emph{space-like} polynomial degrees in the vector~$\pbfFcalhEx \in \Nbb^{\NE}$.
We order the \emph{time-like} facets in~$\Fcalh$ as~$\{\F_j\}_{j=1}^{\card(\Fcalht)}$.

Given~$\ctildeHE$ and~$\nutildeE$ as in~\eqref{ctilde-nutilde},
the corresponding local space on~$\E_j$ reads
\[
\begin{split}
\Vh(\E_j) := \Big\{v \in L^2(\Kfrak_j) \mid\ 
& \ctildeHE \dpt \vh - \nutildeE \Deltax \vh  \in \Pp{\p_j-1}{\Kfrak_j};\;
\vh{}_{|\Kfrakx{}_{,j}} \in \Pp{\p_j^{\text{space}}}{\Kfrakx{}_{,j}};\\[0.2em]
& \nbf_{\F_k, \x} \cdot \nablax \vh{}_{|\F_k} \in \Pp{\p_k^{\text{time}}}{\F_k} \, \forall \F_k \in \FcalEtj \Big\}.
\end{split}
\]
The global space~$\Yh$ is constructed
by the nonconforming coupling of the \emph{time-like} degrees of freedom~\eqref{vertical-dofs}.
An immediate consequence of this maximum strategy
is that $\Pbb_{\p_j}(\Kfrak_j)$ is contained in~$\Vh(\E_j)$.

\begin{figure}[H]
\centering
\begin{minipage}{0.30\textwidth}
\begin{center}
\begin{tikzpicture}[scale=.7]
\draw[black, ->] (-0.5,0) -- (5,0);
\draw[black, ->] (0,-0.5) -- (0,5);
\draw (5,-.2) node[black, left] {\small x}; 
\draw (-.1,4.5) node[black, left] {\small t}; 
\draw[black, very thick, -] (0,0) -- (4,0) -- (4,4) -- (0,4) -- (0,0);
\draw[black, very thick, -] (2,0) -- (2,4);
\draw[black, very thick, -] (0,2) -- (4,2);
\draw [fill=white] (1,1) circle[radius= 0.5 em]; \node[align=left] at (1,1) {\orange{\tiny{$1$}}};
\draw [fill=white] (3,1) circle[radius= 0.5 em]; \node[align=left] at (3,1) {\orange{\tiny{$2$}}};
\draw [fill=white] (1,3) circle[radius= 0.5 em]; \node[align=left] at (1,3) {\orange{\tiny{$3$}}};
\draw [fill=white] (3,3) circle[radius= 0.5 em]; \node[align=left] at (3,3) {\orange{\tiny{$4$}}};
\end{tikzpicture}
\end{center}
\end{minipage}
\begin{minipage}{0.30\textwidth}
\begin{center}
\begin{tikzpicture}[scale=.7]
\draw[black, ->] (-0.5,0) -- (5,0);
\draw[black, ->] (0,-0.5) -- (0,5);
\draw (5,-.2) node[black, left] {\small x}; 
\draw (-.1,4.5) node[black, left] {\small t}; 
\draw[black, very thick, -] (0,0) -- (4,0) -- (4,4) -- (0,4) -- (0,0);
\draw[black, very thick, -] (2,0) -- (2,4);
\draw[black, very thick, -] (0,2) -- (4,2);
\draw [fill=white] (1,0) circle[radius= 0.5 em]; \node[align=left] at (1,0) {\blue{\tiny{$1$}}};
\draw [fill=white] (3,0) circle[radius= 0.5 em]; \node[align=left] at (3,0) {\blue{\tiny{$2$}}};
\draw [fill=white] (1,2) circle[radius= 0.5 em]; \node[align=left] at (1,2) {\blue{\tiny{$3$}}};
\draw [fill=white] (3,2) circle[radius= 0.5 em]; \node[align=left] at (3,2) {\blue{\tiny{$4$}}};
\draw [fill=white] (0,1) circle[radius= 0.5 em]; \node[align=left] at (0,1) {\red{\tiny{$1$}}};
\draw [fill=white] (2,1) circle[radius= 0.5 em]; \node[align=left] at (2,1) {\red{\tiny{$2$}}};
\draw [fill=white] (4,1) circle[radius= 0.5 em]; \node[align=left] at (4,1) {\red{\tiny{$2$}}};
\draw [fill=white] (0,3) circle[radius= 0.5 em]; \node[align=left] at (0,3) {\red{\tiny{$3$}}};
\draw [fill=white] (2,3) circle[radius= 0.5 em]; \node[align=left] at (2,3) {\red{\tiny{$4$}}};
\draw [fill=white] (4,3) circle[radius= 0.5 em]; \node[align=left] at (4,3) {\red{\tiny{$4$}}};
\end{tikzpicture}\end{center}
\end{minipage}
\begin{minipage}{0.30\textwidth}
\begin{center}
\begin{tikzpicture}[scale=.7]
\draw[black, ->] (-0.5,0) -- (5,0);
\draw[black, ->] (0,-0.5) -- (0,5);
\draw (5,-.2) node[black, left] {\small x}; 
\draw (-.1,4.5) node[black, left] {\small t}; 
\draw[black, very thick, -] (0,0) -- (4,0) -- (4,4) -- (0,4) -- (0,0);
\draw[black, very thick, -] (2,0) -- (2,4);
\draw[black, very thick, -] (0,2) -- (4,2);
\draw [fill=blue] (2/3,0) circle[radius= 0.3 em]; \draw [fill=blue] (4/3,0) circle[radius= 0.3 em];
\draw [fill=blue] (2+1/2,0) circle[radius= 0.3 em]; \draw [fill=blue] (2+2/2,0) circle[radius= 0.3 em]; \draw [fill=blue] (2+3/2,0) circle[radius= 0.3 em];
\draw [fill=blue] (2/5,2) circle[radius= 0.3 em]; \draw [fill=blue] (4/5,2) circle[radius= 0.3 em]; \draw [fill=blue] (6/5,2) circle[radius= 0.3 em]; \draw [fill=blue] (8/5,2) circle[radius= 0.3 em];
\draw [fill=blue] (2+2/6,2) circle[radius= 0.3 em]; \draw [fill=blue] (2+4/6,2) circle[radius= 0.3 em]; \draw [fill=blue] (2+6/6,2) circle[radius= 0.3 em]; \draw [fill=blue] (2+8/6,2) circle[radius= 0.3 em]; \draw [fill=blue] (2+10/6,2) circle[radius= 0.3 em];
\draw [fill=red] (0,2/3) circle[radius= 0.3 em]; \draw [fill=red] (0,4/3) circle[radius= 0.3 em];
\draw [fill=red] (2,2/4) circle[radius= 0.3 em]; \draw [fill=red] (2,4/4) circle[radius= 0.3 em]; \draw [fill=red] (2,6/4) circle[radius= 0.3 em]; 
\draw [fill=red] (4,2/4) circle[radius= 0.3 em]; \draw [fill=red] (4,4/4) circle[radius= 0.3 em]; \draw [fill=red] (4,6/4) circle[radius= 0.3 em];
\draw [fill=red] (0,2+2/5) circle[radius= 0.3 em]; \draw [fill=red] (0,2+4/5) circle[radius= 0.3 em]; \draw [fill=red] (0,2+6/5) circle[radius= 0.3 em]; \draw [fill=red] (0,2+8/5) circle[radius= 0.3 em]; 
\draw [fill=red] (2,2+2/6) circle[radius= 0.3 em]; \draw [fill=red] (2,2+4/6) circle[radius= 0.3 em]; \draw [fill=red] (2,2+6/6) circle[radius= 0.3 em]; \draw [fill=red] (2,2+8/6) circle[radius= 0.3 em]; \draw [fill=red] (2,2+10/6) circle[radius= 0.3 em]; 
\draw [fill=red] (4,2+2/6) circle[radius= 0.3 em]; \draw [fill=red] (4,2+4/6) circle[radius= 0.3 em]; \draw [fill=red] (4,2+6/6) circle[radius= 0.3 em]; \draw [fill=red] (4,2+8/6) circle[radius= 0.3 em]; \draw [fill=red] (4,2+10/6) circle[radius= 0.3 em];
\draw [fill=orange] (1,1) circle[radius= 0.3 em];
\draw [fill=orange] (3,4/3) circle[radius= 0.3 em]; \draw [fill=orange] (3-2/6,2/3) circle[radius= 0.3 em]; \draw [fill=orange] (3+2/6,2/3) circle[radius= 0.3 em];
\draw [fill=orange] (1,2+6/4) circle[radius= 0.3 em]; \draw [fill=orange] (1-2/10,2+4/4) circle[radius= 0.3 em]; \draw [fill=orange] (1+2/10,2+4/4) circle[radius= 0.3 em]; \draw [fill=orange] (1-2/4,2+2/4) circle[radius= 0.3 em]; \draw [fill=orange] (1,2+2/4) circle[radius= 0.3 em]; \draw [fill=orange] (1+2/4,2+2/4) circle[radius= 0.3 em];
\draw [fill=orange] (3,2+8/5) circle[radius= 0.3 em]; \draw [fill=orange] (3-2/10,2+6/5) circle[radius= 0.3 em]; \draw [fill=orange] (3+2/10,2+6/5) circle[radius= 0.3 em]; \draw [fill=orange] (3-2/4,2+4/5) circle[radius= 0.3 em]; \draw [fill=orange] (3,2+4/5) circle[radius= 0.3 em]; \draw [fill=orange] (3+2/4,2+4/5) circle[radius= 0.3 em]; \draw [fill=orange] (3-8/10,2+2/5) circle[radius= 0.3 em]; \draw [fill=orange] (3-2/10,2+2/5) circle[radius= 0.3 em]; \draw [fill=orange] (3+2/10,2+2/5) circle[radius= 0.3 em]; \draw [fill=orange] (3+8/10,2+2/5) circle[radius= 0.3 em]; 
\end{tikzpicture}
\end{center}
\end{minipage}
\caption{\textbf{Left panel:} initial distribution of degrees of accuracy over the elements.
\textbf{Central panel:} ``polynomial degrees'' on \emph{space-like} (\blue{\bf blue})
and \emph{time-like} (\red{\bf red}) facets.
\textbf{Right panel:} corresponding degrees of freedom;
the \orange{\bf orange} dots denote the bulk moments;
the \blue{\bf blue} dots denote the \emph{space-like} moments;
the \red{\bf red} dots denote the \emph{time-like} moments.}
\label{figure:hp-meshes}
\end{figure}
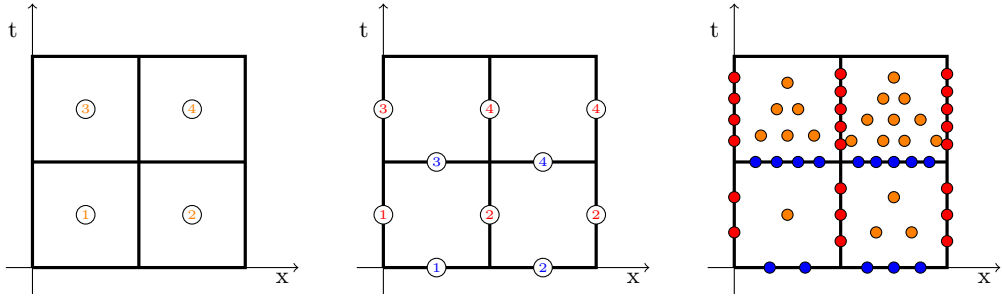

To illustrate the maximum strategy,
we provide an example in Figure~\ref{figure:hp-meshes},
where we consider a uniform Cartesian mesh of~$4$ elements
with different degrees of accuracy.

%%%%%%%%%%%%%%%%%%%%%%%%%%%%%%%%%%%%%%%%%%%%%%%%%%%%%%%%%%%%%%%%%%%%%%%%%%
\section{Numerical investigation: convergence tests} \label{section:numerical-a-priori}
%%%%%%%%%%%%%%%%%%%%%%%%%%%%%%%%%%%%%%%%%%%%%%%%%%%%%%%%%%%%%%%%%%%%%%%%%%
We assess the convergence of the $\h$- and $\h\p$-versions of the virtual element method (VEM) in~\eqref{VEM}.
Since the virtual element solution~$\uh$ is not known in closed form and the error in the $\Xtaun$
norm is not computable,
we report the following error quantities:
given~$\uh$ the solution to~\eqref{VEM},
\footnotesize{\begin{equation} \label{exact-errors}
\begin{split}
& \EcalY
 := \Norm{u-\PiN \uh}_{\Ytaun},
\qquad\qquad
\EcalN := \Norm{\PiN(\Newtonh \Pistar (u - \uh))}_{\Ytaun}, \\
& (\EcalU)^2 :=  \frac{\cH}{2} 
\Bigg(
\Norm{\Pistar(u - \uh)(\cdot, 0)}_{0, \Omega}^2 + \!\!\!\!\!\!\!\!\sum_{
\scriptsize 
\begin{tabular}{c}
$\Ex \in \FcalExE$ \\
$\Ex \not\subset \Omega \times\{0\}$
\end{tabular}
}\!\!\!\!\!\!\!\!
\Norm{\UWEx (\Pistar(u - \uh))}_{0, \Ex}^2  
+ \Norm{\Pistar (u - \uh)(\cdot, T)}_{0, \Omega}^2 
\Bigg),  \\
& \big( \EcalX \big)^2
 := \big( \EcalY \big)^2
 + \big( \EcalN \big)^2
+ \big( \EcalU \big)^2. 
\end{split}
\end{equation}
}\normalsize

%%%%%%%%%%%%%
\subsection{Test cases} \label{subsection:test-cases}
%%%%%%%%%%%%%
We consider test cases with coefficients~$\nu=1$ and~$\cH=1$.
The right-hand side~$f$, and the boundary and initial conditions
are computed accordingly to the exact solutions below.

\paragraph*{Test case 1.}
We define the analytic function
\begin{flalign} \label{test-case-1}
    && u_1(x,t) := \exp(-t) \sin(\pi x)
     &&& \forall (x,t) \in \QT := (0,1) \times (0,1).
\end{flalign}
\paragraph*{Test case 2.}
For~$\alpha > \frac12$, we define the function
\begin{flalign} \label{test-case-2}
    && u_2(x,t):=  \sin(\pi x) \ t^\alpha
    &&& \forall (x,t) \in \QT := (0,1) \times (0,0.1),
\end{flalign}
which belongs to~$H^{\alpha+1/2-\varepsilon}(0,1; \mathcal C^{\infty}(\Omega))$, $\varepsilon>0$. 
\paragraph*{Test case 3.}
We define the function
\small{\begin{flalign} \label{test-case-3}
&& u_3(x,t):= \sum_{n=0}^\infty \frac{4}{(2n+1)\pi}
              \sin ( (2n + 1)\pi \ x) 
              \exp (-(2n + 1)^2 \pi^2 t)
              &&& \forall (x,t) \in \QT := (0,1) \times (0,1),
\end{flalign}}\normalsize
which is the Fourier series of the solution to~\eqref{continuous-strong}
with zero source term~$f$, initial condition~$u_0 = 1$,
and homogeneous Dirichlet boundary conditions~$g$.
In the numerical experiments, the series in~\eqref{test-case-2} is truncated at $n=250$.
The function $u_3$ belongs to $H^s(0, 1; H^1_0 (0, 1))$
for any $s < 1/4$ and to~$H^1(0, 1; H^{-1}(\Omega)) \cap L^2(0, 1; H_0^1(\Omega)) \backslash H^1(\QT)$;
see~\cite{Schotzau_Schwab:2000}.
In particular, $u_3$ is singular at the interface of the (incompatible) initial and boundary conditions.

%%%%%%%%%%%%%
\subsection{$\h$- and $\h\p$-versions for singular solutions} \label{subsection:numerics-hp-version}
%%%%%%%%%%%%%
The performance of the~$\h$-version of the method on smooth solutions
was investigated in~\cite{Gomez-Mascotto-Moiola-Perugia:2022}
and is therefore omitted here.
The test case 1 is used in Section~\ref{section:residual-type} below.

We focus on the convergence of the~$\h$- and~$\h\p$-versions of the method
for singular solutions.
To that aim, we consider the test cases 2 and 3.
Notably, we want to assess exponential convergence
in terms of the cubic root of the number of degrees of freedom for the $\h\p$-version
on certain geometrically refined space--time meshes.

\begin{figure}[h]
\centering
\begin{minipage}{0.2\textwidth}
\begin{center}
\begin{tikzpicture}[scale=0.035]
\draw[fill=red, opacity=0.65] (0,0) -- (100,0) -- (100,100) -- (0,100) -- (0,0);
\draw[black, thick, -] (0,0) -- (100,0) -- (100,100) -- (0,100) -- (0,0);
\draw[black, thick, -] (5,0) -- (5,100);\draw[black, thick, -] (10,0) -- (10,100);\draw[black, thick, -] (15,0) -- (15,100);\draw[black, thick, -] (20,0) -- (20,100);\draw[black, thick, -] (25,0) -- (25,100);\draw[black, thick, -] (30,0) -- (30,100);\draw[black, thick, -] (35,0) -- (35,100);\draw[black, thick, -] (40,0) -- (40,100);\draw[black, thick, -] (45,0) -- (45,100);\draw[black, thick, -] (50,0) -- (50,100);\draw[black, thick, -] (55,0) -- (55,100);\draw[black, thick, -] (60,0) -- (60,100);\draw[black, thick, -] (65,0) -- (65,100);\draw[black, thick, -] (70,0) -- (70,100);\draw[black, thick, -] (75,0) -- (75,100);\draw[black, thick, -] (80,0) -- (80,100);\draw[black, thick, -] (85,0) -- (85,100);\draw[black, thick, -] (90,0) -- (90,100);\draw[black, thick, -] (95,0) -- (95,100);
\end{tikzpicture}
\end{center}
\end{minipage}
\qquad\qquad\qquad
\begin{minipage}{0.2\textwidth}
\begin{center}
\begin{tikzpicture}[scale=0.035]
\draw[fill=blue, opacity=0.65] (0,10) -- (100,10) -- (100,100) -- (0,100) -- (0,10);
\draw[fill=red, opacity=0.65] (0,0) -- (100,0) -- (100,10) -- (0,10) -- (0,0);
\draw[black, thick, -] (0,0) -- (100,0) -- (100,100) -- (0,100) -- (0,0);
\draw[black, thick, -] (5,0) -- (5,100);\draw[black, thick, -] (10,0) -- (10,100);\draw[black, thick, -] (15,0) -- (15,100);\draw[black, thick, -] (20,0) -- (20,100);\draw[black, thick, -] (25,0) -- (25,100);\draw[black, thick, -] (30,0) -- (30,100);\draw[black, thick, -] (35,0) -- (35,100);\draw[black, thick, -] (40,0) -- (40,100);\draw[black, thick, -] (45,0) -- (45,100);\draw[black, thick, -] (50,0) -- (50,100);\draw[black, thick, -] (55,0) -- (55,100);\draw[black, thick, -] (60,0) -- (60,100);\draw[black, thick, -] (65,0) -- (65,100);\draw[black, thick, -] (70,0) -- (70,100);\draw[black, thick, -] (75,0) -- (75,100);\draw[black, thick, -] (80,0) -- (80,100);\draw[black, thick, -] (85,0) -- (85,100);\draw[black, thick, -] (90,0) -- (90,100);\draw[black, thick, -] (95,0) -- (95,100);
\draw[black, thick, -] (0,10) -- (100,10);
\end{tikzpicture}
\end{center}
\end{minipage}
\qquad\qquad\qquad
\begin{minipage}{0.2\textwidth}
\begin{center}
\begin{tikzpicture}[scale=0.035]
\draw[fill=green, opacity=0.65] (0,10) -- (100,10) -- (100,100) -- (0,100) -- (0,10);
\draw[fill=blue, opacity=0.65] (0,1) -- (100,1) -- (100,10) -- (0,10) -- (0,1);
\draw[fill=red, opacity=0.65] (0,0) -- (100,0) -- (100,1) -- (0,1) -- (0,0);
\draw[black, thick, -] (0,0) -- (100,0) -- (100,100) -- (0,100) -- (0,0);
\draw[black, thick, -] (5,0) -- (5,100);\draw[black, thick, -] (10,0) -- (10,100);\draw[black, thick, -] (15,0) -- (15,100);\draw[black, thick, -] (20,0) -- (20,100);\draw[black, thick, -] (25,0) -- (25,100);\draw[black, thick, -] (30,0) -- (30,100);\draw[black, thick, -] (35,0) -- (35,100);\draw[black, thick, -] (40,0) -- (40,100);\draw[black, thick, -] (45,0) -- (45,100);\draw[black, thick, -] (50,0) -- (50,100);\draw[black, thick, -] (55,0) -- (55,100);\draw[black, thick, -] (60,0) -- (60,100);\draw[black, thick, -] (65,0) -- (65,100);\draw[black, thick, -] (70,0) -- (70,100);\draw[black, thick, -] (75,0) -- (75,100);\draw[black, thick, -] (80,0) -- (80,100);\draw[black, thick, -] (85,0) -- (85,100);\draw[black, thick, -] (90,0) -- (90,100);\draw[black, thick, -] (95,0) -- (95,100);
\draw[black, thick, -] (0,10) -- (100,10);
\draw[black, thick, -] (0,1) -- (100,1);
\end{tikzpicture}
\end{center}
\end{minipage}
\caption{First three meshes employed in the $\h\p$ refinements for the test case~2 with exact solution~$u_2$ in~\eqref{test-case-2}.
The space--time domain is~$\QT= (0,1)\times(0,0.1)$.
For a better understanding of the figure, we scale the $t$-coordinates by~10.
In colours, we represent the local degrees of accuracy:
\red{\bf red:} $\p=1$; \blue{\bf blue:} $\p=2$; \darkgreen{\bf green:} $\p=3$.}
\label{figure:test-case-2-hp-meshes}
\end{figure}
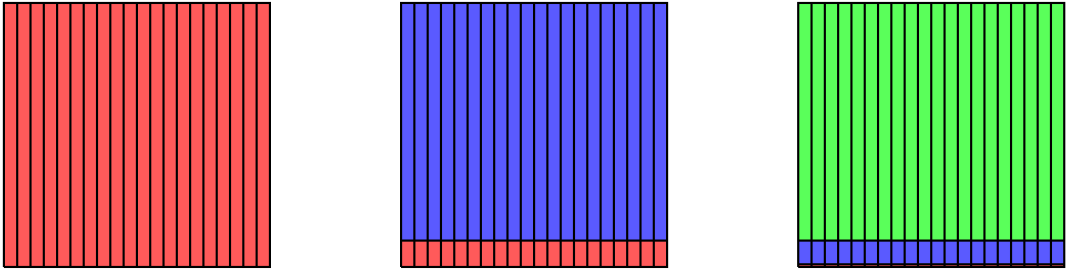

First, we consider the test case 2.
For the~$h$-version of the method, we consider uniform degree of accuracy~$p = 1$
and a sequence of uniform Cartesian space--time meshes with~$h_t = 2h_x =  0.2\times 2^{-i}, \ i = 1, \ldots, 6$.
For the~$hp$-version of the method, we proceed similarly as in~\cite[Example 2]{Cangiani-Dong-Georgoulis:2017}:
we fix a partition of the spatial domain with~$h_x = 0.05$
and consider a sequence of  temporal meshes geometrically graded towards~$t = 0$ with grading factor~$\sigma_t = 0.1$.
In addition, the degree of accuracy~$p$ is increased by~$1$ from one time slab to the next one. 
In Figure~\ref{figure:test-case-2-hp-meshes},
we depict the first three meshes with varying degrees of accuracy. 

\begin{figure}[!ht]
    \centering
    \includegraphics[width = 0.45\textwidth]{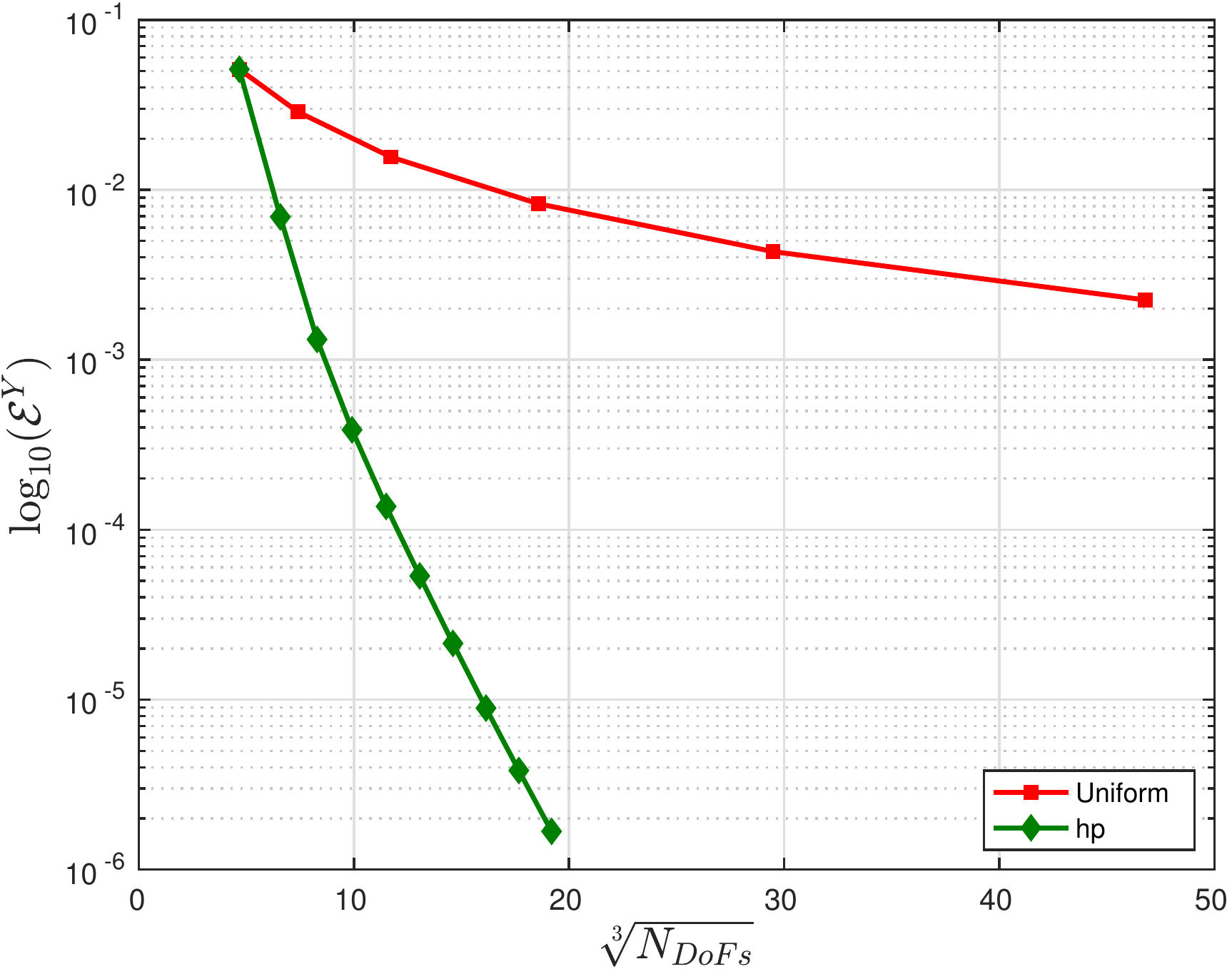}
    \hspace{0.2in}
    \includegraphics[width = 0.45\textwidth]{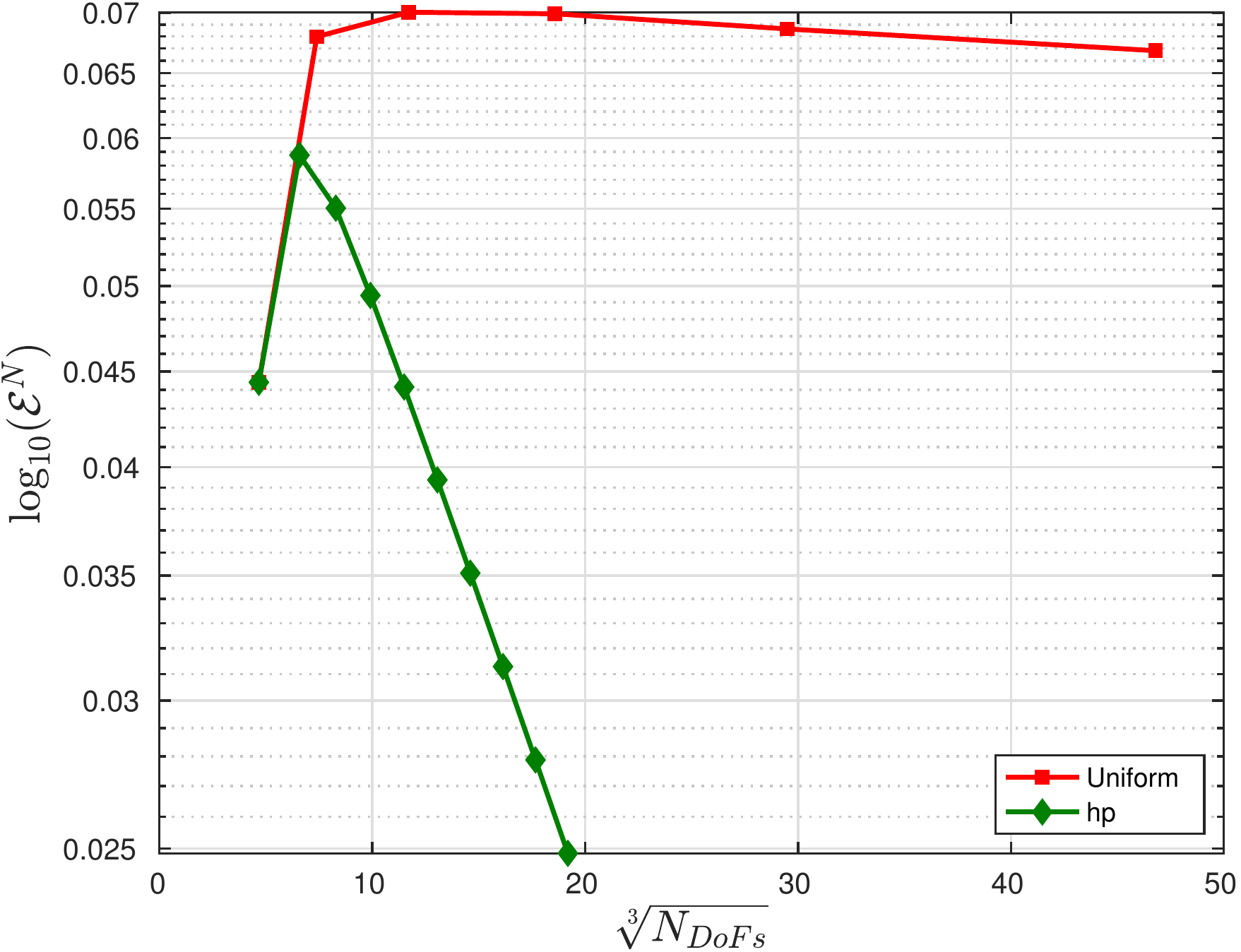}\\
    \includegraphics[width = 0.45\textwidth]{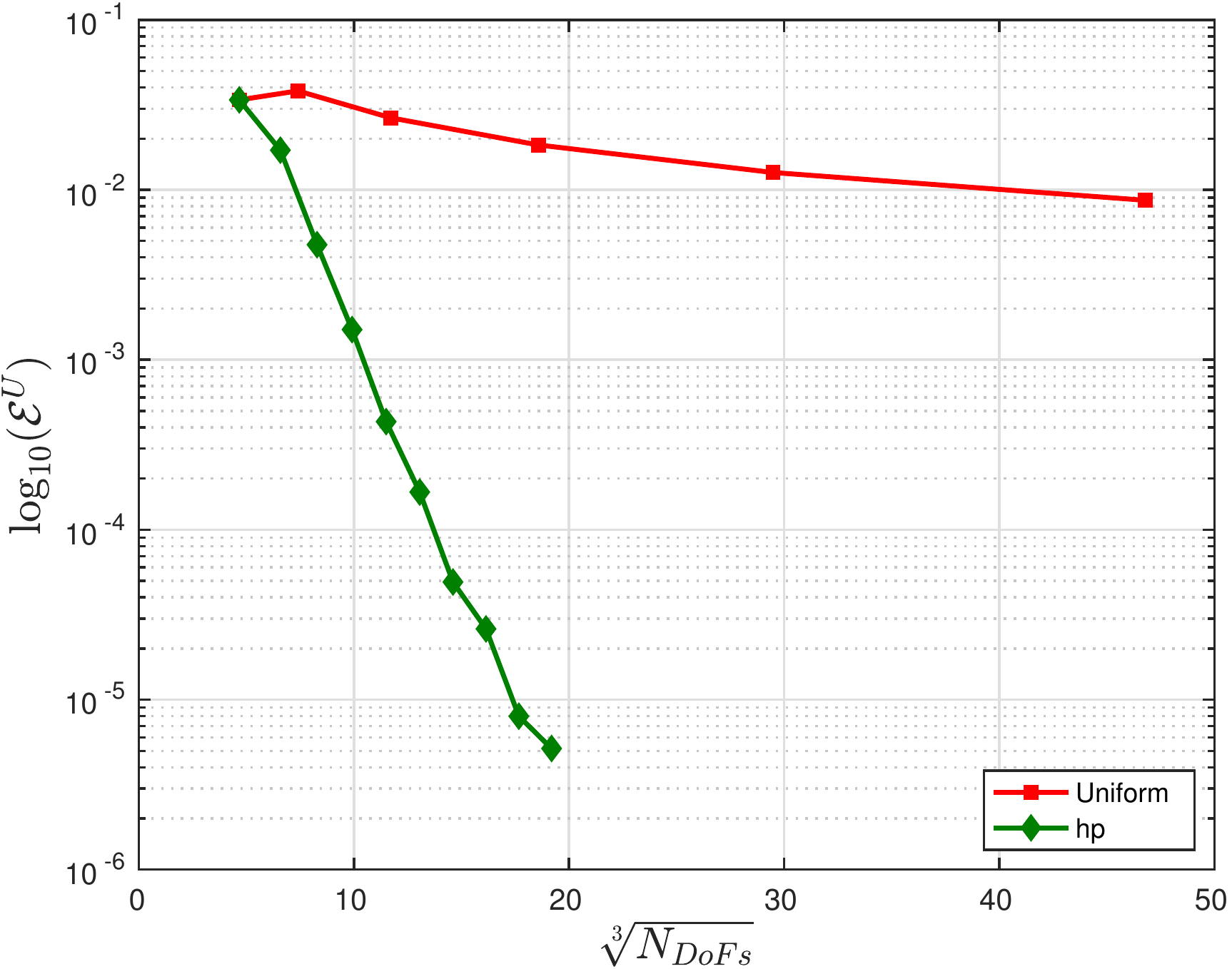}
    \hspace{0.2in}
    \includegraphics[width = 0.45\textwidth]{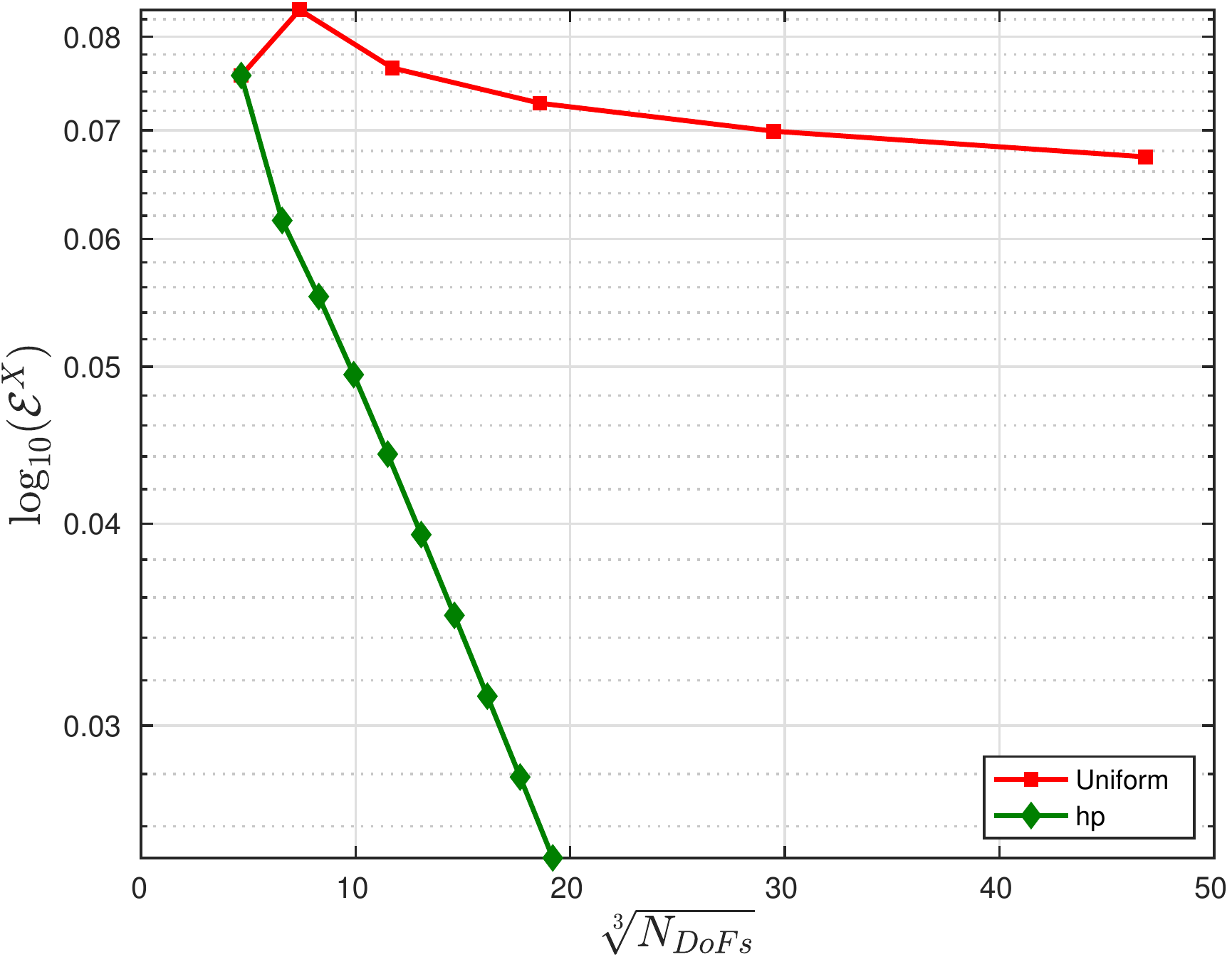}
    \caption{Convergence of the errors in~\eqref{exact-errors}
    for the~$h$- and~$hp$-versions of the method.
    We consider the test case~2 with exact solution~$u_2$ in~\eqref{test-case-2}, $\alpha = 0.55$.
    \label{FIG::hp-Singularity-t0-a55}}
    \end{figure}

\begin{figure}[!ht]
    \centering
    \includegraphics[width = 0.45\textwidth]{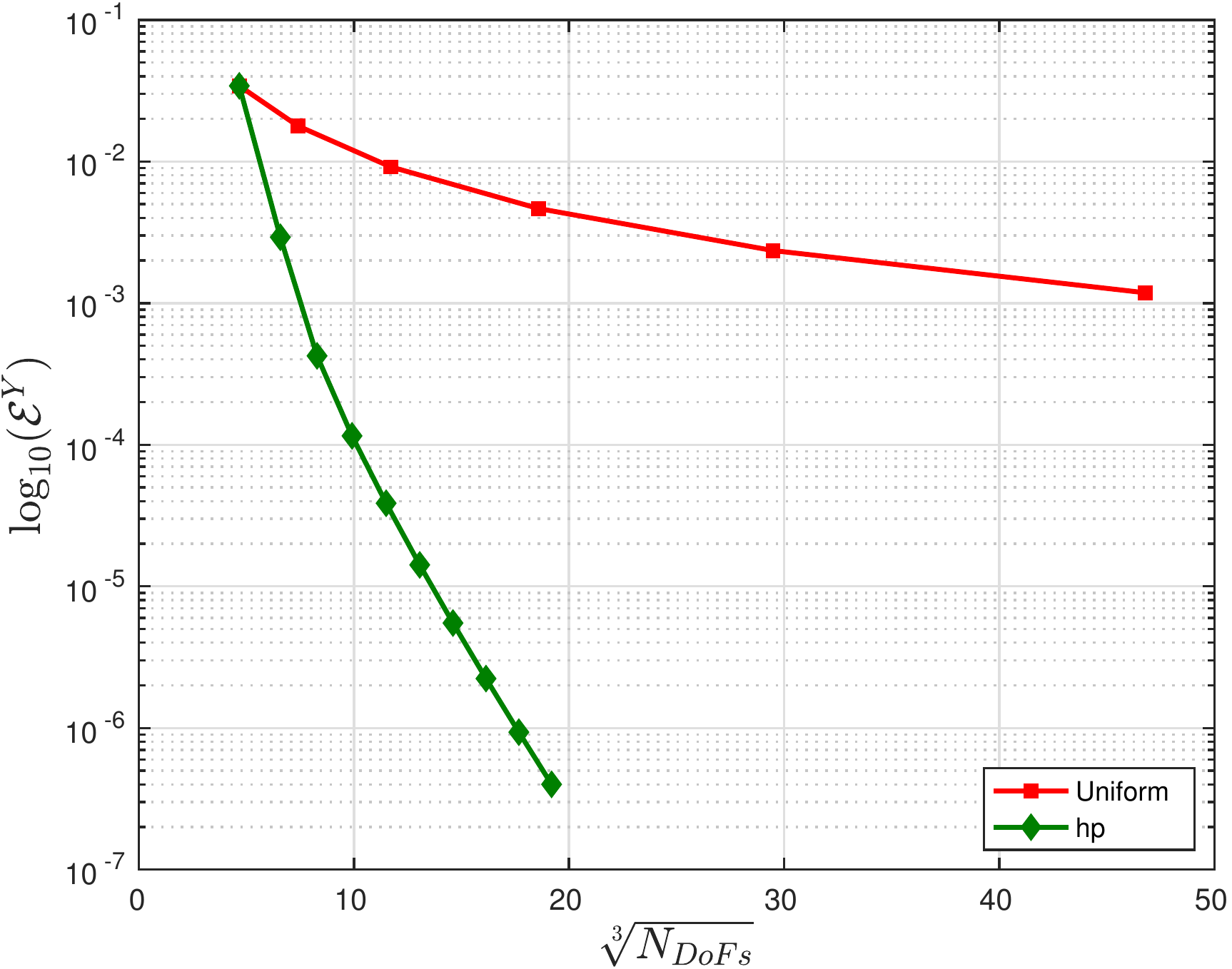}
    \hspace{0.2in}
    \includegraphics[width = 0.45\textwidth]{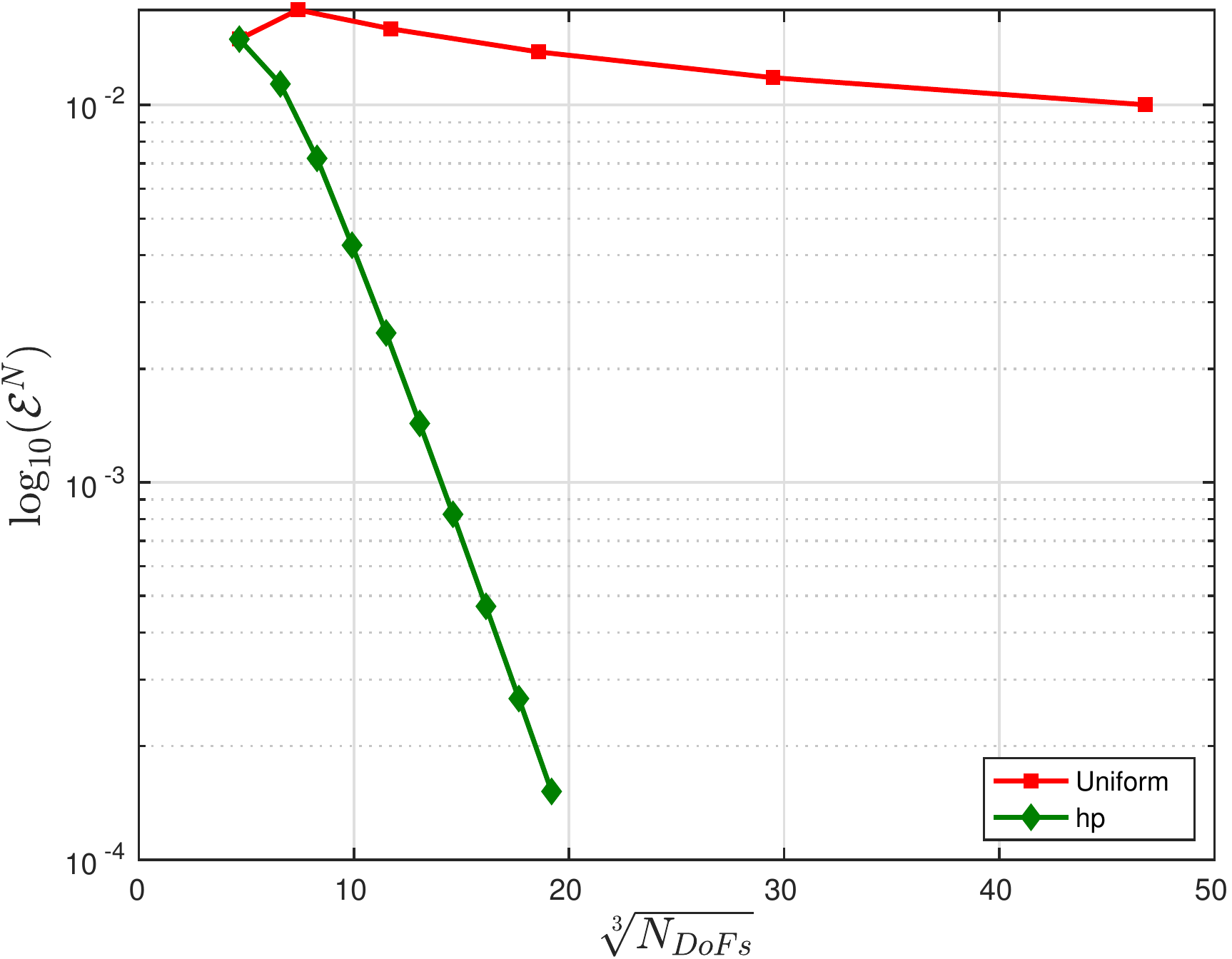}\\
    \includegraphics[width = 0.45\textwidth]{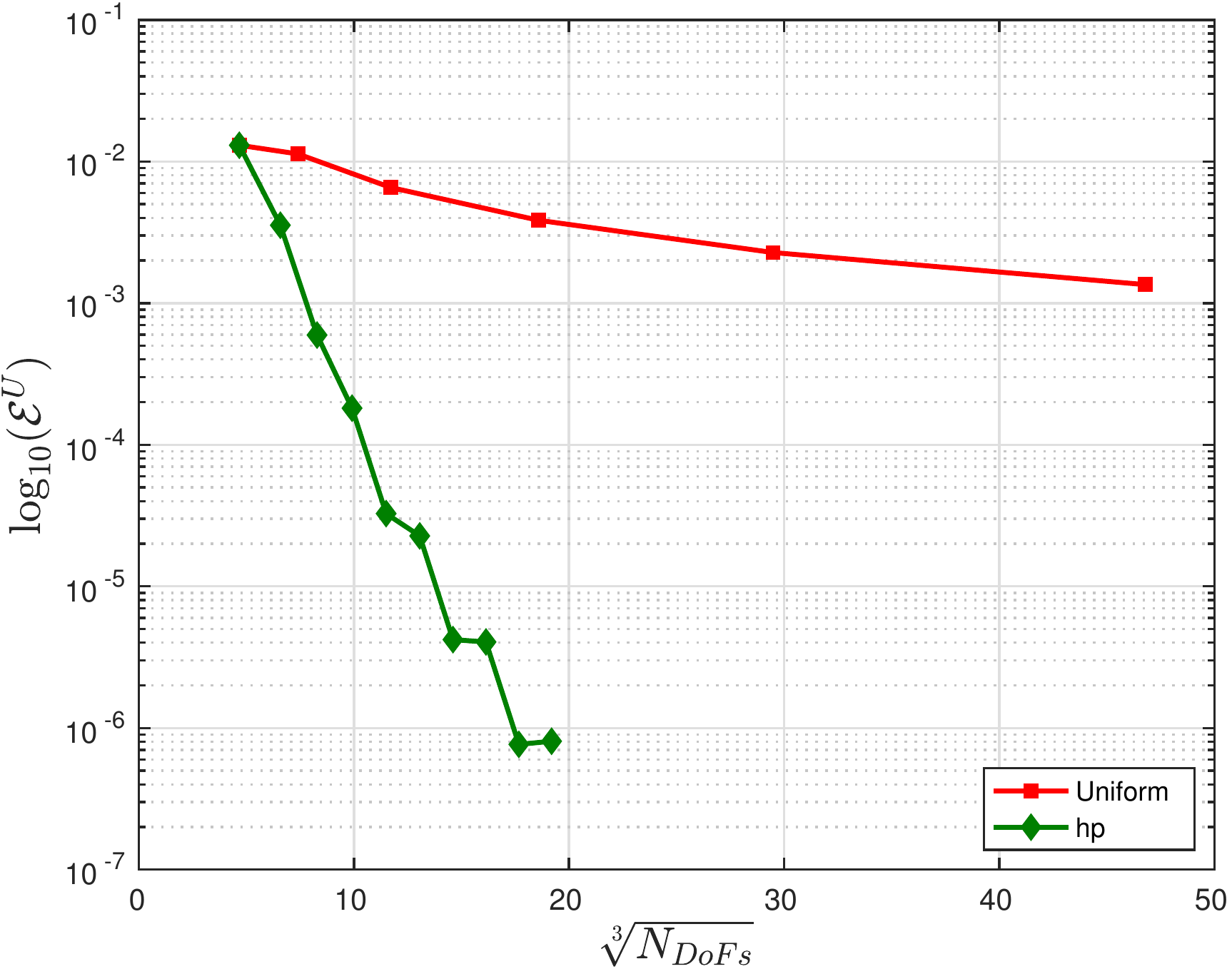}
    \hspace{0.2in}
    \includegraphics[width = 0.45\textwidth]{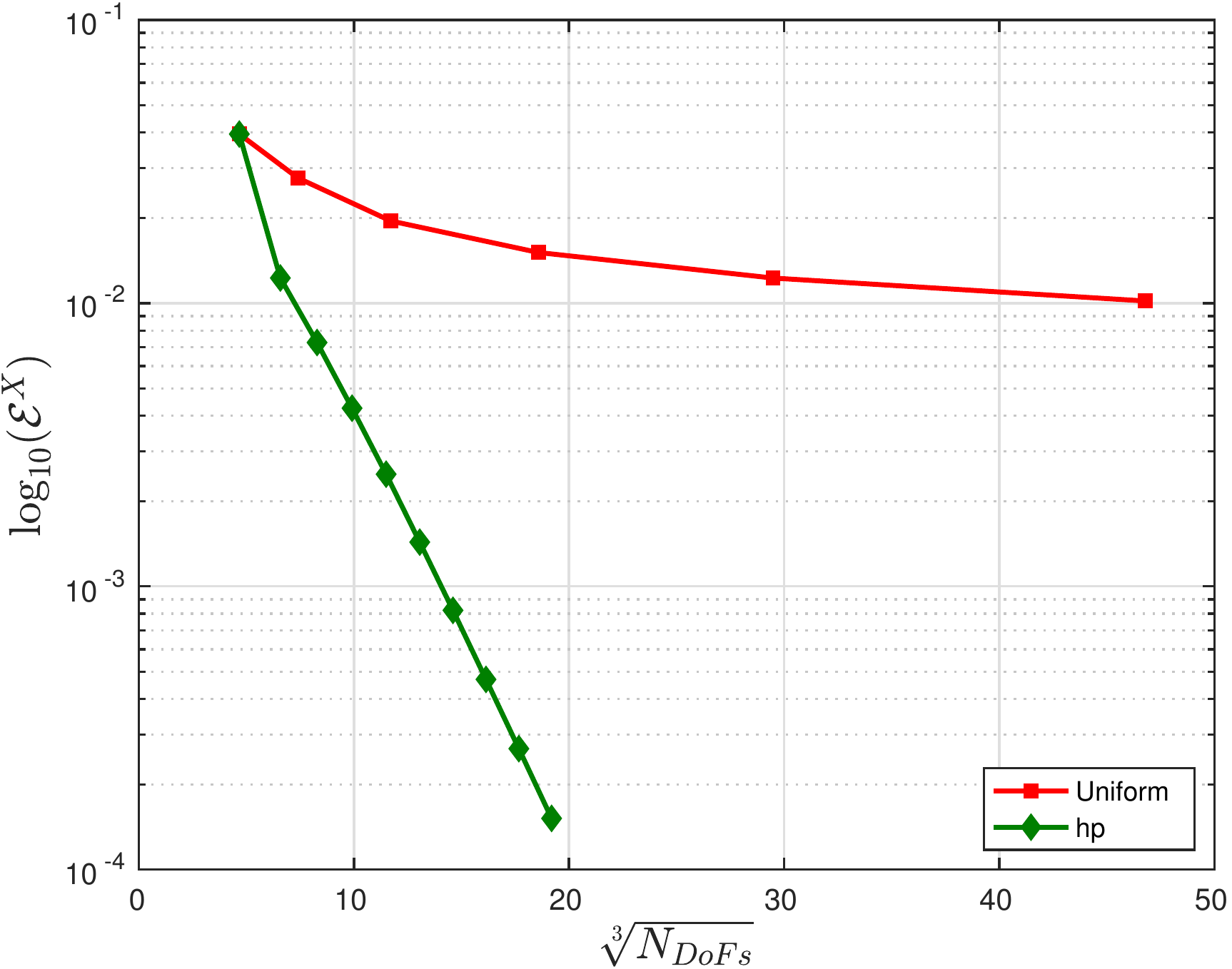}
    \caption{Convergence of the errors in~\eqref{exact-errors}
    for the~$h$- and~$hp$-versions of the method.
    We consider the test case~2 with exact solution~$u_2$ in~\eqref{test-case-2}, $\alpha = 0.75$.
    \label{FIG::hp-Singularity-t0-a75}}
\end{figure}

In Figures~\ref{FIG::hp-Singularity-t0-a55} and~\ref{FIG::hp-Singularity-t0-a75},
we show the errors in~\eqref{exact-errors} in \emph{semilogy} scale for~$\alpha = 0.55$ and~$\alpha = 0.75$, respectively.
Exponential convergence in terms of the cubic root of the number of degrees of freedom
is observed for the~$hp$-version of the method and both values of~$\alpha$.
In all cases, the $\h$-version is outperformed and displays only an algebraic decay of the error.

\medskip 

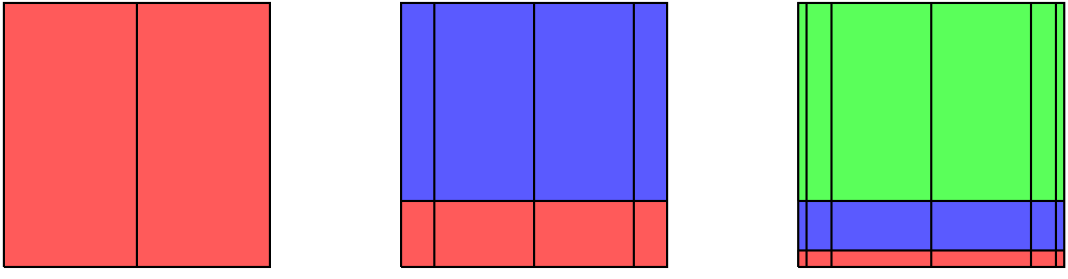
\begin{figure}[h]
\centering
\begin{minipage}{0.2\textwidth}
\begin{center}
\begin{tikzpicture}[scale=0.035]
\draw[fill=red, opacity=0.65] (0,0) -- (100,0) -- (100,100) -- (0,100) -- (0,0);
\draw[black, thick, -] (0,0) -- (100,0) -- (100,100) -- (0,100) -- (0,0);
\draw[black, thick, -] (50,0) -- (50,100);
\end{tikzpicture}
\end{center}
\end{minipage}
\qquad\qquad\qquad
\begin{minipage}{0.2\textwidth}
\begin{center}
\begin{tikzpicture}[scale=0.035]
\draw[fill=red, opacity=0.65] (0,0) -- (100,0) -- (100,100) -- (0,100) -- (0,0);
\draw[fill=white] (0,25) -- (100,25) -- (100,100) -- (0,100) -- (0,25);
\draw[fill=blue, opacity=0.65] (0,25) -- (100,25) -- (100,100) -- (0,100) -- (0,25);
\draw[black, thick, -] (0,0) -- (100,0) -- (100,100) -- (0,100) -- (0,0);
\draw[black, thick, -] (50,0) -- (50,100);
\draw[black, thick, -] (12.5,0) -- (12.5,100);
\draw[black, thick, -] (100-12.5,0) -- (100-12.5,100);
\draw[black, thick, -] (0,25) -- (100,25);
\end{tikzpicture}
\end{center}
\end{minipage}
\qquad\qquad\qquad
\begin{minipage}{0.2\textwidth}
\begin{center}
\begin{tikzpicture}[scale=0.035]
\draw[fill=blue, opacity=0.65] (0,0) -- (100,0) -- (100,100) -- (0,100) -- (0,0);
\draw[fill=white] (0,25) -- (100,25) -- (100,100) -- (0,100) -- (0,25);
\draw[fill=green, opacity=0.65] (0,25) -- (100,25) -- (100,100) -- (0,100) -- (0,25);
\draw[fill=white] (0,0) -- (100,0) -- (100,25/4) -- (0,25/4) -- (0,0);
\draw[fill=red, opacity = 0.65] (0,0) -- (100,0) -- (100,25/4) -- (0,25/4) -- (0,0);
\draw[black, thick, -] (0,0) -- (100,0) -- (100,100) -- (0,100) -- (0,0);
\draw[black, thick, -] (50,0) -- (50,100);
\draw[black, thick, -] (12.5,0) -- (12.5,100);
\draw[black, thick, -] (100-12.5,0) -- (100-12.5,100);
\draw[black, thick, -] (0,25) -- (100,25);
\draw[black, thick, -] (12.5/4,0) -- (12.5/4,100);
\draw[black, thick, -] (100-12.5/4,0) -- (100-12.5/4,100);
\draw[black, thick, -] (0,25/4) -- (100,25/4);
\end{tikzpicture}
\end{center}
\end{minipage}
\caption{First three meshes employed in the $\h\p$ refinements for the test case~2 with exact solution~$u_2$ in~\eqref{test-case-2}.
The space--time domain is~$\QT= (0,1)\times (0,1)$.
In colours, we represent the local degrees of accuracy:
\red{\bf red:} $\p=1$; \blue{\bf blue:} $\p=2$; \darkgreen{\bf green:} $\p=3$.}
\label{figure:test-case-3-hp-meshes}
\end{figure}

Next, we focus on the test case~$3$.
For the~$h$-version of the method, we consider uniform degree of accuracy~$p = 1$,
and a sequence of uniform Cartesian space--time meshes with~$h_x = 0.5h_t = 2^{-i}, \ i = 1, \ldots, 8$.
For the~$hp$-version of the method, we proceed similarly as in~\cite[Section~7.4]{Schotzau_Schwab:2000}:
we consider a sequence of space--time meshes geometrically graded towards~$x = 0,\ x = 1$ and~$t = 0$ with grading factors~$\sigma_x = \sigma_t = 0.25$.
In Figure~\ref{figure:test-case-3-hp-meshes},
we depict the first three meshes with varying degrees of accuracy.
\begin{figure}[!ht]
    \centering
    \includegraphics[width = 0.48\textwidth]{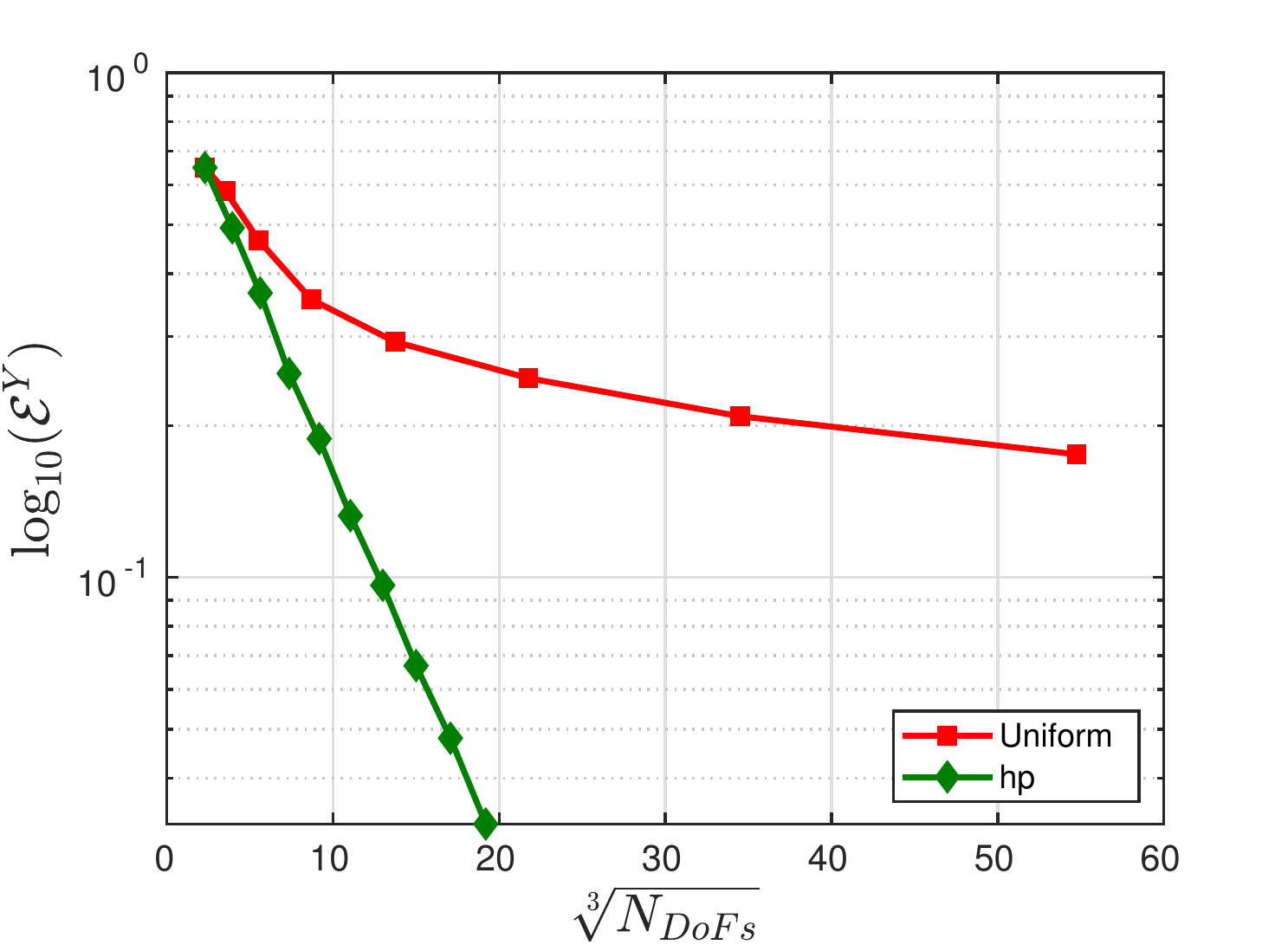}
    \hspace{0.1in}
    \includegraphics[width = 0.48\textwidth]{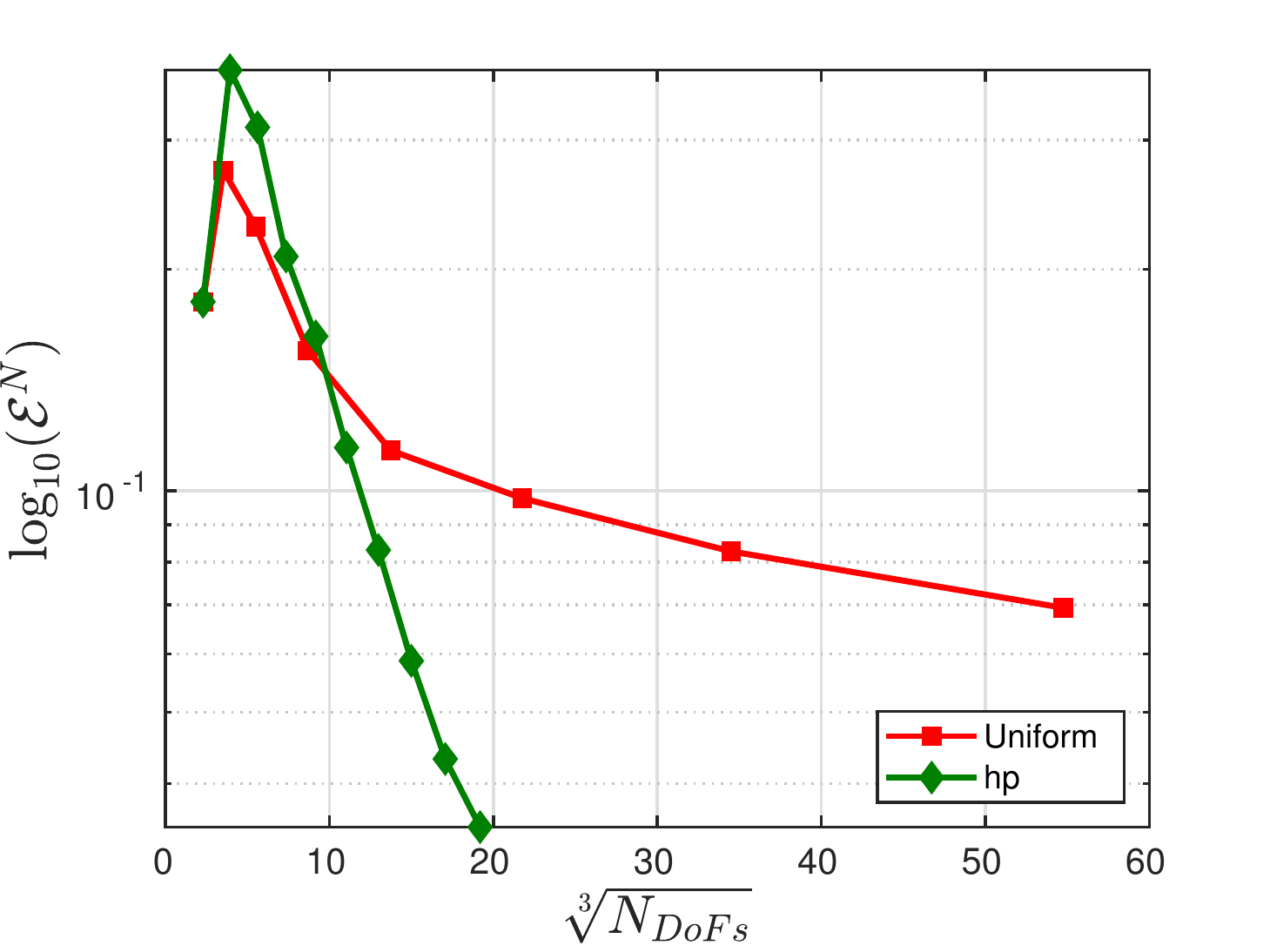}\\ 
    \includegraphics[width = 0.48\textwidth]{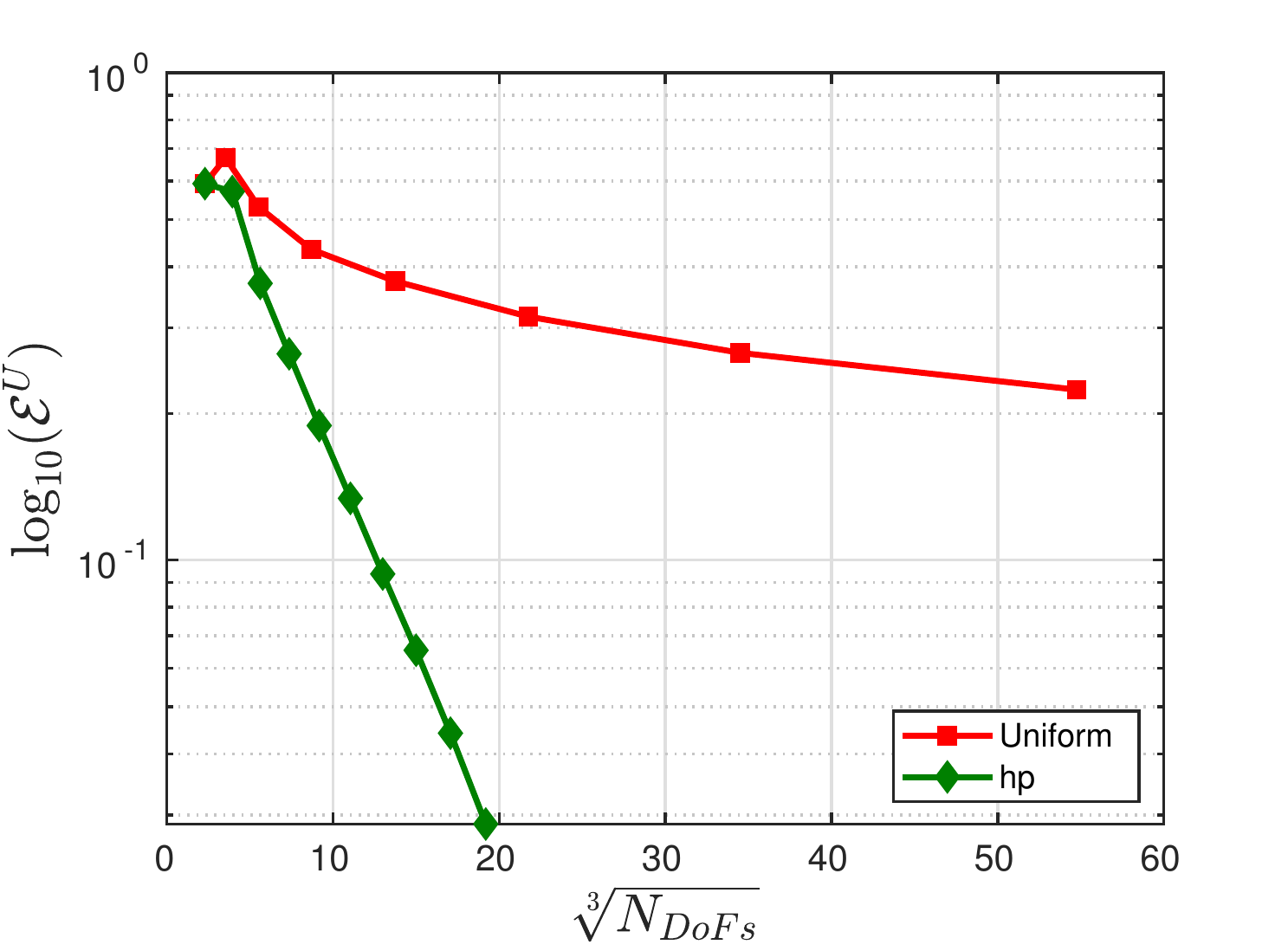} 
    \hspace{0.1in}
    \includegraphics[width = 0.48\textwidth]{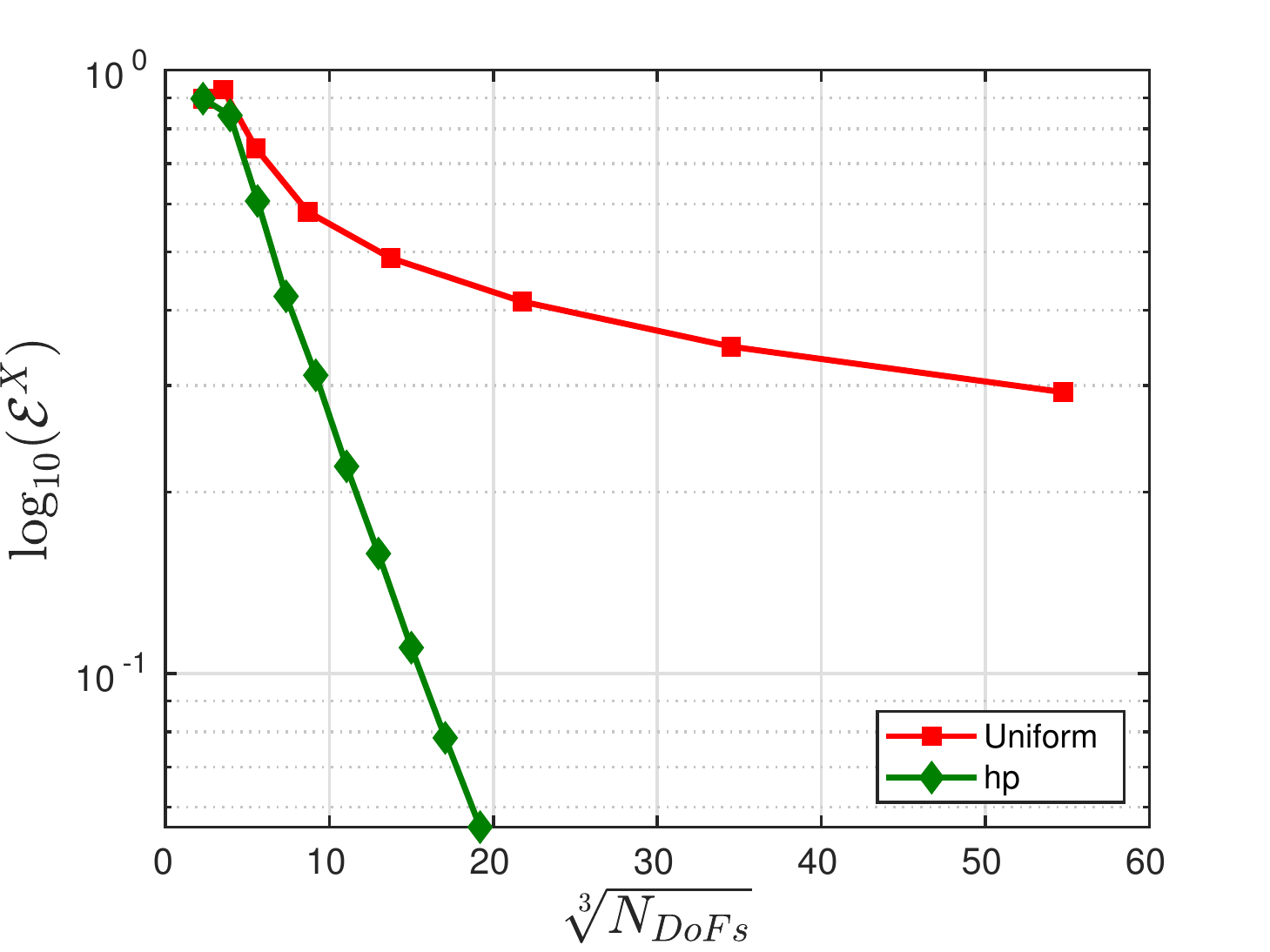}
    \caption{Convergence of the errors in~\eqref{exact-errors}
    for the~$h$- and~$hp$-versions of the method.
    We consider the test case~3 with exact solution~$u_3$ in~\eqref{test-case-3}.
    \label{FIG::hp-Incompatible}}
\end{figure}

In Figure~\ref{FIG::hp-Incompatible}, we show the errors in~\eqref{exact-errors} in \emph{semilogy} scale.
Exponential convergence in terms of the cubic root of the number of degrees of freedom
is observed for the~$hp$-version of the method;
only algebraic convergence is observed for the $\h$-version.

%%%%%%%%%%%%%%%%%%%%%%%%%%%%%%%%%%%%%%%%%%%%%%%%%%%%%%%%%%%%%%%%%%%%%%%%%%
\section{Numerical investigation: an adaptive procedure} \label{section:residual-type}
%%%%%%%%%%%%%%%%%%%%%%%%%%%%%%%%%%%%%%%%%%%%%%%%%%%%%%%%%%%%%%%%%%%%%%%%%%

We consider a standard adaptive algorithm of the form
\begin{equation} \label{adaptive:algorithm}
\textbf{SOLVE}
\quad \Longrightarrow \quad
\textbf{ESTIMATE}
\quad \Longrightarrow \quad
\textbf{MARK}
\quad \Longrightarrow \quad
\textbf{REFINE}.
\end{equation}
We base the \textbf{MARK} step on a D\"orfler marking strategy~\cite{Dorfler:1996}
with parameter~$\theta$
that we shall specify at each occurrence:
\[
\text{mark a set } \mathcal A 
\qquad
\text{ such that }
\qquad
\sum_{\E\in \mathcal A } \etaE^2
\ge \theta \eta^2.
\]
The \textbf{REFINE} step
involves only space--time mesh element refinements
as described in Section~\ref{subsection:mesh-refinement}.
For the \textbf{ESTIMATE} step,
we use the following local, computable residual-type error indicator:
given an element~$\Kfrak \in \taun$,
\[
\etaE^2:=\sum_{i=1}^5\eta_{\E,i}^2,
\]
where
\begin{equation} \label{local-error-indicator}
\begin{split}
\eta_{\E,1}^2& := \nu^{-1}\frac{\hfrakx^2}{\p^2} \Norm{
f + \nu \Deltax \PiN \uh - \cH \dpt \Pistar \uh}_{0, \Kfrak}^2,  \\
\eta_{\E,2}^2&:= \frac12 \nu^{-1}
\!\!\!\!\!\!\!\!\!\!\!\!\!
\sum_{
\scriptsize
\begin{tabular}{c}
$\F \in \FcalEt,$\\ 
$\F \not \subset \partial\Omega\times (0,T)$
\end{tabular}
}
\!\!\!\!\!\!\!\!\!\!\!\!\!\!
\frac{\hFx}{\p}\Norm{\nu \jump{\nablax \PiN \uh}}_{0,\F}^2,  \\
\eta_{\E,3}^2&:=  \nu\!\!\!\!\!\!\!\!\sum_{
\scriptsize
\begin{tabular}{c}
$\F \in \FcalEt,$\\ 
$\F  \subset \partial\Omega\times (0,T)$
\end{tabular}
}\!\!\!\!\!\!\!\!
\p \hFx^{-1}  \Norm{ \PiN \uh}_{0,\F}^2
+ \frac12 \nu \!\!\!\!\!\!\!\!\sum_{
\scriptsize
\begin{tabular}{c}
$\F \in \FcalEt,$\\ 
$\F \not \subset \partial\Omega\times (0,T)$
\end{tabular}
}\!\!\!\!\!\!\!\!
\p \hFx^{-1} \Norm{ \jump{\PiN \uh}}_{0,\F}^2, \\
\eta_{K,4}^2&:= \sum_{\Ex \in \FcalExE}
\cH^{-1}\Norm{\UWEx(\uh)}^2_{0,\Ex}, \\ 
\eta_{\E,5}^2&:=\nu \SE( (I-\PiN)\uh , (I-\PiN)\uh)  .
\end{split}
\end{equation}
The above local residual-type error indicator consists of five terms:
$\eta_{\E,1}$ is the internal residual of the projected discrete solution;
$\eta_{\E,2}$ is the boundary residual involving the normal trace of the gradient of the projected discrete solution;
$\eta_{\E,3}$ is a term due to the nonconformity in space
involving the jump of traces on time-like facets;
$\eta_{\E,4}$ is related to the upwind terms in the scheme;
$\eta_{\E,5}$ is a correction term due to the virtual element stabilization of the method.

The global error indicator and its parts read
\begin{equation} \label{global-error-indicator}
\eta^2 := \sum_{i = 1}^5 \eta_i^2,
\qquad\qquad
\eta_i^2 :=  \sum_{\Kfrak \in \taun} \eta_{\E, i}^2. 
\end{equation} 
%%%%%%%%%%%%
\subsection{Assessment of the reliability and efficiency of the error indicator}
\label{subsection:efficiency-h}
\begin{figure}[!ht]
\centering  
\includegraphics[width = 0.48\textwidth]{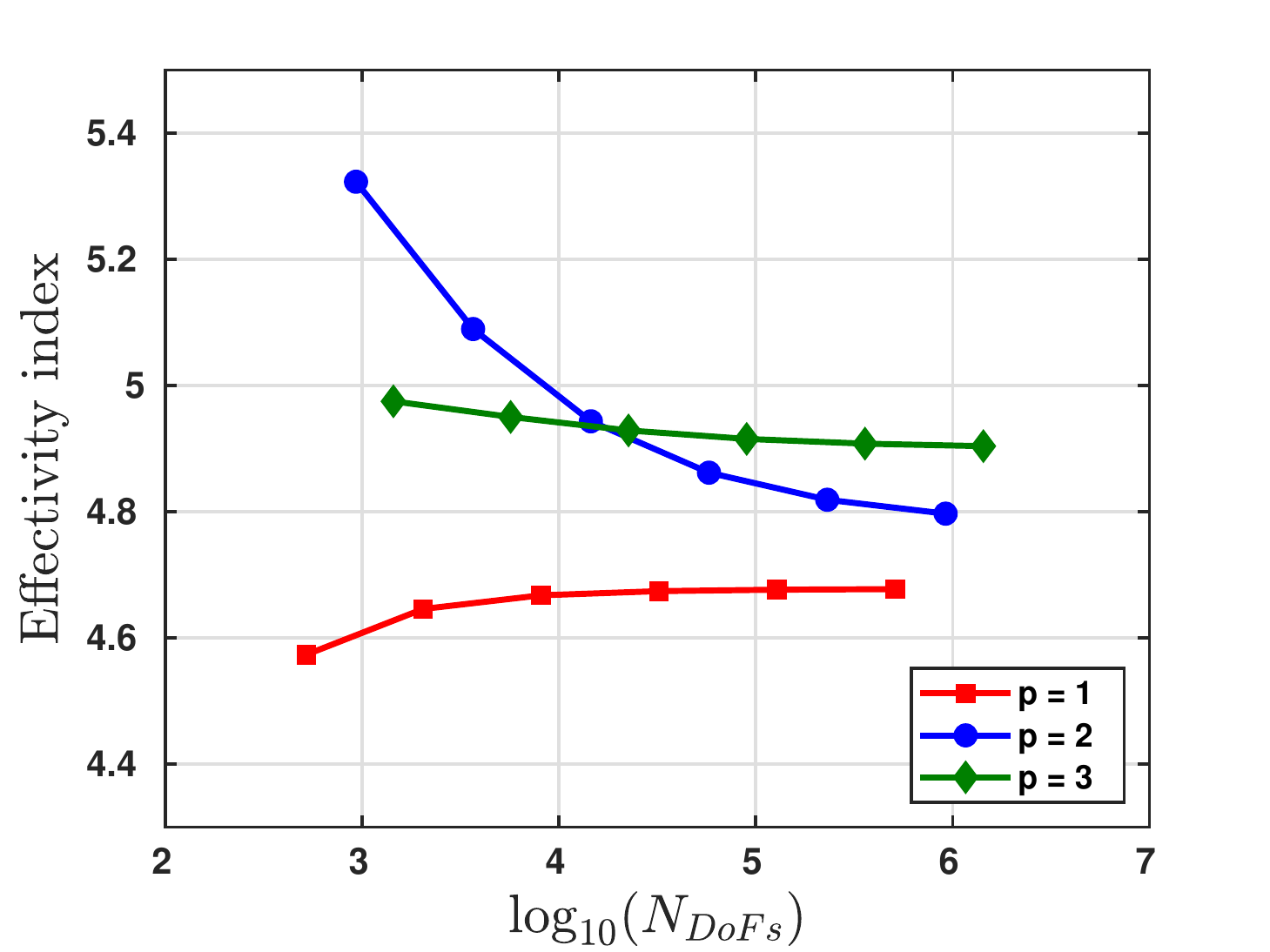} 
\hspace{0.1in}
\includegraphics[width = 0.46\textwidth]{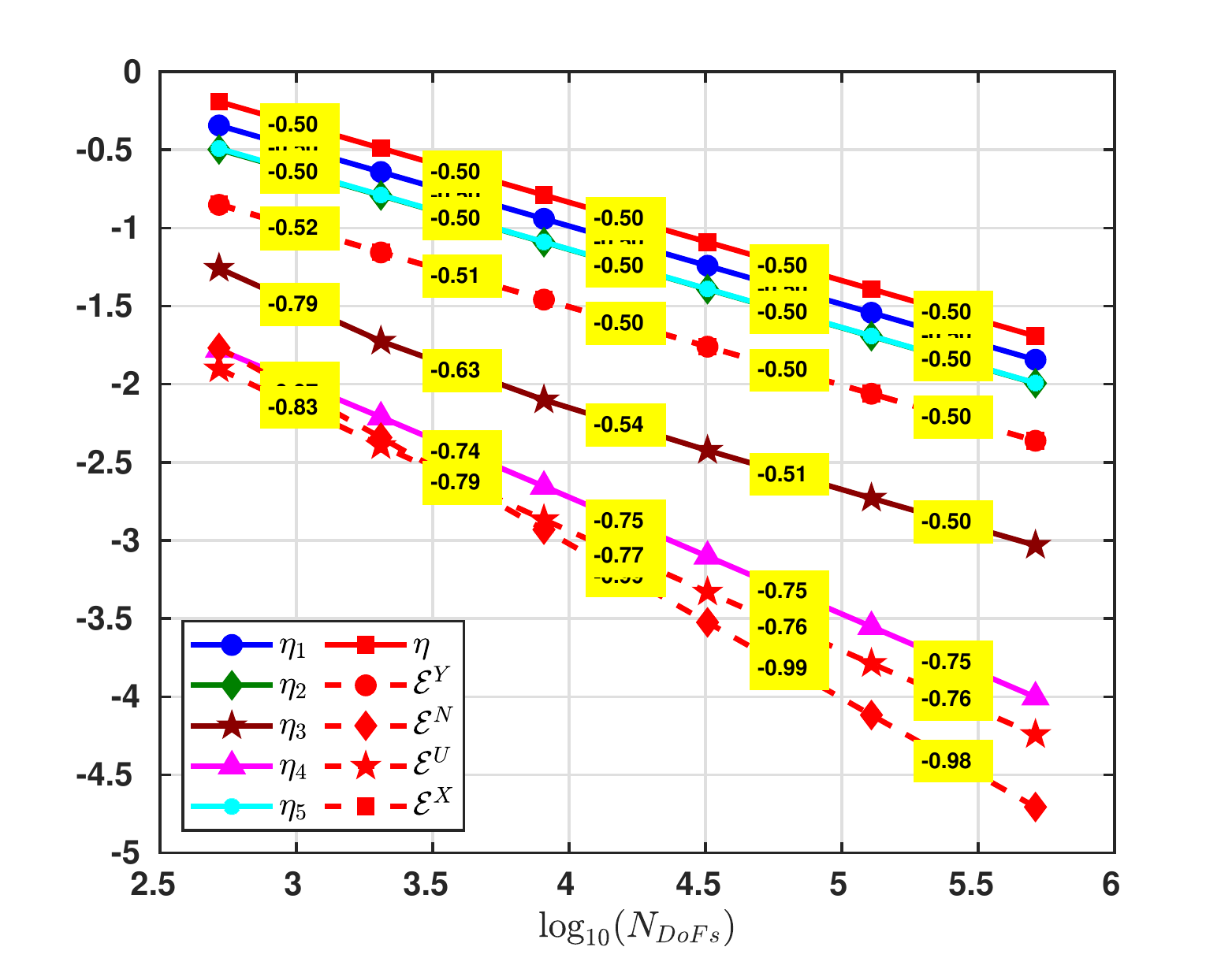}
\caption{The test case 1 with exact solution~$u_1$ in~\eqref{test-case-1}. \textbf{Left panel:} Effectivity index.
\textbf{Right panel:} Comparison of the errors in~\eqref{exact-errors} with the terms appearing in the error indicator~\eqref{global-error-indicator} for~$p = 1$.}
\label{FIG::EFFECTIVITY-SMOOTH}
\end{figure}

%%%%%%%%%%%%
We test the performance of the proposed error indicator~$\eta$ in~\eqref{global-error-indicator}.
To this aim, we introduce the effectivity index
related to the computable error~$\EcalY$:
\begin{equation} \label{effectivity-index}
\text{effectivity index}:= \frac{\eta}{\EcalY}.
\end{equation}
In Figure~\ref{FIG::EFFECTIVITY-SMOOTH} (left panel),
we assess numerically the reliability and efficiency of~$\eta$ for the test case~$1$
with exact solution~$u_1$ in~\eqref{test-case-1}
under uniform mesh refinements,
starting from a uniform space--time Cartesian mesh with~$h_x = h_t = 0.1$. The effectivity indices for~$p = 1, 2, 3$, tend to constant values,
whence the error indicator~$\eta$ appears to be efficient and reliable with respect to the error~$\EcalY$.

In Figure~\ref{FIG::EFFECTIVITY-SMOOTH} (right panel), we show all the errors in~\eqref{exact-errors} and the five terms appearing in the error indicator~\eqref{global-error-indicator} for degree of accuracy~$p = 1$.
The error indicator~$\eta$ decays with order~$\mathcal{O}(N_{DoFs}^{-\frac12})$),
i.e., slower than that for the error~$\EcalN$ in~\eqref{exact-errors},
which decays with order~$\mathcal{O}(N_{DoFs}^{-1})$.
This suggests that the error indicator~$\eta$ is not efficient with respect to the error~$\EcalN$.
%%

\begin{figure}[!ht]
\centering
\includegraphics[width = 0.46\textwidth]{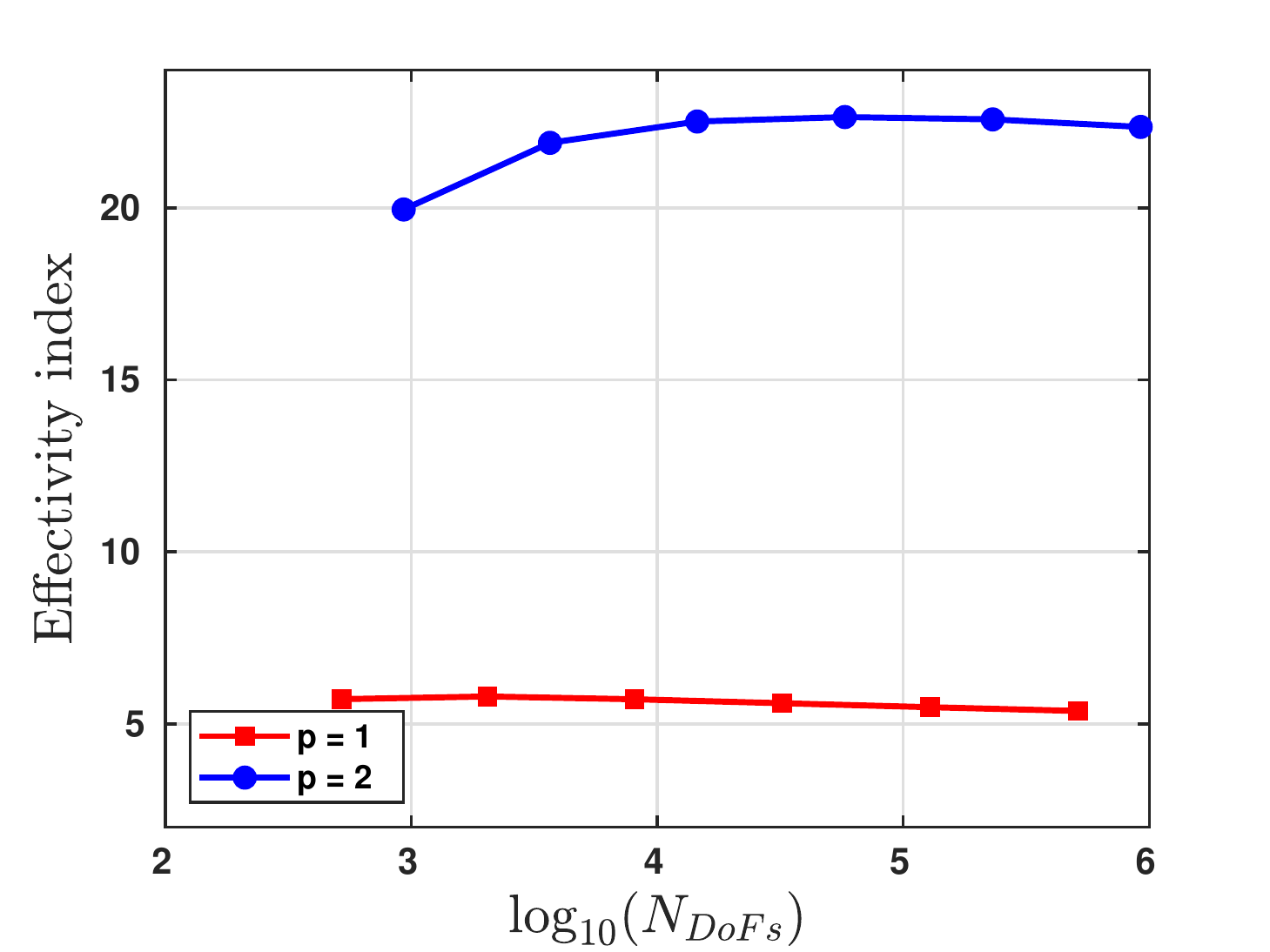}
\hspace{0.1in}
\includegraphics[width = 0.46\textwidth]{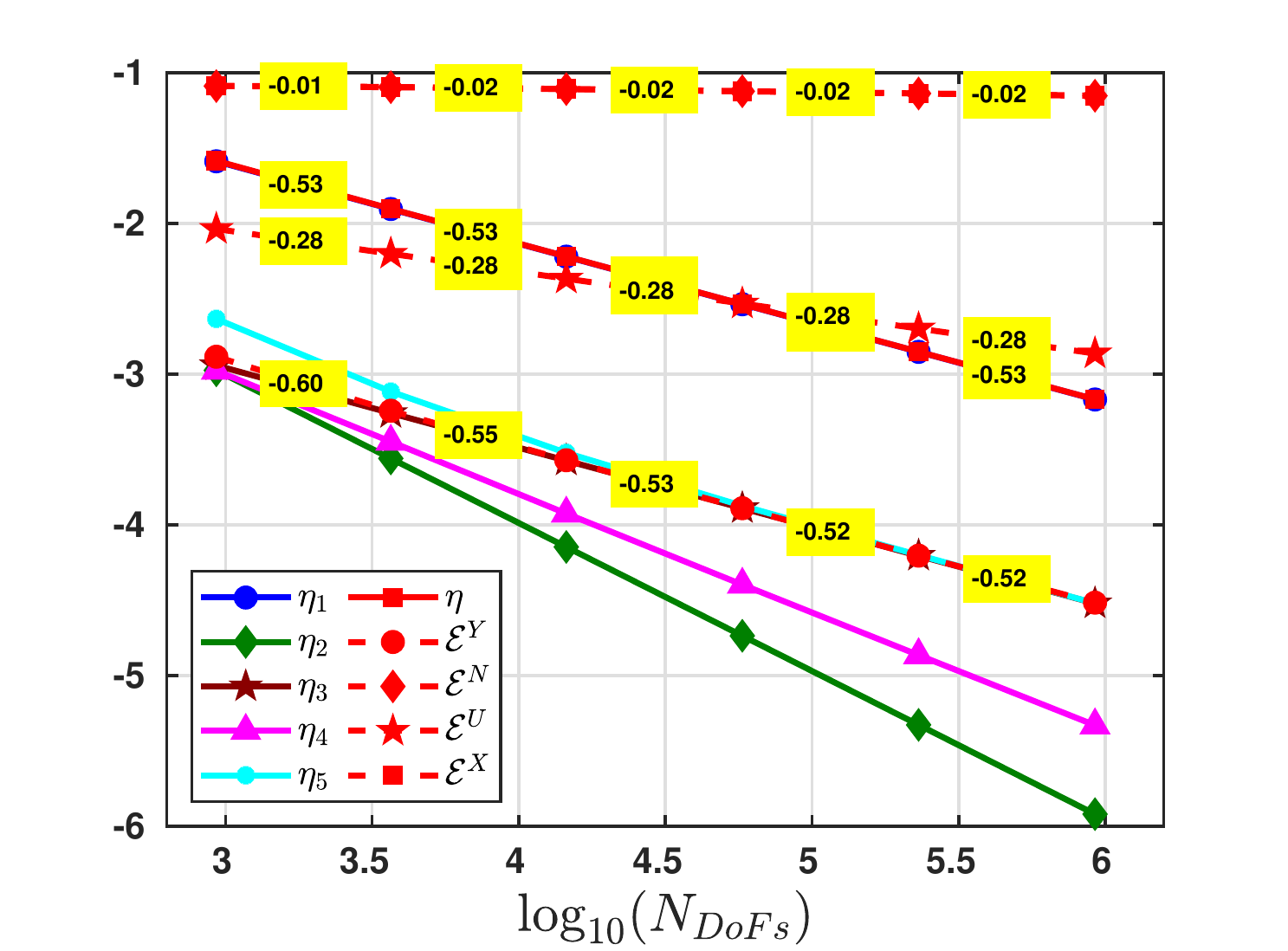}
\caption{
The test case 2 with exact solution~$u_2$ in~\eqref{test-case-2} and~$\alpha = 0.55$. 
\textbf{Left panel:} Effectivity index. \textbf{Right panel:} Comparison of the errors in~\eqref{exact-errors} with the terms appearing in the error indicator~\eqref{global-error-indicator} for~$p = 2$}.
\label{FIG::EFFECTIVITY-Singular-T0-A55}
\end{figure}

\begin{figure}[!ht]
\centering
\includegraphics[width = 0.46\textwidth]{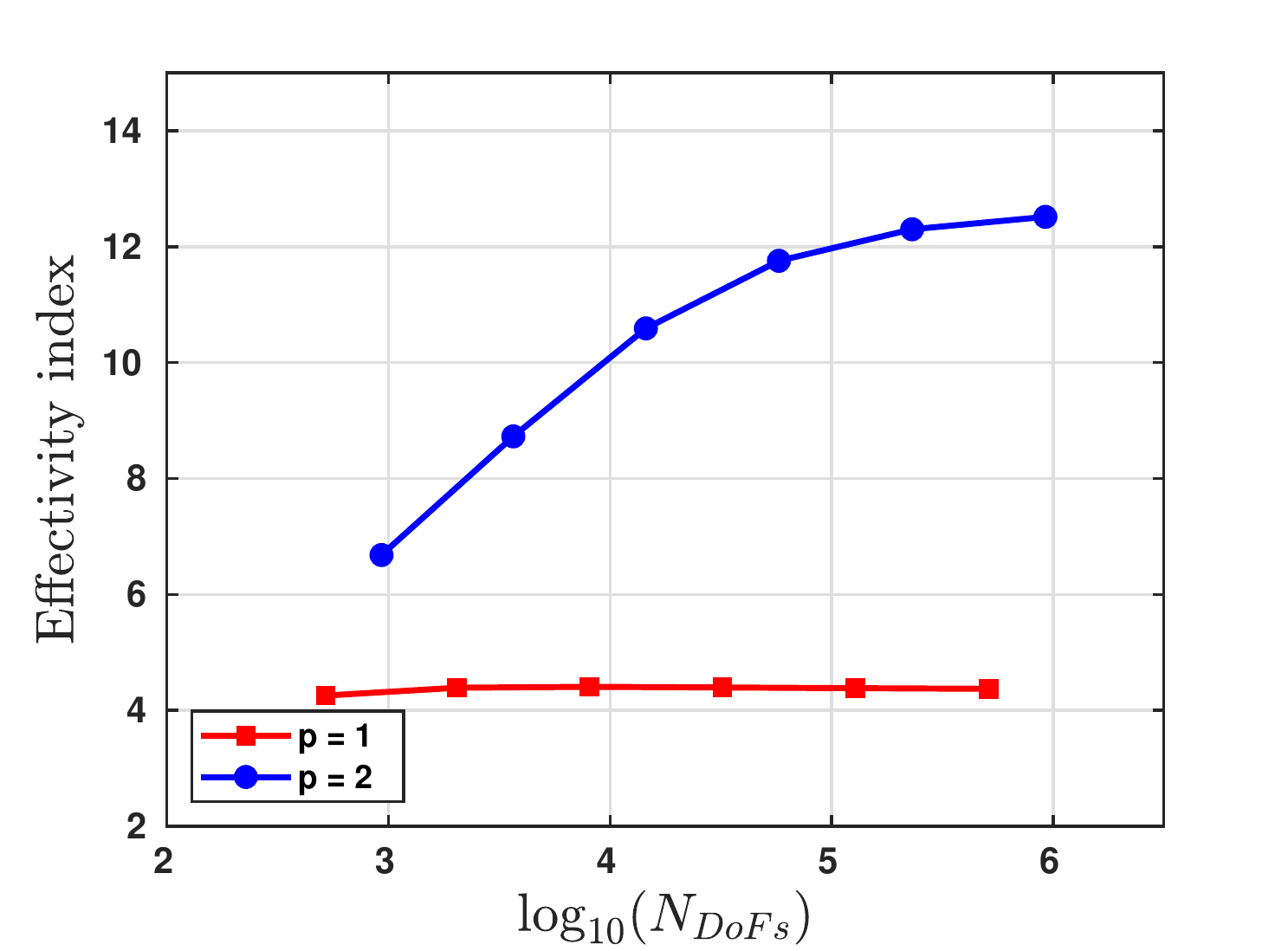}
\hspace{0.1in}
\includegraphics[width = 0.46\textwidth]{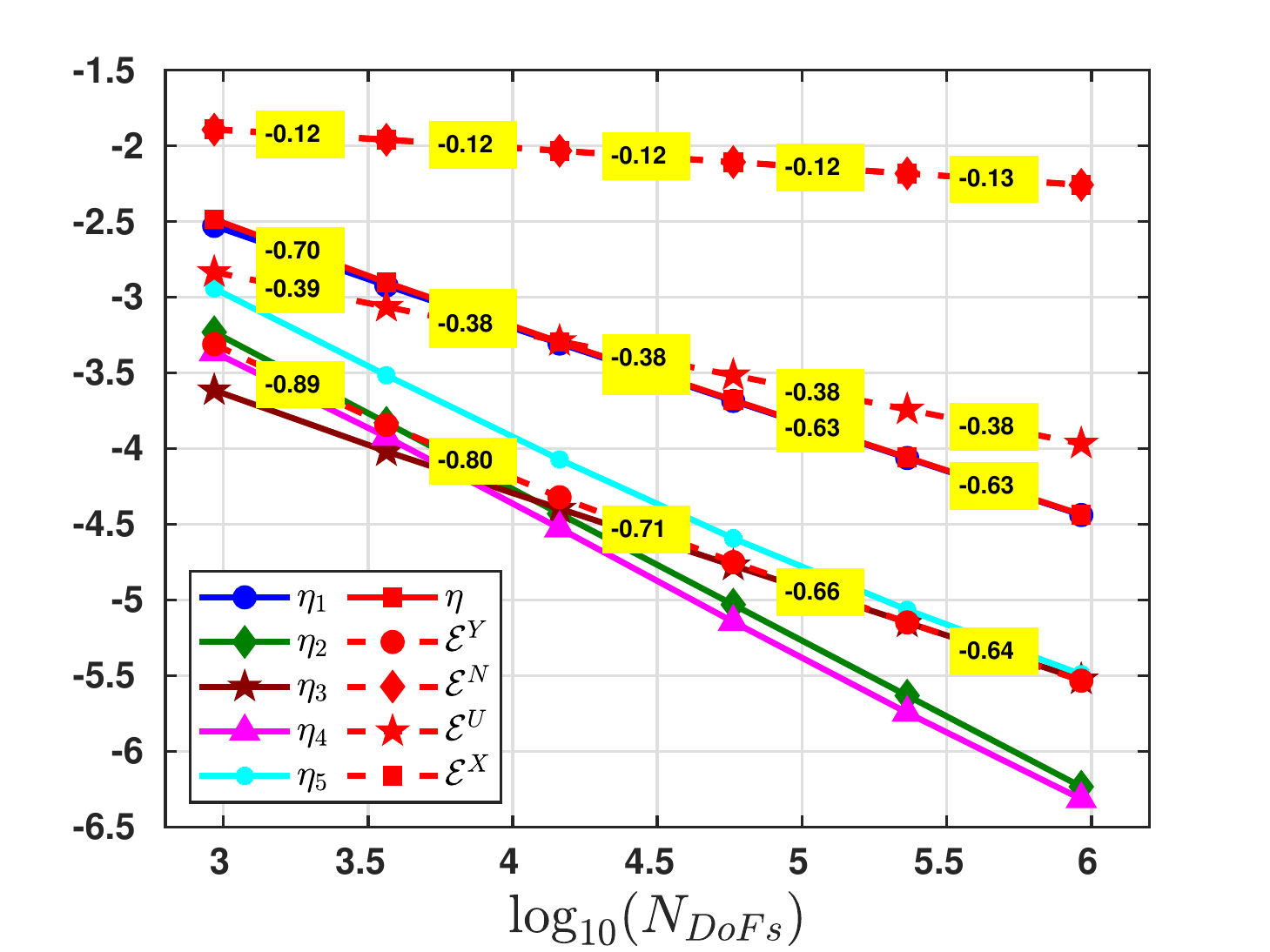}
\caption{The test case 2 with exact solution~$u_2$ in~\eqref{test-case-2} and~$\alpha = 0.75$. \textbf{Left panel:} Effectivity index. \textbf{Right panel:} Comparison of the errors in~\eqref{exact-errors} with the terms appearing in the error indicator~\eqref{global-error-indicator} for~$p = 2$}.
\label{FIG::EFFECTIVITY-Singular-T0-A75}
\end{figure}

\begin{figure}[!ht]
\centering
\includegraphics[width = 0.46\textwidth]{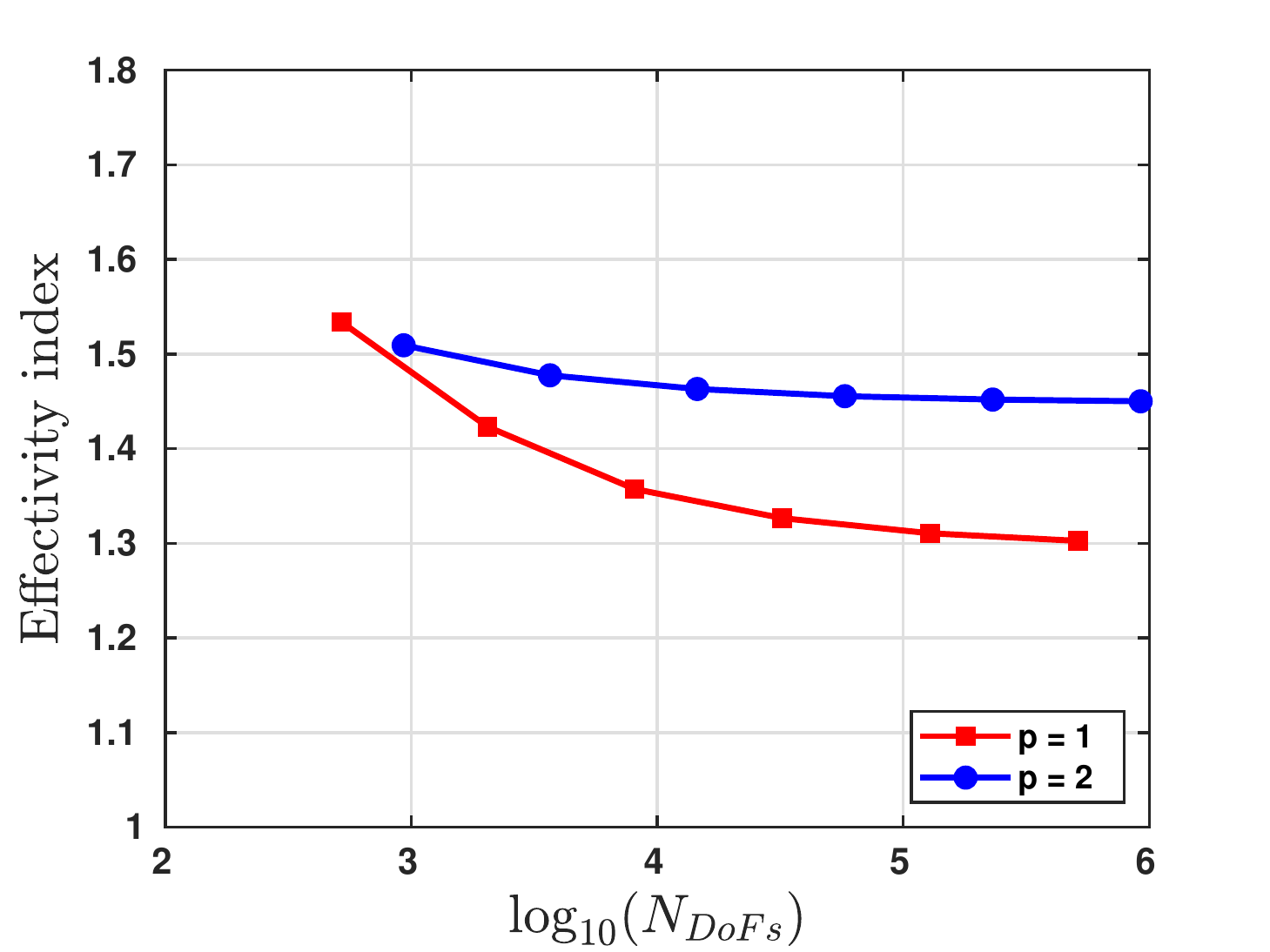}
\hspace{0.1in}
\includegraphics[width = 0.46\textwidth]{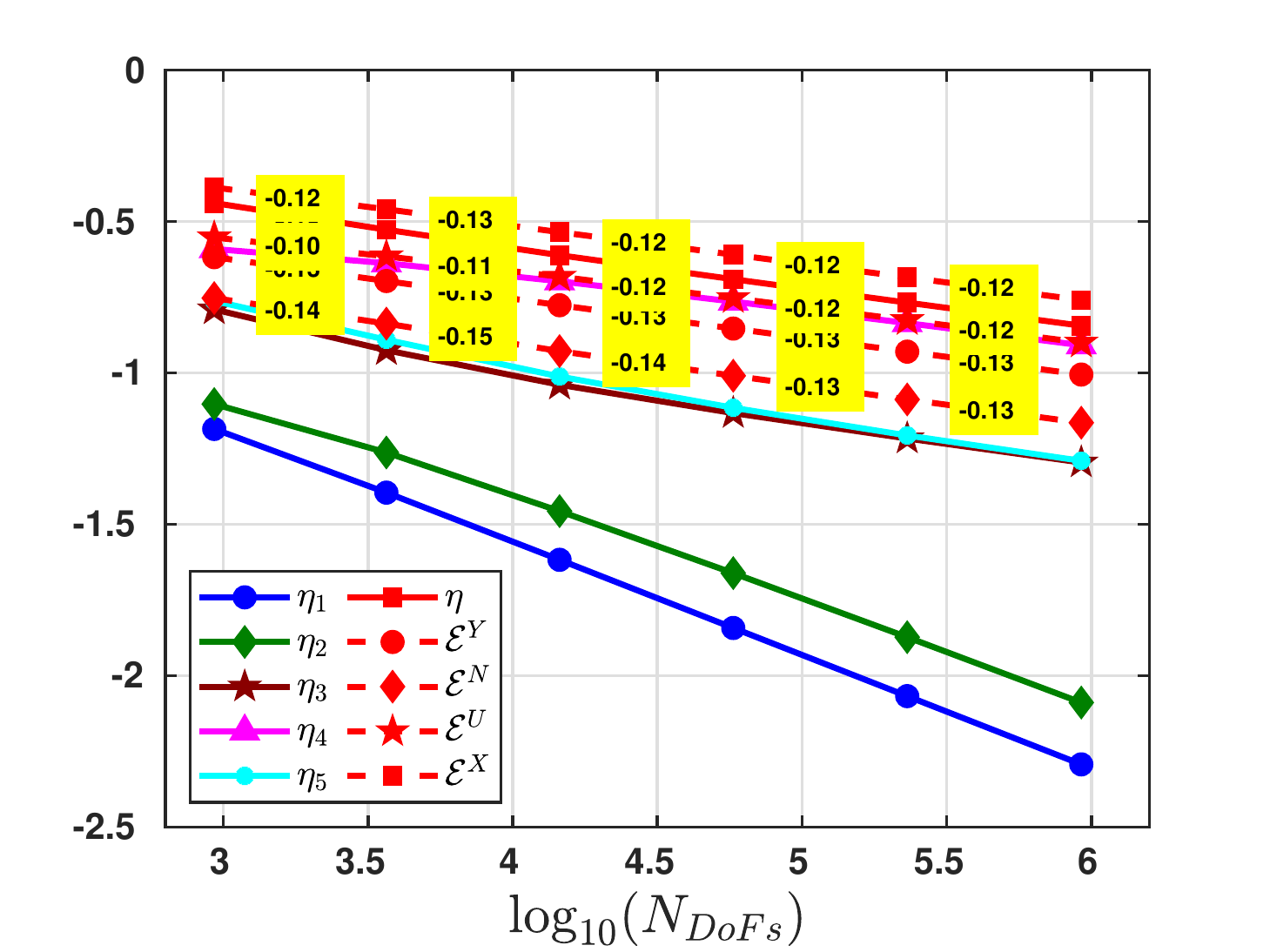}
\caption{The test case 3 with exact solution~$u_3$ in~\eqref{test-case-3}. \textbf{Left panel:} Effectivity index. \textbf{Right panel:} Comparison of the errors in~\eqref{exact-errors} with the terms appearing in the error indicator~\eqref{global-error-indicator} for~$p = 2$.}
\label{FIG::EFFECTIVITY-INCOMPATIBLE}
\end{figure}

Next, we focus on the test cases~2 and~3 with singular solutions~$u_2$ ($\alpha=0.55$ and~$0.75$) in~\eqref{test-case-2}
and~$u_3$ in~\eqref{test-case-3}
on the space--time domains~$\QT = (0,1)\times(0,1)$ and~$\QT = (0,1)\times(0,0.1)$,
respectively.
In Figures~\ref{FIG::EFFECTIVITY-Singular-T0-A55}--\ref{FIG::EFFECTIVITY-INCOMPATIBLE},
we show the decay of the errors in~\eqref{exact-errors}
and of the five terms appearing in the error indicator
under uniform mesh refinements,
starting with Cartesian meshes with~$h_x = 10 h_t = 0.1$ (for the test case 2)
and~$h_x = h_t = 0.1$ (for the test case 3).
In all cases, the effectivity indices shown on the left panels tend to constant values,
which suggests that the error indicator~\eqref{global-error-indicator} is efficient and reliable with respect to the error~$\EcalY$ also for singular solutions.

For the test case 2, the errors~$\EcalY$ and~$\EcalN$ decay with orders~$\mathcal{O}(N_{DoFs}^{-\frac12(\alpha + \frac12)})$ and~$\mathcal{O}(N_{DoFs}^{-\frac12(\alpha - \frac12)})$, respectively; see Figures~\ref{FIG::EFFECTIVITY-Singular-T0-A55} and~\ref{FIG::EFFECTIVITY-Singular-T0-A75}.
The error indicator~$\eta$ decays with the same order as that of the error~$\EcalY$,
while it it not reliable with respect to the error~$\EcalN$. 

For the test case 3, the error indicator~$\eta$ and all the errors in~\eqref{exact-errors} decay
with the same order, namely~$\mathcal{O}(N_{DoFs}^{-\frac{s}{2}})$ with~$s = \frac14$;
see Figure~\ref{FIG::EFFECTIVITY-INCOMPATIBLE}. 

In summary, the above experiments seem to indicate that the error indicator~$\eta$ is reliable and efficient for the error~$\EcalY$ but not for the error~$\EcalN$.

  
%%%%%%%%%%%%
\subsection{Adaptive mesh refinements}
\label{subsection:uni-vs-ada-h}
%%%%%%%%%%%%
We test the performance of method~\eqref{VEM} under adaptive mesh refinements
as described in~\eqref{adaptive:algorithm}.
We consider the test cases 2 and 3 with exact solutions~$u_2$ ($\alpha=0.55$) in~\eqref{test-case-2}
and~$u_3$ in~\eqref{test-case-3}, respectively.
The marking step is dictated by the error indicator in~\eqref{global-error-indicator}.

We are also interested in comparing the results
with those obtained with an adaptive procedure for the continuous finite element method (FEM) of~\eqref{continuous-weak} in~\cite{Steinbach:2015}:
\begin{equation} \label{FEM}
\begin{cases}
\text{find } \utildeh \in \Xtildeh \text{ such that}\\
b(\utildeh,\vtildeh)=(f,\vtildeh)_{0,\QT} \quad \quad \forall \vtildeh \in \Ytildeh.
\end{cases}
\end{equation}
Above, $\Ytildeh$ is a space of continuous piecewise polynomials
over a space--time simplicial tessellation of~$\QT$
and~$\Xtildeh$ is the subspace of~$\Ytildeh$
of functions with zero initial condition.
We recall the residual-type error indicator introduced in~\cite{Steinbach-Yang:2018, Steinbach-Yang:2019}:
\begin{equation} \label{FEM-error-indicator}
\begin{split}
\etaFEM^2
:= \sum_{i = 1}^2 \etaFEM_i^2, 
& \qquad  \etaFEM_i^2 = \sum_{\E \in \taun} \etaFEM_{\E, i}^2, \\
\etaFEM_{\E, 1}^2
:= \frac{\hE^2}{\p^2} \Norm{f + \nu \Deltax \utildeh - \cH \dpt \utildeh}_{0, \E}^2, & \qquad 
\etaFEM_{\E, 2}^2 : = \frac12 \sum_{\F \in \FcalEt}
\frac{\hF}{\p}\Norm{\nu \jump{\nablax \utildeh}}_{0,\F}^2,
\end{split}
\end{equation}
and the error quantity
\begin{equation}
\EY = \Norm{u - \utildeh}_{Y}.
\end{equation}

For the VEM, we start with a mesh with~$1$ element;
for the continuous FEM, we start with a structured simplicial mesh with~$2$ elements.

\begin{figure}[!ht] 
\centering
\includegraphics[width = 0.32\textwidth, height = 1.57in, clip, trim = 0 0 15 10]{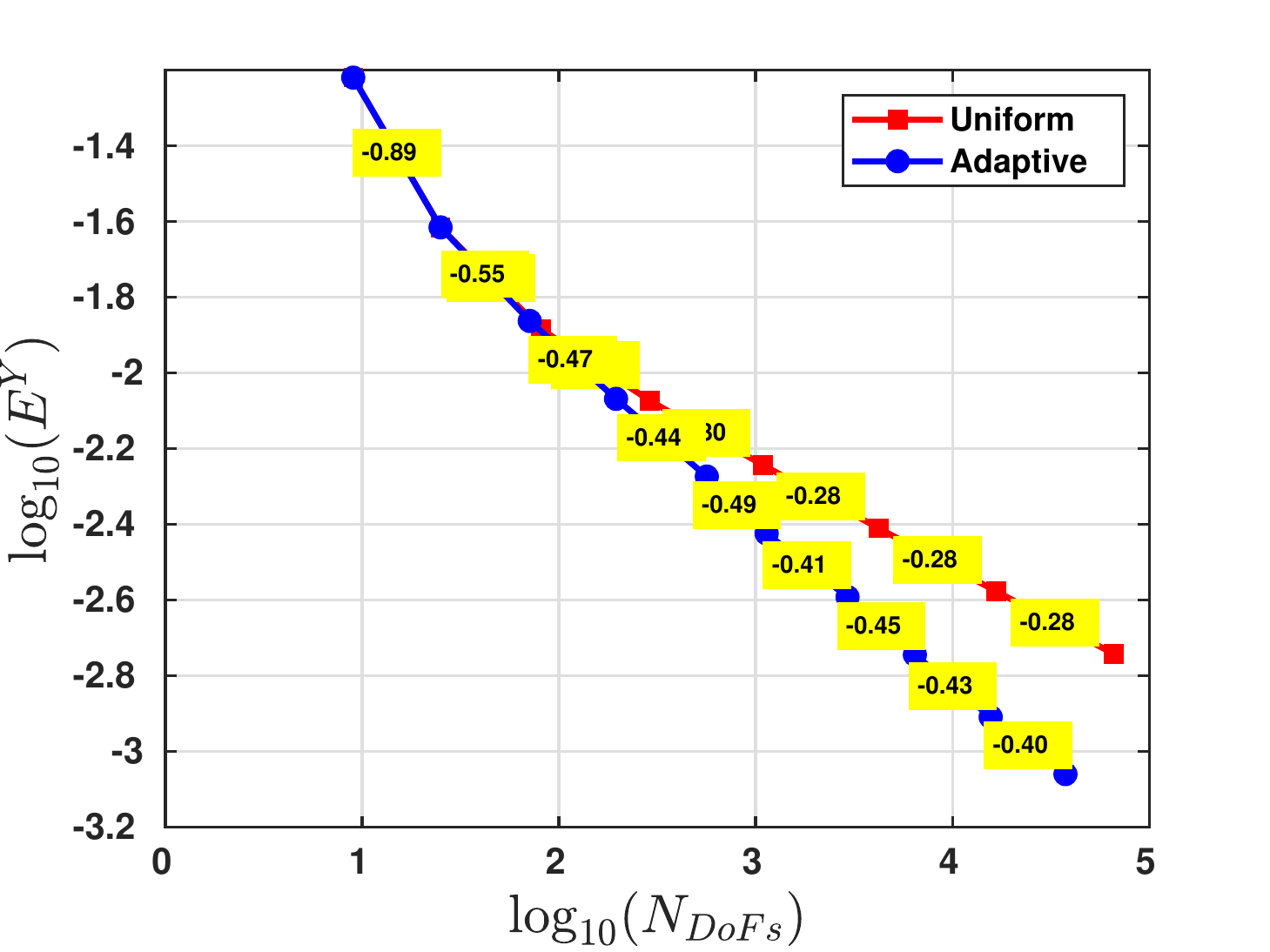} 
\includegraphics[width = 0.32\textwidth, height = 1.57in, clip, trim = 0 0 15 10]{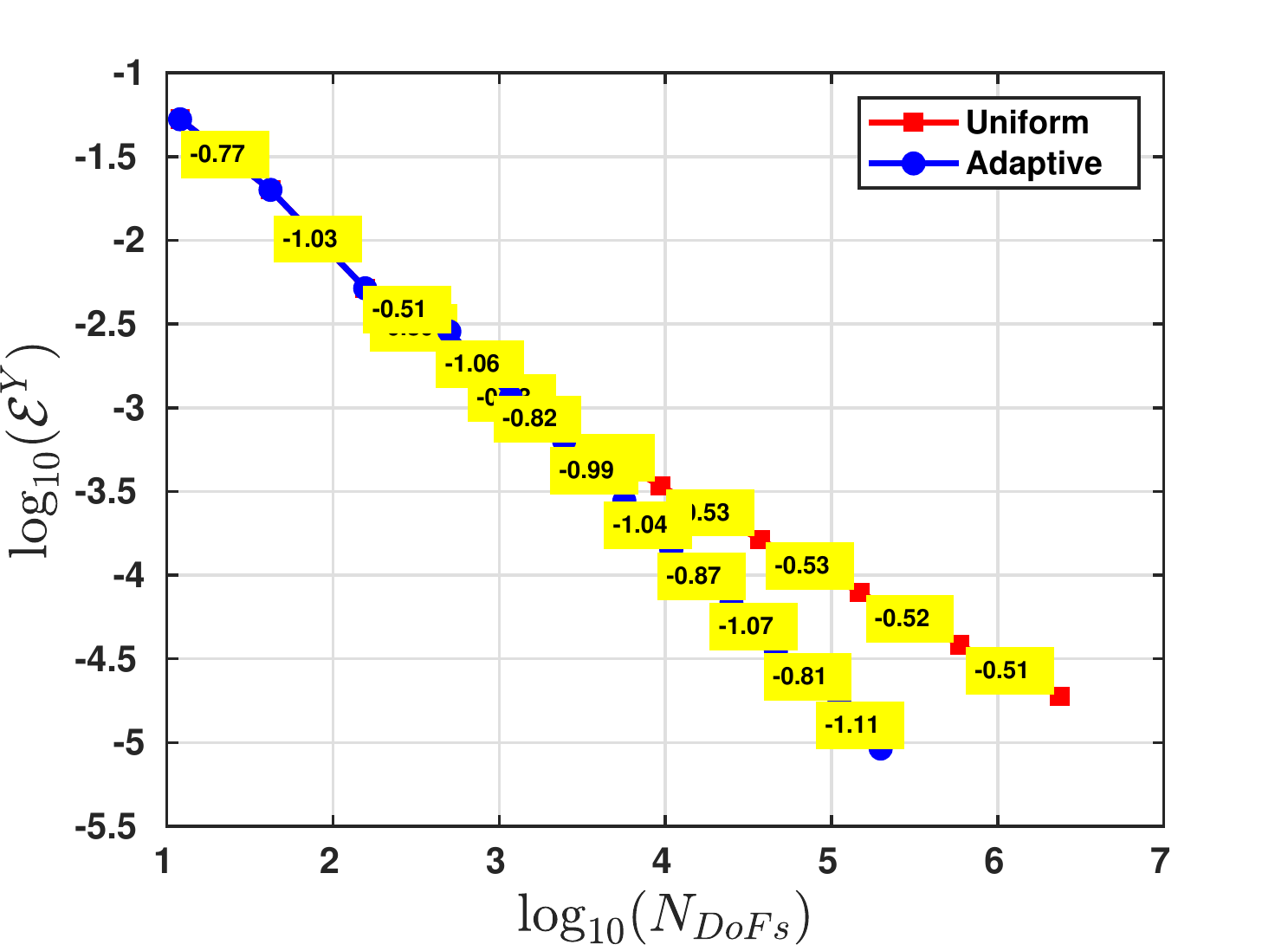}
\includegraphics[width = 0.32\textwidth, height = 1.57in, clip, trim = 0 0 15 10]{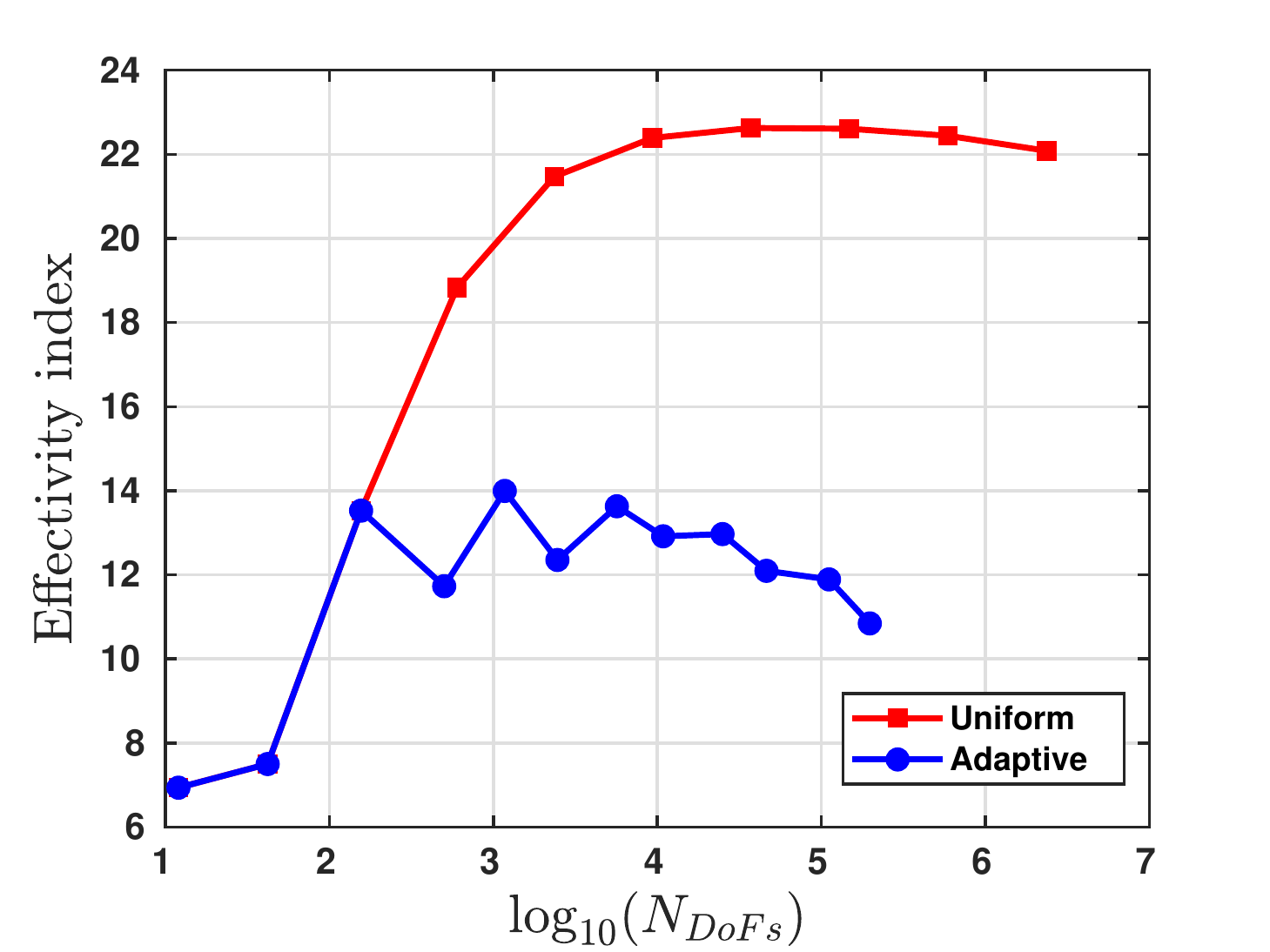} 
\caption{The test case 2 with exact solution~$u_2$ in~\eqref{test-case-2} and~$\alpha = 0.55$.
\textbf{Left panel:} Error~$E^Y$ for the continuous FEM with~$p = 2$ under uniform and adaptive mesh refinements, with convergence rates of approximately~$\mathcal{O}(N_{DoFs}^{-0.28})$ and $\mathcal{O}(N_{DoFs}^{-0.40})$,   respectively.
\textbf{Central panel:} Error~$\EcalY$ for the VEM with~$p = 2$ under uniform and adaptive mesh refinements, with convergence rates of approximately~$\mathcal{O}(N_{DoFs}^{-0.52})$ and $\mathcal{O}(N_{DoFs}^{-1})$,   respectively. \textbf{Right panel:} The effectivity index for the adaptive VEM.}
\label{FIG::h-ADAPTIVE-TEST-CASE-2}
\end{figure}

In Figure~\ref{FIG::h-ADAPTIVE-TEST-CASE-2}, we show the errors of both methods under uniform and adaptive mesh refinements for the test case~$2$ with~$\alpha = 0.55$.
Uniform and adaptive mesh refinements for the VEM in~\eqref{VEM}
lead to higher convergence rates than those for the continuous FEM in~\eqref{FEM}.
For the D\"orfler marking strategy, we set~$\theta = 0.99$ for the VEM and~$\theta = 0.9$ for the FEM.

\begin{figure}[!ht]
\centering
\includegraphics[width = 0.32\textwidth, height =1.57in, clip, trim = 0 0 15 10]{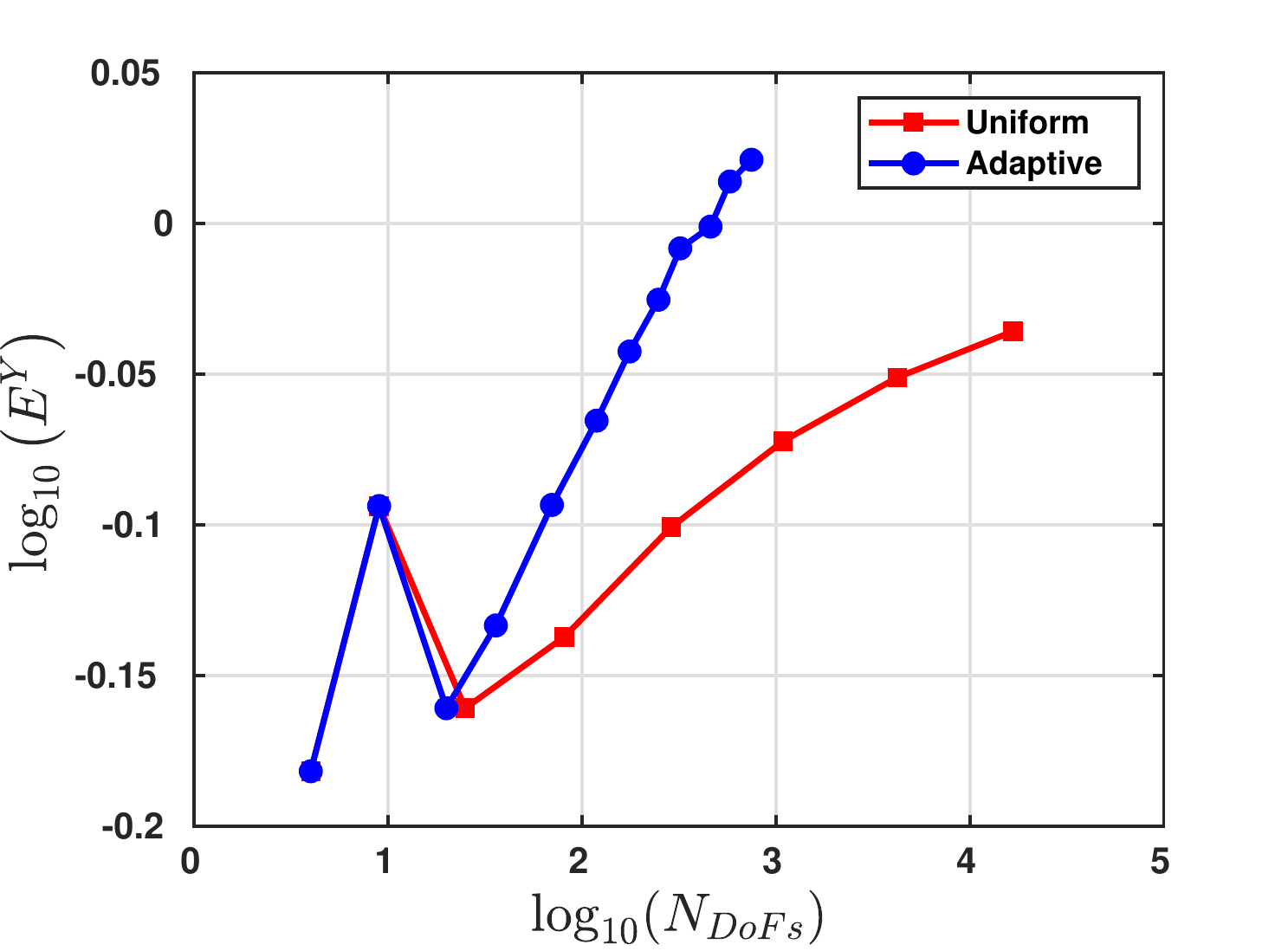}
\includegraphics[width = 0.33\textwidth, height = 1.57in, clip, trim = 0 0 15 10]  {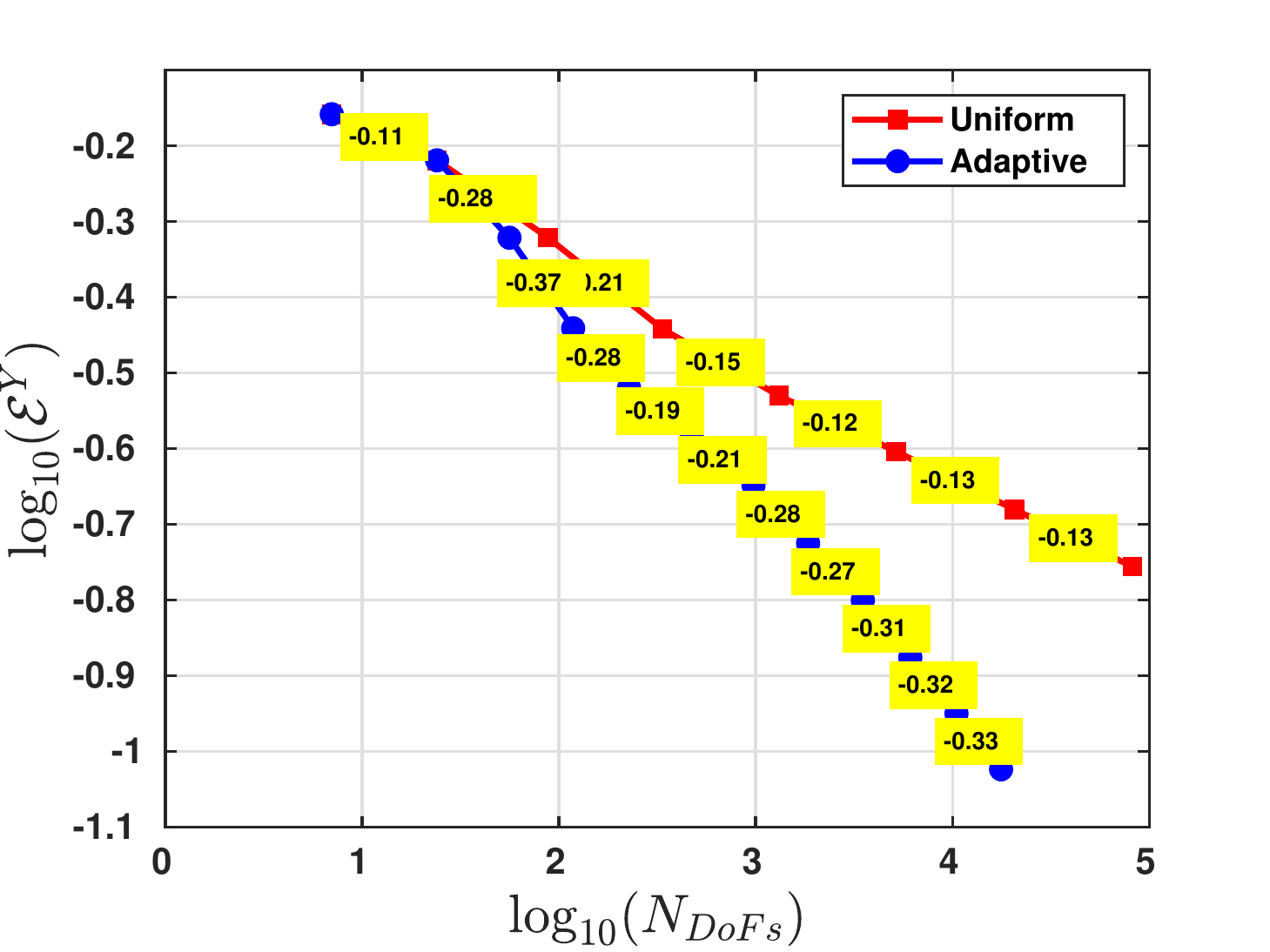}
\includegraphics[width = 0.32\textwidth, height =1.57in, clip, trim = 0 0 15 10] {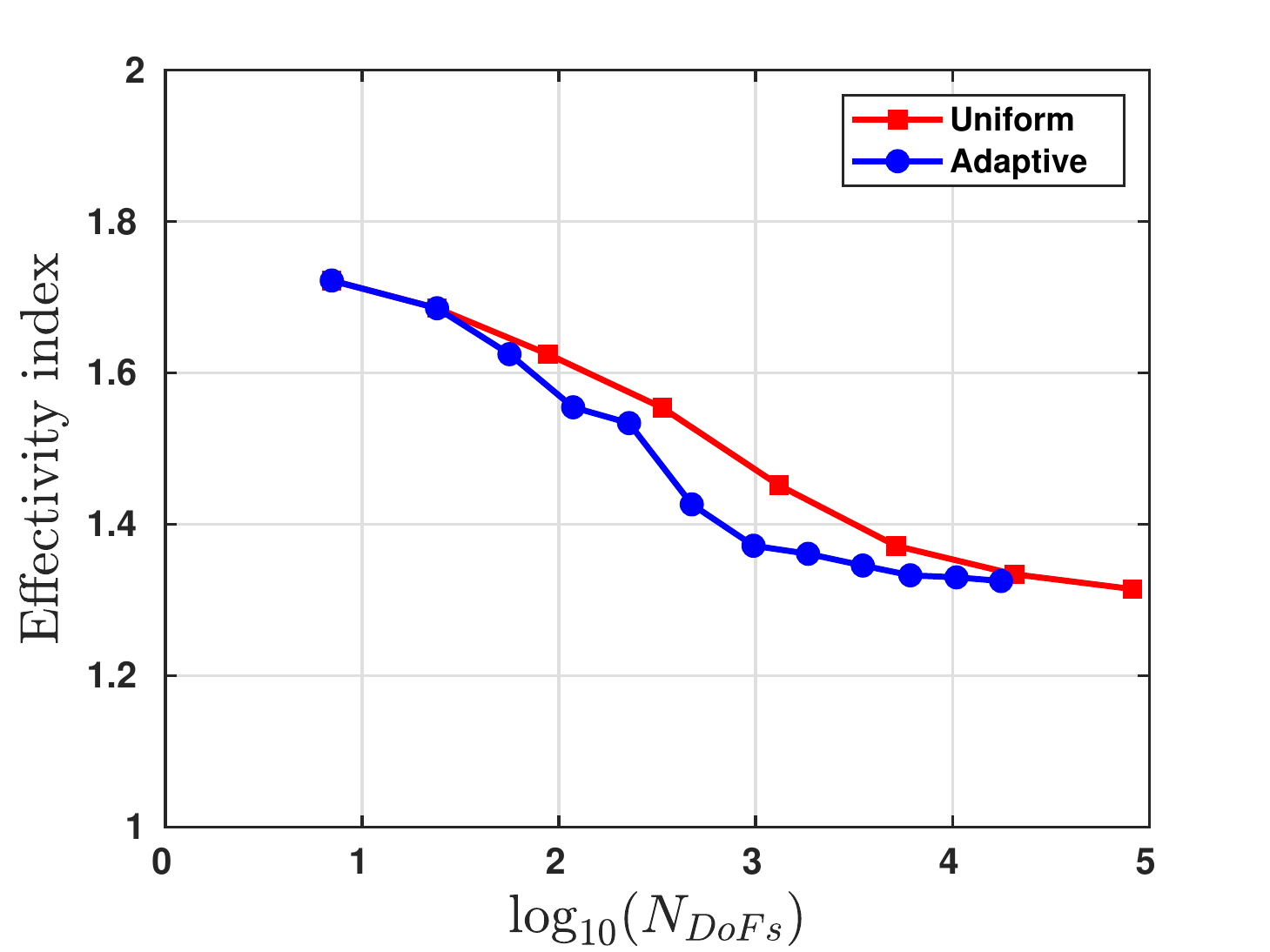}
\caption{The test case 3 with exact solution~$u_3$ in~\eqref{test-case-3}.
\textbf{Left panel:} Error~$E^Y$ for the continuous finite element with~$p = 1$ under uniform and adaptive refinements.
\textbf{Central panel:} Error~$\EcalY$ for the VEM with~$p = 1$ under uniform and adaptive mesh refinements, with convergence rates of approximately~$\mathcal{O}(N_{DoFs}^{-0.13})$ and $\mathcal{O}(N_{DoFs}^{-33})$, respectively. \textbf{Right panel:} The effectivity index for the adaptive VEM.
\label{FIG::h-ADAPTIVE-TEST-CASE-3}}
\end{figure}

In Figure \ref{FIG::h-ADAPTIVE-TEST-CASE-3},
we show the~$Y$-type errors for both methods under uniform and adaptive mesh refinements with D\"orfler marking parameter~$\theta = 0.9$ for the test case~$3$ and~$p=1$.
The adaptive procedure for the VEM~\eqref{VEM} leads to higher convergence rates compared to those obtained for uniform refinements.
For the continuous FEM~\eqref{FEM}, although the adaptive procedure produces meshes that are refined towards the bottom corners
as in Figure~\ref{FIG::SEQUENCE-MESHES-INCOMPATIBLE} (right panels),
the error~$E^Y$ does not converge to zero.
This experiment seems to suggest that the presented nonconforming approach is able to capture correctly singularities arising from incompatible initial and boundary data.

%

\begin{figure}[!ht]
\centering
\subfloat[VEM (mesh 7)]{
\includegraphics[width = 0.4\textwidth]{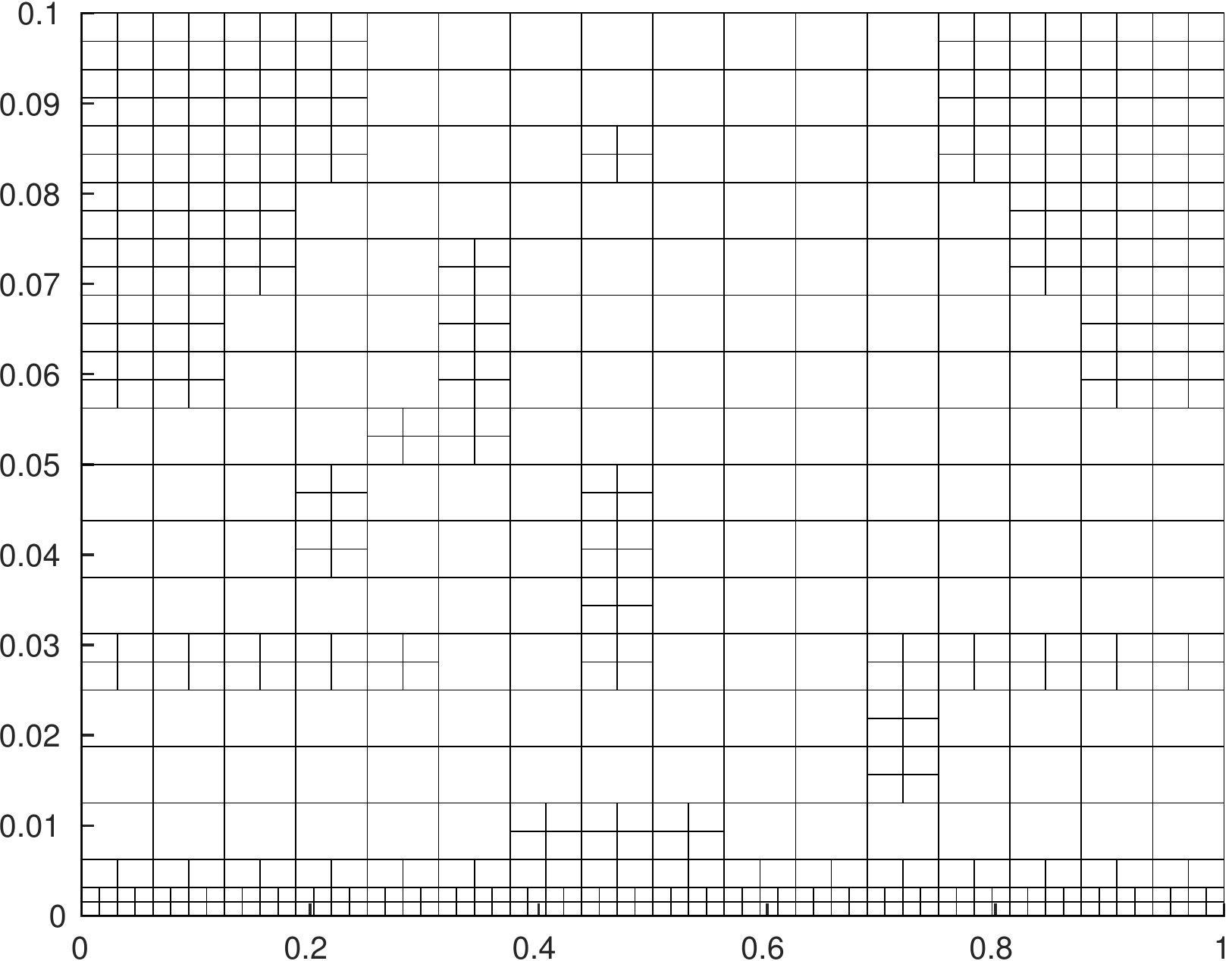}
}
\hspace{0.1in}
\subfloat[FEM (mesh 7)]{
\includegraphics[width = 0.4\textwidth]{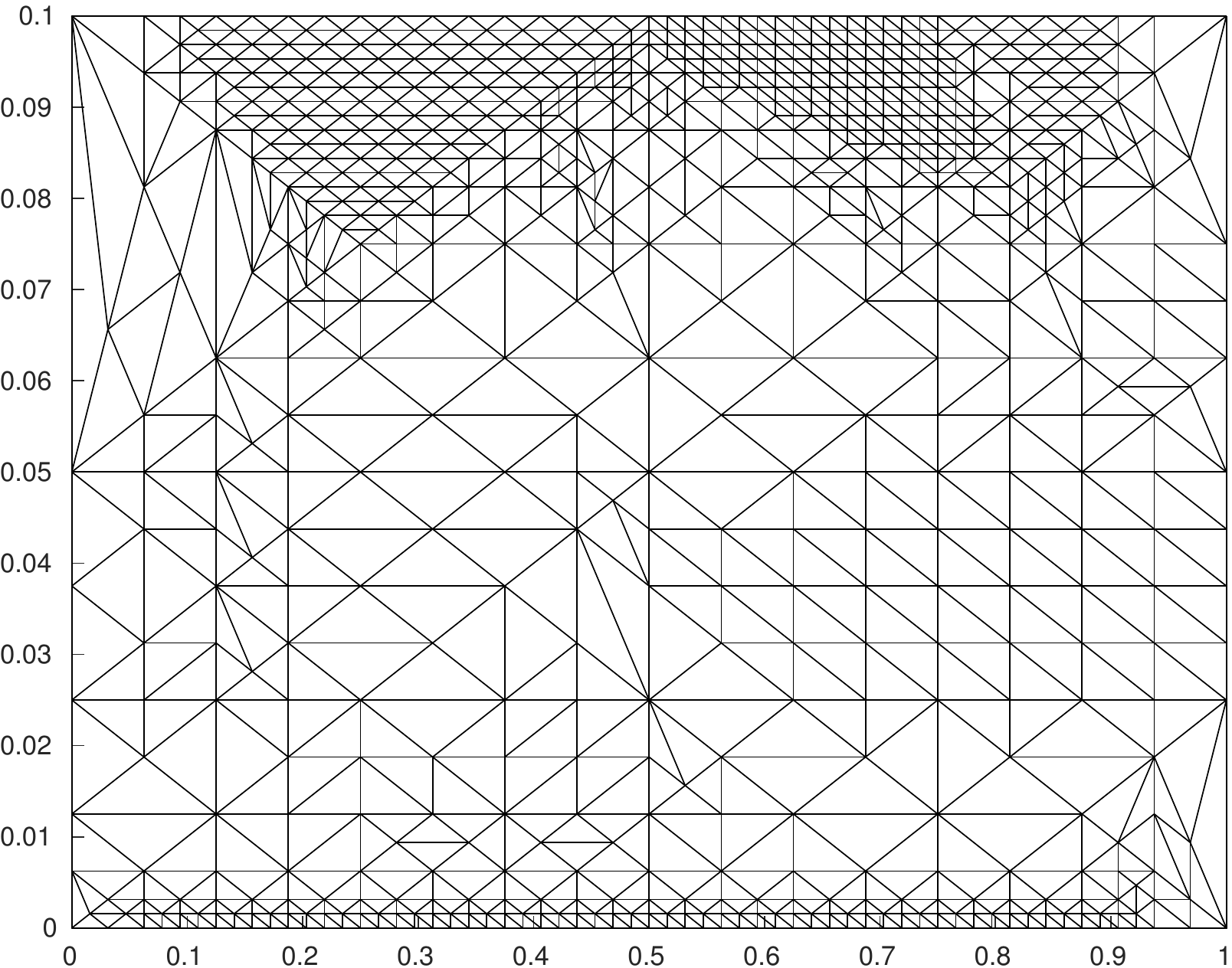}
}\\
\subfloat[VEM (mesh 10)]{ 
\includegraphics[width = 0.4\textwidth]{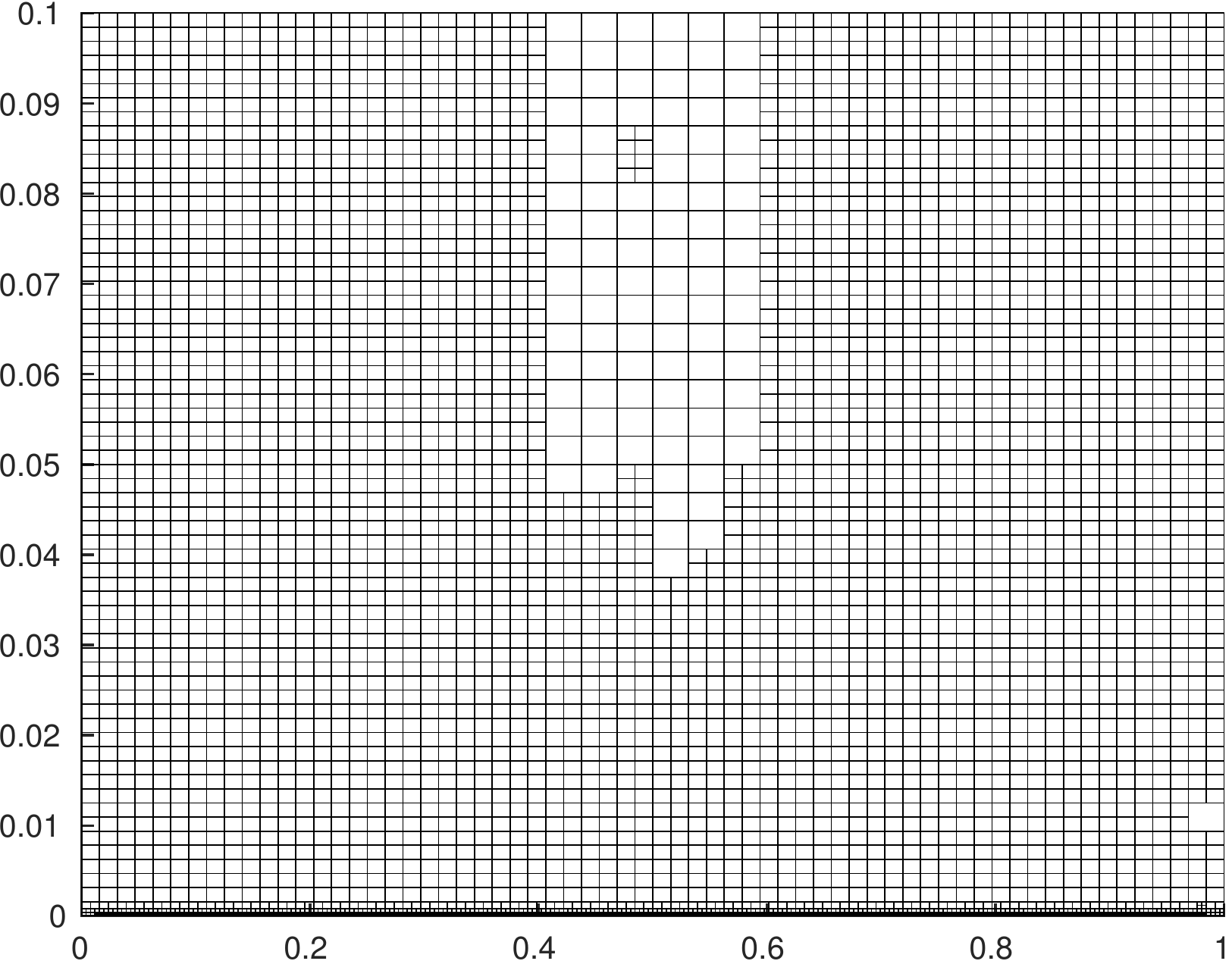} 
}
\hspace{0.1in}
\subfloat[FEM (mesh 10)]{
\includegraphics[width = 0.4\textwidth]{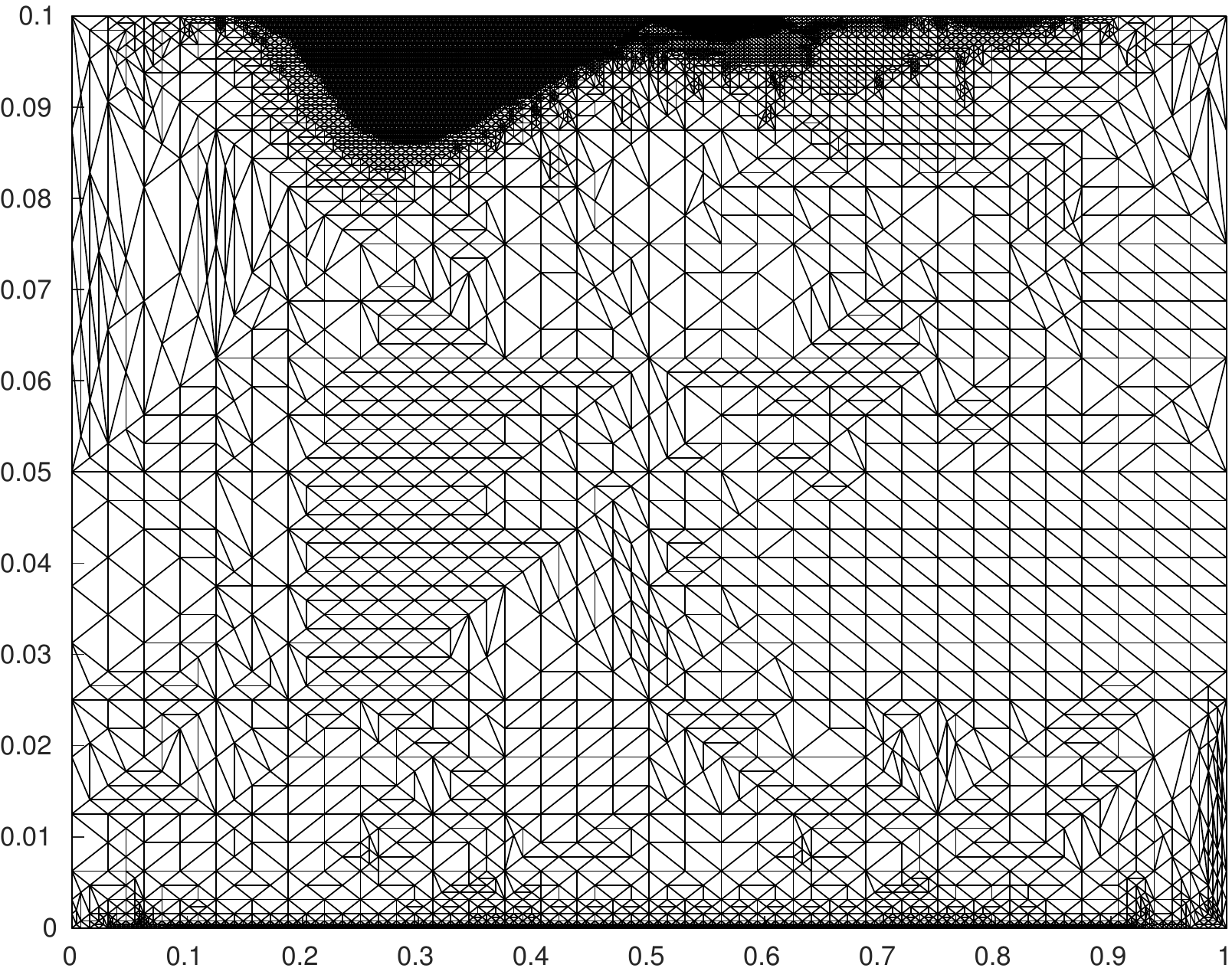}
}
\caption{Meshes generated by the adaptive schemes
driven by the VEM error indicator~$\eta$ in~\eqref{global-error-indicator} (left panel) and the continuous FEM error indicator~$\etaFEM$ (right panel)
for the test case with exact solutions~$u_2$ ($\alpha = 0.55$) in~\eqref{test-case-2}.}
\label{FIG::SEQUENCE-MESHES-t-alpha} 
\end{figure}

\begin{figure}[!ht]
\centering
\subfloat[VEM (mesh 7)]{
\includegraphics[width = 0.38\textwidth]{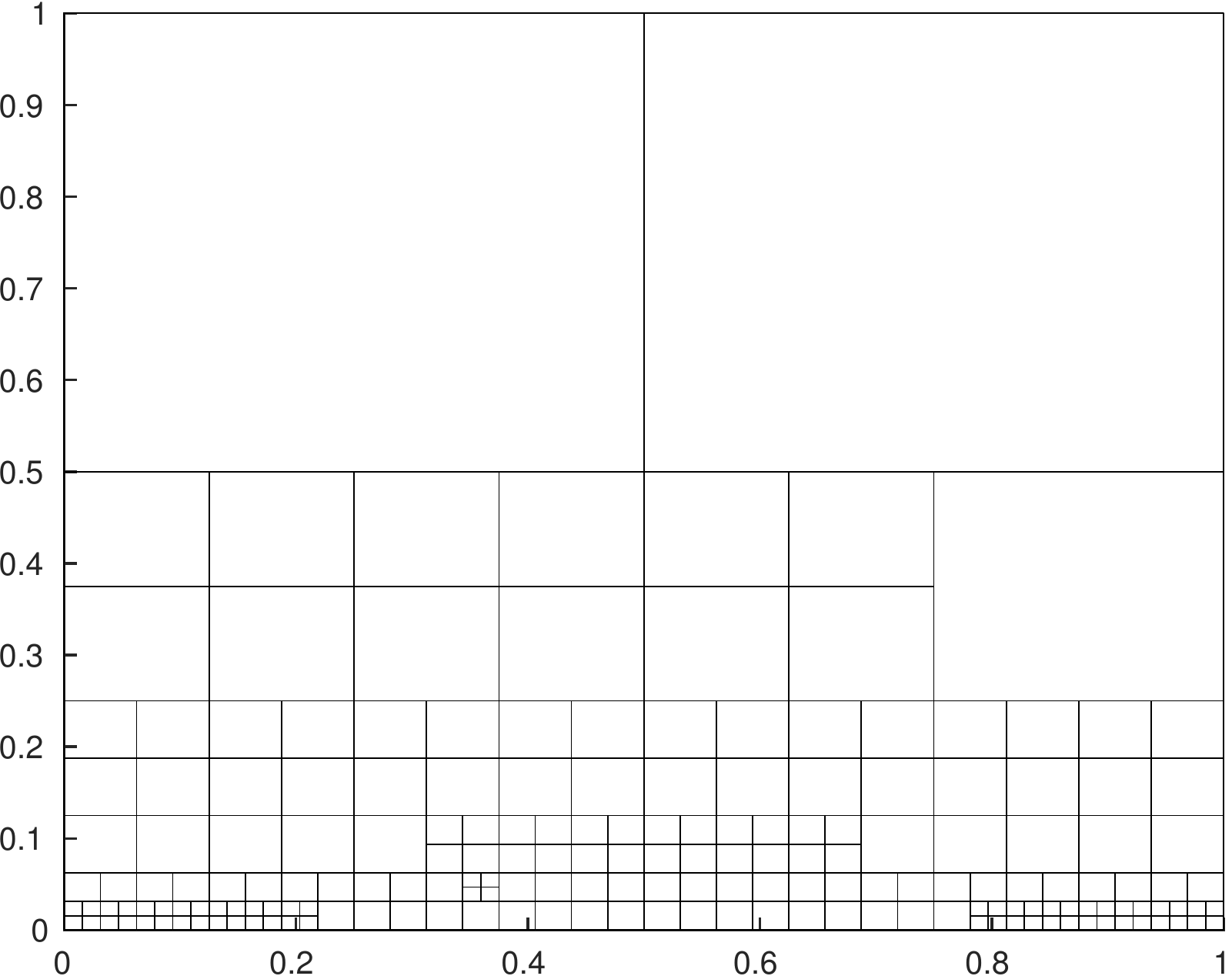}
}
\hspace{0.1in}
\subfloat[FEM (mesh 7)]{
\includegraphics[width = 0.38\textwidth]{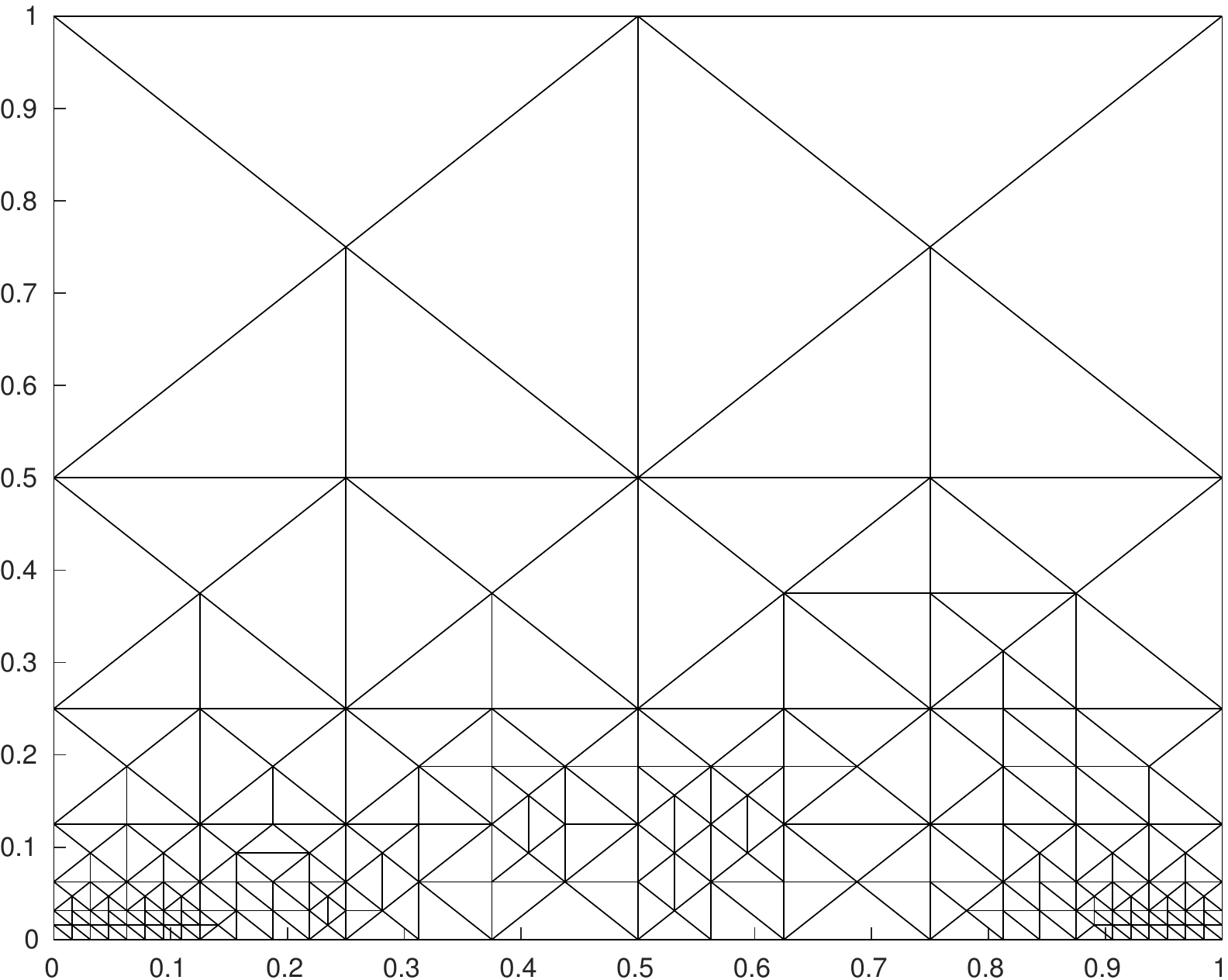}
}\\
\subfloat[VEM (mesh 12)]{
\includegraphics[width = 0.38\textwidth]{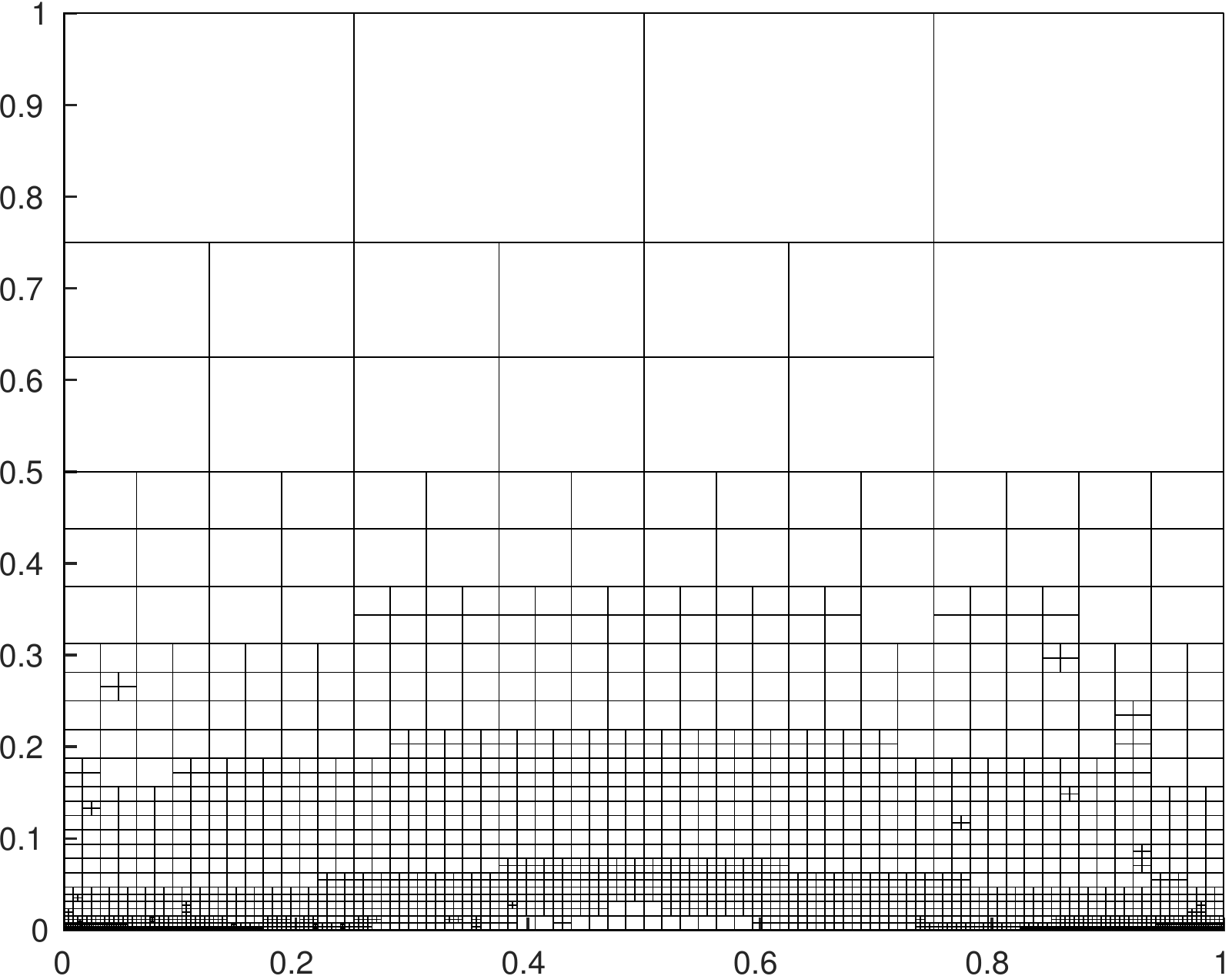}
}
\hspace{0.1in}
\subfloat[FEM (mesh 12)]{
\includegraphics[width = 0.38\textwidth]{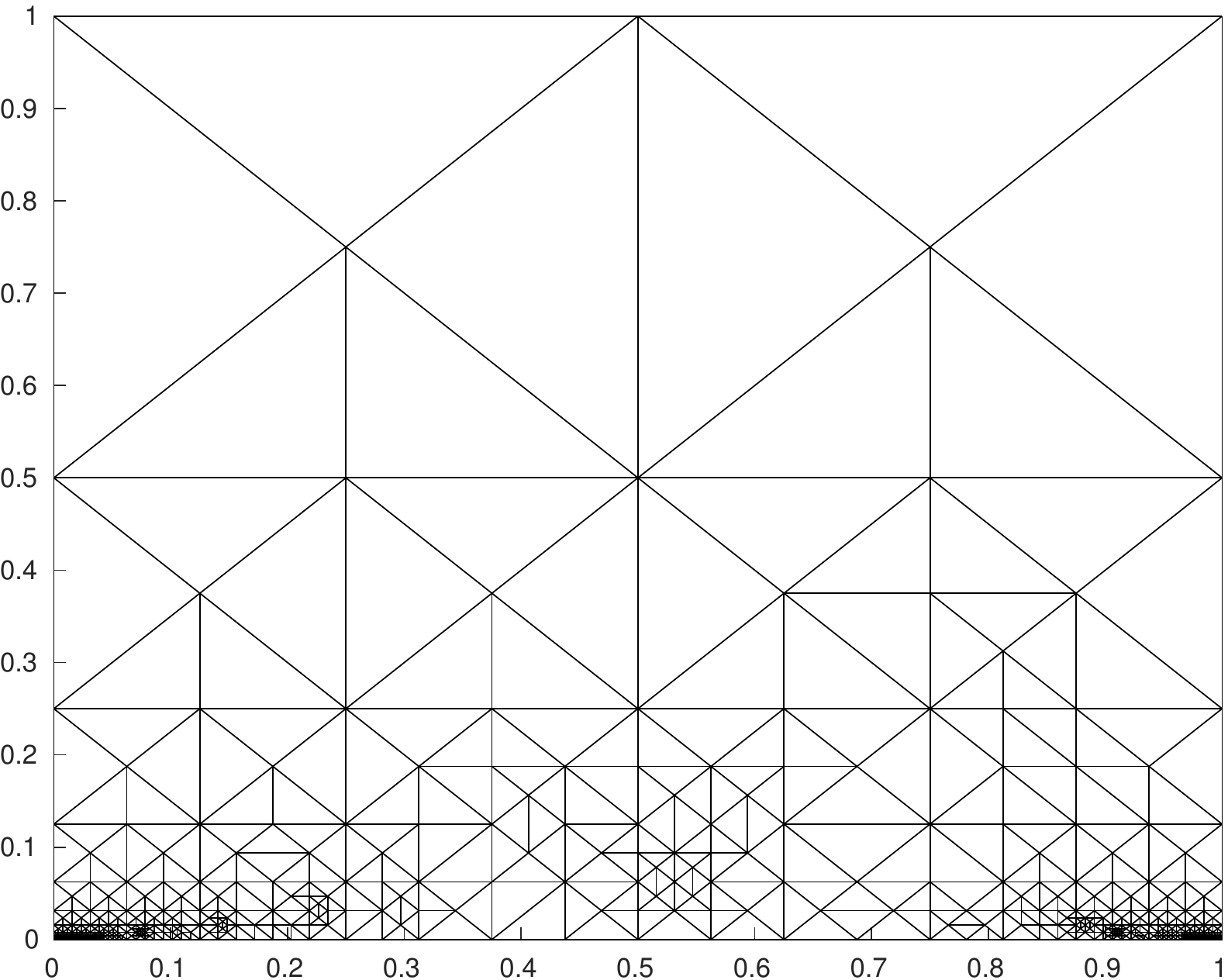}
}
\caption{Meshes generated by the adaptive schemes
driven by the VEM error indicator~$\eta$ in~\eqref{global-error-indicator} (left panel) and the continuous FEM error indicator~$\etaFEM$ in~\eqref{FEM-error-indicator} (right panel)
for the test case with exact solutions~$u_3$ in~\eqref{test-case-3}.}
\label{FIG::SEQUENCE-MESHES-INCOMPATIBLE}
\end{figure}

%%%%
In Figures~\ref{FIG::SEQUENCE-MESHES-t-alpha} and~\ref{FIG::SEQUENCE-MESHES-INCOMPATIBLE},
we plot some meshes
for the test cases with exact solutions~$u_2$ ($\alpha=0.55$) and~$u_3$, respectively, produced by the adaptive procedure driven by the VEM error indicator
in~\eqref{global-error-indicator} and the continuous FEM error indicator in~\eqref{FEM-error-indicator}. 

Next, in Tables~\ref{table:u2-0.55} and~\ref{table:u3},
we focus on the adaptive mesh refinements
driven by the VEM error indicator in~\eqref{global-error-indicator}
and report the number of time-slabs and ``reference elements''
as in Sections~\ref{subsection:flag} and~\ref{subsection:flagging}
for some adaptively generated meshes;.
The time-slab and the element-topology strategies 
allow us to hasten considerably
the assembling and the solving time of method~\eqref{VEM}.

\begin{table}[!ht] 
\centering
\begin{tabular}{ccccc} 
Mesh & Number of reference elements & Total number of elements & Number of time-slabs\\
\hline
$m_2$ & 1 & 4 & 2  \\
$m_4$ & 1 & 52 & 5 \\
$m_6$ & 16 & 265 &  10 \\
$m_8$ & 19 & 1189 & 21 \\
$m_{10}$ & 24 & 5110 & 45 \\
$m_{12}$ & 24 & 21883 & 103
\end{tabular}
\caption{Performance of flagging strategies for the test case~$2$ with~$\alpha = 0.55$, $p=2$ and adaptive refinements.}
\label{table:u2-0.55}
\end{table}
%%%%%%%%
%%%%%%%%
\begin{table}[!ht] 
\centering
\begin{tabular}{ccccc}
Mesh & Number of reference elements & Total number of elements & Number of time-slabs\\
\hline
$m_2$ & 1 & 4 & 2 \\
$m_4$ & 3 & 22 &  4 \\
$m_6$ & 9 & 91 & 5 \\
$m_8$ & 21 & 361 & 9 \\
$m_{10}$ & 33 & 1204 & 13 \\
$m_{12}$ & 35 & 3493 & 21
\end{tabular}
\caption{Performance of flagging strategies for the test case~$3$ with~$p=1$ and adaptive refinements.}
\label{table:u3}
\end{table}

%%%%%%%%%%%%%%%%%%%%%%%%%%%%%%%%%%%%%%%%%%%%%%%%%%%%%%%%%%%%%%%%%%%%%%%%%%
\section{Conclusions} \label{section:conclusions}
%%%%%%%%%%%%%%%%%%%%%%%%%%%%%%%%%%%%%%%%%%%%%%%%%%%%%%%%%%%%%%%%%%%%%%%%%%

We extended the virtual element framework of~\cite{Gomez-Mascotto-Moiola-Perugia:2022}
to the case of general prismatic space--time meshes with possible hanging nodes
and nonuniform degrees of accuracy, useful in~$\h\p$-adaptive procedures.
We discussed flagging strategies to handle the space--time meshes data structure;
this improves the performance of the method.
We investigated numerically the $\h\p$-version of the method
and demonstrated the expected exponential convergence in terms of suitable roots of the number of degrees of freedom for singular solutions.
A residual-type error indicator was introduced
and an~$h$-adaptive refinement procedure was tested based on it.
This error indicator appears to be reliable and efficient for the error~$\EcalY$
but not for the error~$\EcalN$ in~\eqref{exact-errors}.
We investigated numerically the overall virtual element adaptive procedure.
The results obtained for certain singular solutions show that it outperforms a corresponding one for a continuous finite element method.

%%%%%%%%%%%%%%%%%%%%%%%%%%%%%%%%%%%%%%%%%%%%%%%%%%%%%%%%%%%%%%%%%%%%%%%%%%
{\footnotesize
\bibliography{bibliogr.bib}}
\bibliographystyle{plain}

\end{document}